\documentclass[11pt,a4paper,reqno]{amsart}
\usepackage{amsmath,amssymb,mathtools,mathrsfs}
\usepackage[T1]{fontenc}
\usepackage[utf8]{inputenc}
\usepackage{lmodern}
\usepackage[textwidth=162mm,textheight=237mm,centering,headsep=6mm]{geometry}
\usepackage{microtype}
\usepackage{longtable,booktabs,array,calc}
\usepackage{fvextra,needspace,etoolbox}
\usepackage{xurl}
\usepackage[unicode,hidelinks]{hyperref}
\usepackage{bookmark}
\AtBeginEnvironment{longtable}{\footnotesize}
\DefineVerbatimEnvironment{verbatim}{Verbatim}{breaklines=true,breakanywhere=true,fontsize=\footnotesize}

\makeatletter
\patchcmd\longtable{\par}{\if@noskipsec\mbox{}\fi\par}{}{}
\newcommand{\ResultHeading}[1]{\if@noskipsec\leavevmode\par\fi\par\addvspace{5pt}\noindent\textbf{#1}}
\newcommand{\ProofHeading}[1]{\needspace{4\baselineskip}\if@noskipsec\leavevmode\par\fi\par\addvspace{3pt}\noindent\textit{#1}}
\makeatother
\title[Canonical volumes of threefolds]{Canonical volumes of threefolds\\ with positive geometric genus}
\author{Pinxian Bie}
\author{Peien Du}
\author{Zhengjie Yu}
\date{10 November 2026}
\hypersetup{pdftitle={Canonical volumes of threefolds with positive geometric genus},pdfauthor={Pinxian Bie, Peien Du, Zhengjie Yu},pdflang={en}}
\begin{document}
\begin{abstract}
We prove that a complex projective threefold of general type with canonical
singularities and positive geometric genus has canonical volume at least
\(1/30\), with strict inequality when its Euler characteristic is positive.
We describe the equality cases as canonical hypersurfaces of degree \(28\)
in \(\mathbb P(1,3,4,5,14)\). These canonical models and their smooth
projective resolutions are simply connected. Their coarse moduli space is
irreducible and rational of dimension \(103\), with finite quotient
singularities; all admissible \(\mathbb Q\)-Gorenstein deformations are
unobstructed and arise by varying the equation. 
\end{abstract}
\maketitle
\markboth{Canonical volumes of threefolds}{Canonical volumes of threefolds}
{\small\tableofcontents}
\medskip\normalfont
\needspace{4\baselineskip}

\hypertarget{introduction}{%
\needspace{10\baselineskip}\hypertarget{section-1}{}\section{1. Introduction}\label{introduction}}

The canonical volume measures the asymptotic growth of pluricanonical sections.
For threefolds of general type, its smallest values reflect both the geometry
of low-degree pluricanonical maps and the singularities of a minimal model.
We study threefolds with positive geometric genus, and determine the geometry
and moduli of those attaining the smallest volume.

We work over \(\mathbb C\). For a projective threefold \(X\) with canonical
singularities, we write
\[
P_m(X)=h^0(X,\omega_X^{[m]}),\qquad p_g(X)=P_1(X),\qquad
\operatorname{Vol}(X)=\limsup_{m\to\infty}\frac{6P_m(X)}{m^3}.
\]
Here \(\omega_X^{[m]}=(\omega_X^{\otimes m})^{**}\), and \(\varphi_m\)
denotes the map defined by the complete \(m\)-canonical system. On a terminal
minimal model \(Y\), the volume is \(K_Y^3\). Our first result is the following.

\hypertarget{result-1.1}{}\ResultHeading{Theorem 1.1.} Let \(X\) be a complex projective threefold of
general type with canonical singularities and \(p_g(X)>0\). Then
\[
\operatorname{Vol}(X)\ge\frac1{30}.
\]
The bound is sharp. If \(\chi(\mathcal O_X)>0\), the inequality is strict.

We also describe all canonical models for which equality holds.

\hypertarget{result-1.2}{}\ResultHeading{Theorem 1.2.} Let \(X\) be as in \hyperlink{result-1.1}{Theorem 1.1}, with
\(\operatorname{Vol}(X)=1/30\), and let \(Z\) be its canonical model.
Then \(p_g(X)=1\), \(q(X)=h^{2,0}(X)=\chi(\mathcal O_X)=0\), and
\[
R(X,K_X)\simeq
\mathbb C[s,x,y,z,w]/(w^2-F_{28}(s,x,y,z)),
\qquad \deg(s,x,y,z,w)=(1,3,4,5,14).
\]
In particular, \(Z\simeq X_{28}\subset\mathbb P(1,3,4,5,14)\) and
\(\omega_Z\simeq\mathcal O_Z(1)\). Conversely, every normal hypersurface
of this form with canonical singularities has \(p_g=1\) and volume \(1/30\).
The projection forgetting \(w\) is an intrinsic finite double cover.
The fourth-canonical map has connected general fibre of genus two;
\(\varphi_m\) has degree two for \(5\le m\le13\), and is birational
for every \(m\ge14\).

\hypertarget{result-1.3}{}\ResultHeading{Theorem 1.3.} Every canonical model in \hyperlink{result-1.2}{Theorem 1.2} admits an
admissible \(\mathbb Q\)-Gorenstein deformation to a quasismooth hypersurface
of the same degree and weights. Both the canonical model and each of its
smooth projective resolutions are topologically simply connected.

Let \(\mathcal M_{1/30,1}\) be the coarse moduli space of normal canonical
threefolds with ample canonical divisor, \(K^3=1/30\), and \(p_g=1\).
We use Koll\'ar families: the reflexive powers of the relative dualizing
sheaf are flat and commute with base change. These are the admissible
\(\mathbb Q\)-Gorenstein families considered in this paper.

\hypertarget{result-1.4}{}\ResultHeading{Theorem 1.4.} The space \(\mathcal M_{1/30,1}\) is irreducible
and rational of dimension \(103\), and has finite quotient singularities.
For every member \(Z\), the semiuniversal admissible
\(\mathbb Q\)-Gorenstein deformation space is smooth of dimension \(103\).
Every such deformation is induced by varying the weighted equation.
A very general member has automorphism group \(\mathbb Z/2\), generated
by the covering involution. This involution acts trivially on the deformation
space, and the coarse moduli space has a dense smooth open subset.

The last three theorems are counterparts, for \(p_g=1\) and volume \(1/30\),
of the small-volume descriptions of Chen--Hu--Jiang
\cite[Theorems 1.4--1.6]{ref9}. Their hypersurface recognition theorem
\cite[Theorem 8.1]{ref9} supplies the canonical-ring description once we
prove that \(\varphi_5\) is generically finite. We also follow their
specialization argument for simple connectedness
\cite[Section 9.3]{ref9}. The deformation statements above include
nonquasismooth canonical members. We distinguish the smooth deformation
space from its finite quotient, which is the local coarse moduli space.

The study of small canonical volumes in dimension three combines
restriction to surfaces and curves with Reid's Riemann--Roch formula and
the analysis of baskets of singularities. Chen established the sharp
bounds \(1/3\), \(1\), and \(2\) for \(p_g\ge2\), \(p_g\ge3\), and
\(p_g\ge4\), respectively \cite[Theorems 1.4--1.5]{ref8}.
Chen--Hu--Jiang completed the Noether inequality
\(\operatorname{Vol}(X)\ge \frac43p_g(X)-\frac{10}{3}\)
\cite[Theorem 1.1]{ref9}, whose right-hand side gives no positive
bound when \(p_g=1\).

For \(p_g=1\), Chen--Chen proved the sharp bound \(1/30\) under the
assumption \(\chi(\mathcal O_X)\le0\), and determined the numerical data
at equality \cite[Theorem 1.1]{ref7}. Without this assumption, their
general estimate was \(1/75\) \cite[Corollary 1.7]{ref5}.
Jiang's Severi-type inequalities give the stronger bound \(3/8\) when
\(q(X)>0\) \cite[Theorem 1.5]{ref16}. Thus the remaining volume problem
concerns \(p_g=1\), \(q=0\), and positive Euler characteristic. In this
case \(h^{2,0}(X)=\chi(\mathcal O_X)\), so global holomorphic two-forms
are necessarily present.

The birationality problem developed in parallel. For \(p_g=1\),
Chen--Chen obtained birationality in every degree at least eighteen
\cite[Corollary 1.7]{ref5}; Chen--Hu--Penegini proved birationality
in degree seventeen and conjectured the bound fourteen
\cite[Theorem 1.3 and Conjecture B]{ref10}. Yu has established the
sharp fourteenth-canonical theorem for positive geometric genus
\cite[Theorem 1.1]{ref1}, independently of the present volume and
equality theorems. On the marked genus-two curves occurring below,
it gives \(\xi\ge5/14\).

The role of holomorphic forms in volume estimates is also reflected
in the Noether-type inequalities of Chen--Jiang
\cite[Theorem 1.1]{refCJForms}. Our argument uses the global polar
inequality of \cite[Proposition 3.5]{ref1}, together with the
whole-member trace construction of \cite[Theorem 8.2]{refM}.
Normalized polar squares also occur in the proof of the general
\(1/420\) bound \cite[Lemma 9.1]{refM}.
We recall the precise statements needed in
\hyperlink{section-2}{Section 2}; the volume and equality arguments
are then carried out here.

Suppose first that \(V<1/30\). We use canonical restriction and
Reid's Riemann--Roch formula to reduce to a four-canonical pencil with
minimal general fibre of type \((2,3,0)\), \((2,2,0)\), or \((1,2,0)\),
where the entries denote \((K^2,p_g,q)\). The special pencil model
of Chen--Hu--Jiang \cite[Proposition 3.1]{ref9} keeps track of actual
divisor coefficients. This is needed because a rational linear equivalence
alone does not determine the marked poles or the local characters.

The global polar obstruction first bounds \(h^{2,0}\), supplies a reduced visible canonical component, and excludes repeated low-degree derivative branches. We then analyze only the surviving normalized canonical sections. Estimates for graded subrings turn their polynomial degrees into lower bounds for the volume, while small spaces of sections constrain the primitive denominators. The basket lists record
necessary numerical possibilities; their enumeration does not assume
that these possibilities are geometrically realizable.

At equality we revisit the restriction inequalities. The resulting
adjunction and residue calculations exclude positive Euler characteristic.
When \(\chi=0\), the known basket determines the Hilbert series. We prove
that the fifth-canonical map is generically finite by excluding its two
possible quadric surface images. We may then apply the recognition theorem
cited above. Finally, lifting the generators and their relation over
Artinian bases identifies the admissible deformations with variations of
the equation. Counting coefficients and coordinate changes gives dimension
\(103\); the triangular coordinate group and Miyata's invariant-field
theorem give rationality.

\hyperlink{section-2}{Section 2} collects the preliminary tools and proves the multiplicity and
degree estimates. \hyperlink{section-3}{Section 3} reduces the strict volume problem to the three
fibre types. Sections 4--6 treat these types, using the canonical-ring
arguments of \hyperlink{section-7}{Section 7}. \hyperlink{section-8}{Section 8} completes the volume proof and excludes
positive Euler characteristic at equality. \hyperlink{section-9}{Section 9} proves the geometric
and moduli statements. The appendices contain the exact polynomial
calculation and the basket computations, including the closed endpoint.
Executable copies of their Python programs are supplied as ancillary files.

\medskip
\noindent\textbf{AI assistance.}
The authors have long been interested in this problem before the rise of AI.The authors’ principal contributions lie in the birational arguments involving canonical restriction, actual divisor markings and local indices; the framework using holomorphic two-forms and global polar obstructions; the pluriform trace and fixed-part estimates; and the deformation and moduli analysis. During this research, ChatGPT (OpenAI) provided a differential-geometric perspective and suggested several key ideas for the analysis of canonical subrings and function fields, the local differential and Darboux arguments, and the treatment of equality cases, including the route to applying the cited hypersurface recognition theorem. These suggestions informed the development of the proofs and their integration with the geometric framework. The authors take full responsibility for the mathematical content, references and final presentation.

\noindent\textbf{Acknowledgements.}
We are deeply grateful to our supervisor, Meng Chen, whose teaching and insights into explicit birational geometry, especially pluricanonical maps and canonical volumes, have laid the foundation for this work. We thank him in particular for drawing our attention to the canonical volume problem at an early stage of our research.

\needspace{4\baselineskip}

\hypertarget{preliminaries}{%
\needspace{10\baselineskip}\hypertarget{section-2}{}\section{2. Preliminaries}\label{preliminaries}}

\needspace{4\baselineskip}

\hypertarget{models-sections-and-baskets}{%
\hypertarget{section-2.1}{}\subsection{2.1. Models, sections, and baskets}\label{models-sections-and-baskets}}

Throughout the paper we work over \(\mathbb C\), and by a fibration we mean a surjective morphism with connected fibres. By the minimal model programme and finite generation \cite{ref2,ref18}, we may choose a terminal \(\mathbb Q\)-factorial minimal model \(Y\) and a canonical model \(Z\) of \(X\). On a smooth common model \(\pi:W\to Y\), choose a nonzero canonical section \(s\) and put
\[
Q=\pi^*\operatorname{div}_Y(s),\qquad V=Q^3,\qquad
h=h^0(W,\Omega_W^2),\qquad q=h^1(W,\mathcal O_W).
\]
With this notation, \(Q\) is nef and big, \(K_W=Q+E_\pi\) with \(E_\pi\ge0\), and \(\chi(\mathcal O_W)=1-q+h-p_g\). We use
\[
R=\bigoplus_{m\ge0}R_m,\quad R_m=H^0(Z,\omega_Z^{[m]})=H^0(W,mK_W),
\quad E_m=s^{-m}R_m,\quad L_m=\mathbb C(E_m).
\]
Multiplication by \(s\) gives \(E_m\subset E_{m+1}\); \(\varphi_m\) is birational precisely when \(L_m=\mathbb C(W)\). Throughout, \(D\ge_{\mathbb Q}D'\) means that \(D-D'\) is \(\mathbb Q\)-linearly equivalent to an effective divisor. An inequality without this subscript compares actual coefficients. Rounding is always applied to an actual divisor.

For a terminal basket we use \((b,r)=\frac1r(1,-1,b)\), with \(0<b\le r/2\) and \((b,r)=1\), retaining multiplicities. Reid's formula gives, for \(m\ge2\),
\[
P_m=\frac{m(m-1)(2m-1)}{12}V-(2m-1)\chi
 +\sum_{(b,r)\in B}\sum_{j=1}^{m-1}
 \frac{\overline{jb}(r-\overline{jb})}{2r},                 \tag{2.1}\label{eq:2.1}
\]
where \(\overline{jb}\) is the least nonnegative residue \cite{ref21,ref6}. In particular,
\[
V=2P_2+6\chi-\sum_{(b,r)\in B}\frac{b(r-b)}r.             \tag{2.2}\label{eq:2.2}
\]
A prime packing replaces \((b,r),(c,u)\), with \(|bu-cr|=1\), by \((b+c,r+u)\) and decreases \(V\) by \(1/(ru(r+u))\). A basket constituent is not necessarily an actual singular point. A local-index assertion always concerns an actual point of \(Y\).

\needspace{4\baselineskip}

\hypertarget{restriction-and-marked-curves}{%
\hypertarget{section-2.2}{}\subsection{2.2. Restriction and marked curves}\label{restriction-and-marked-curves}}

We now fix a general moving surface \(F\). Let \(\sigma:F\to S\) be its smooth minimal-model contraction and put \(D=Q|_F\), \(A=\sigma^*K_S\). We recall the restriction results in \cite[Propositions 2.1--2.3]{ref1}. If \(Q\ge_{\mathbb Q}\mu F\) and \(F\) is semiample, then
\[
D\ge_{\mathbb Q}\frac{\mu}{1+\mu}A,\qquad V\ge\mu D^2.
\tag{2.3}\label{eq:2.3}
\]
If \(Q-\mu F\) is only pseudoeffective, use any rational \(0<p<\mu\) in the first inequality and then let \(p\to\mu\). If \(C\) is a general member of a free curve system on \(F\), \(D\ge_{\mathbb Q}\beta C\) and \(\xi=DC\), put
\[\alpha(m)=(m-1-\mu^{-1}-\beta^{-1})\xi.\]
For \(\alpha(m)>1\),
\[m\xi\ge2g(C)-2+\lceil\alpha(m)\rceil.\tag{2.4}\label{eq:2.4}\]
With separation of the surface and curve systems, \(\alpha(m)>2\) gives birationality. We also use the rounding consequence \cite[Lemma 2.4]{ref1}: if \(c>n\ge1\) are integers, \(x>0\) is rational, \(\gcd(c,n+1)>1\), and \(mx\ge n+\lceil(m-c)x\rceil\) for every sufficiently large \(m\), then \(x\ge(n+1)/c\).

To compare actual coefficients, we use the special pencil model of \cite[Proposition 3.1]{ref9}, in the form proved in \cite[Proposition 2.2]{ref1}. For the complete connected four-canonical pencil this gives
\[
0\le D\le K_F\le5D,\qquad
5D\ge\sigma^*A_0,\quad A_0=\sigma_*K_F\in|K_S|.
\tag{2.5}\label{eq:2.5}
\]
If we subsequently blow up this model, we pull back both sides of the comparison. In particular, the symbol \(K_F\) in \eqref{eq:2.5} continues to denote the pullback of the canonical divisor on the special model. For \(m\ge6\), vanishing gives
\[|mK_W||_F\succeq|K_F+\lceil(m-5)D\rceil|\succeq|K_F|.\]
The system separates the outer fibres and the canonical curves on \(F\). Changes of surface or curve pencil use \cite[Propositions 3.6--3.7]{ref10}, with the precise separation conventions of \cite[Section 4.4]{ref1}.

For a \((1,2)\) minimal fibre, a general resolved canonical curve \(C\) has a marked Weierstrass point \(p\) with
\[
Q|_C=\xi p,\qquad
(a/s^m)|_C\in H^0(C,\mathcal O_C(\lfloor m\xi\rfloor p))
\quad(a\in R_m).
\tag{2.6}\label{eq:2.6}
\]
These are actual divisor statements \cite[Lemmas 2.5--2.6]{ref1}. The semigroup at \(p\) is \(\langle2,5\rangle\), so \cite[Theorem 1.1]{ref1} gives \(\xi\ge5/14\). The actual local index \(r_P\) at \(P=\pi(p)\) clears \(\xi\) and occurs among the basket denominators. For a smooth cyclic index-one cover and reduced \(\xi=a/e\), the canonical coefficient has multiplicity at most \(a\). We use its character as well as its multiplicity. If \(\xi<2/5\) in the quartic pencil, the minimal-model map is an isomorphism over the canonical base point and \(AD=\xi\) \cite[Lemma 5.3]{ref1}. On a minimal \((1,2)\) surface, an integral curve through a very general point with canonical degree one belongs to the canonical pencil \cite[Lemma 5.4]{ref1}. The impossible square-two endpoint and the canonical fixed-cycle calculation are used with their precise hypotheses in \cite[Proposition 4.3 and Proposition 4.5--Lemma 4.6]{ref1}.

\needspace{4\baselineskip}

\hypertarget{two-forms-and-the-canonical-quotient}{%
\hypertarget{section-2.3}{}\subsection{2.3. Two-forms and the canonical quotient}\label{two-forms-and-the-canonical-quotient}}

We next recall the differential tools from \cite[Section 3]{ref1}. Throughout this subsection we assume \(q(W)=0\). A mixed Hodge structure has a weight filtration; cohomology of a smooth projective variety is pure of its cohomological degree, and morphisms are strict for this filtration \cite[Theorem 2.3.5 and Corollaries 3.2.13(ii), 3.2.14, and 3.2.17]{ref12}. If \(q(W)=0\), these facts imply that a holomorphic form expressible as a sum of wedges of closed rational one-forms with holomorphic two-forms is zero. The analogous statement holds for a holomorphic wedge of two closed rational one-forms \cite[Lemma 3.1]{ref1}.

Let \(0\ne\alpha\in H^0(W,\Omega_W^2)\), and define the meromorphic vector field \(v\) by \(\iota_vs=\alpha\). We then have \(df\wedge\alpha=v(f)s\) and \(L_vs=0\). If \(v(f)\ne0\), the functions \(f,v(f)\) are algebraically independent. A finite-dimensional subspace of \(\mathbb C(W)\) stable under \(v\) is annihilated by \(v\) \cite[Lemma 3.2]{ref1}. We refer to this as the Hodge stability lemma. For a regular general fibre over \(\mathbb P^1\), restriction of absolute two-forms is injective and \(h\le p_g(F)-1\). For \(p_g(F)=1\) all such forms restrict to zero and their evaluation rank is at most \(q(F)\).

If a fibre has \(q(F)=1\) and a nonzero absolute two-form restricts to zero, then \(AD\ge1\) \cite[Lemma 3.16]{ref1}.

Let \(\mathcal A=R/(s)\), and let \(\mathcal S\) be its graded saturation, the ring of regular functions on \(\operatorname{Spec}\mathcal A\) minus its vertex. Nilpotents are retained. For \(q=0\), \cite[Proposition 3.3]{ref1} gives an injection \(\mathcal A\hookrightarrow\mathcal S\), with quotient \(\mathbb C^h\) in degree two and zero in every other degree. There are a degree-two derivation \(\Delta:R\to R\) and \(\theta=\theta_\alpha\in\mathcal S_2\) such that
\[
\Delta_m(b)=s^{m+2}v(b/s^m),\qquad
\overline{\Delta_m(b)}=-m\theta\bar b.\tag{2.7}\label{eq:2.7}
\]
If \(P_1=P_2=P_3=1\) and \(P_4=2\), write \(t=x/s^4\) for \(x\) independent of \(s^4\). Then \(\theta^2=\lambda\bar x\) and
\[v^2(b/s^m)-m(m+2)\lambda t(b/s^m)\in E_{m+3}.\tag{2.8}\label{eq:2.8}\]
In a local canonical frame \(\nu\), with \(s=g\nu\) and \(b=f\nu^m\), the formula behind these assertions is
\[\Delta_m(b)=\nu^{m+1}\otimes((g\,df-mf\,dg)\wedge\alpha).\tag{2.9}\label{eq:2.9}\]
The local expression \(g\,df-mf\,dg\) is the negative of the weighted Gaussian operator \(\gamma_{1,m}(g\otimes f)\); see \cite[formula (1)]{refBF}. Wedge product with \(\alpha\), with the indicated powers of the canonical trivialization, gives the operator used here.

It also gives \(qb\Delta_p(a)-pa\Delta_q(b)\in sR_{p+q+1}\). These formulas will be called the Euler identities. They apply in the full generic Artin rings of the canonical divisor; cancellation modulo \(s\) is checked there, not merely on its reduced support.

The following global polar inequality is \cite[Proposition 3.5]{ref1}.
The normalized polar-square construction also appears in the proof of the
general volume bound \cite[Lemma 9.1]{refM}; we use the formulation of
\cite{ref1}, which applies directly to a nonzero canonical section.

\hypertarget{prop-volume-polar}{}\ResultHeading{Proposition 2.1 (Global polar inequality).}
Let \(h=h^{2,0}(W)\). For every \(m\ge1\) and every nonzero \(a\in R_m\), there is a homogeneous quadratic map
\[
 \mathfrak P_a:H^0(W,\Omega_W^2)\longrightarrow R_{2m+2}/a^2R_2
\]
whose inverse image of zero is \(\{0\}\). In particular,
\[
 h\le P_{2m+2}-P_2.
\]

If \(P_1>0\), this gives \(h\le P_4-P_2\).

\ProofHeading{Proof.}
The construction and its nonvanishing are proved in \cite[Proposition 3.5]{ref1}. The common zero locus of fewer than \(h\) homogeneous equations of positive degree in \(h\) variables is positive-dimensional. Applying this to the coordinates of \(\mathfrak P_a\) gives the inequality. \(\square\)

\hypertarget{cor-volume-canonical-polar}{}\ResultHeading{Corollary 2.2 (Canonical polar consequences).}
Suppose \(q=0\), \(P_1=P_3=1\), and \(h>0\). For every nonzero absolute two-form \(\alpha\), the associated Euler residue \(\theta_\alpha\) is a unit at the generic point of a reduced prime component \(T\) of the canonical divisor on the terminal model. This component survives on the canonical model. Consequently
\[
 P_{2r}>P_{2r-1}\qquad(r\ge2).
\]
If also \(P_4=2\), choose \(x\in R_4\setminus sR_3\), and put \(t=x/s^4\). Then
\[
 \theta_\alpha^2=\lambda\bar x,\quad \lambda\ne0,\qquad
 x|_T\ne0,\qquad t'\notin E_5.
\]
In particular no \(W\in R_7\) satisfies \(sW=x^2\).

\ProofHeading{Proof.}
The reduced visible component is supplied by \cite[Corollary 3.6]{ref1}. Every power of its Euler residue remains nonzero there. Since saturation changes \(\mathcal A\) only in degree two, \(\theta_\alpha^r\ne0\) lies in \(\mathcal S_{2r}=\mathcal A_{2r}\) for \(r\ge2\), proving the strict even increments. If \(P_4=2\), then \(\mathcal A_4=\mathbb C\bar x\), so the nonzero square of \(\theta_\alpha\) is \(\lambda\bar x\), with \(\lambda\ne0\); both \(\theta_\alpha\) and \(x\) are units at \(T\). If \(t'\in E_5\), then \(\Delta_4(x)\in sR_5\), and the Euler identity gives \(\theta_\alpha\bar x=0\), contradicting its value at \(T\). The relation \(sW=x^2\) is excluded at the same component. This argument uses the reduced visible component, rather than an inference about powers of an arbitrary element of a nonreduced ring. \(\square\)

\needspace{4\baselineskip}

\hypertarget{local-expansions-and-canonical-subrings}{%
\hypertarget{section-2.4}{}\subsection{2.4. Local expansions and canonical subrings}\label{local-expansions-and-canonical-subrings}}

We call a component \(T\) of the canonical divisor defined by \(s\) transverse if \(\alpha|_T\ne0\), and tangent otherwise. At a general point of a transverse component of multiplicity \(\ell\), \cite[Lemma 3.7]{ref1} supplies analytic coordinates with
\[\alpha=dy_1\wedge dy_2,\qquad
s=d\tau\wedge dy_1\wedge dy_2,\qquad
\tau=\rho^{\ell+1},\qquad v=\partial_\tau.\]
We use analytic straightening for a nonvanishing holomorphic vector field, in the form of analytic dependence of its local flow \cite[Theorem 1]{refCarrillo}. The coefficient and resonance comparisons below belong to the classical Kowalevski--Painlev\'e method \cite[\S 1]{refKowalevski}, \cite[\S 3]{refConteMusette}. Here all expansions come from meromorphic functions on the given variety, so their convergence is already available; we do not assume existence or convergence of arbitrary formal solutions, nor a global Painlev\'e property.

We call \(\tau\) the actual time. A meromorphic function has an existing convergent Puiseux expansion in \(\tau\), with exponents in \((\ell+1)^{-1}\mathbb Z\) and analytic transverse coefficients. At a tangent component, \(\operatorname{ord}_T(vf)\ge\operatorname{ord}_T(f)-\ell\). At least one transverse component exists on a smooth projective model.

For algebraically independent \(f_1,f_2,f_3\), if \(K=\mathbb C(f_1,f_2,f_3)\) has prime index in \(L=\mathbb C(W)\), then
\[J=(df_1\wedge df_2\wedge df_3)/s\notin K,\qquad L=K(J).\]
The three functions cannot all have integral Laurent expansions in actual time at a transverse component \cite[Lemma 3.8]{ref1}. If \(v(K)\subset K\), then \(v(J)/J\) is the coordinate divergence of \(v\) and has at most a simple time pole. A first nonintegral exponent of an existing expansion must solve the corresponding homogeneous indicial equation: earlier integral exponents cannot force it. We use the resulting local tests exactly as in \cite[Lemma 3.10]{ref1}. An indicial equation here means the coefficient equation for that first unknown Laurent coefficient, not an assertion of existence of formal solutions. In the quadratic case the involution negates the absolute two- and three-forms and preserves \(v\) \cite[Corollary 3.9]{ref1}.

For a finitely generated graded algebra \(C\) of dimension \(d\), write
\[e(C)=\lim_{q\to1^-}(1-q)^d\sum_{m\ge0}(\dim C_m)q^m.\]
A polynomial ring of weights \(w_i\) has multiplicity \(1/\prod w_i\); a hypersurface of degree \(a\) has multiplicity \(a/\prod w_i\). Also \(e(R)=V\). For a graded subring \(B\subset R\) of dimension four and field index \(\delta\),
\[V\ge\delta e(B),\qquad \dim(UV)\ge\dim U+\dim V-1\tag{2.10}\label{eq:2.10}\]
for nonzero finite-dimensional homogeneous subspaces of a domain \cite[Lemma 3.11]{ref1}. Here we do not need \(R\) to be finite over \(B\). If \(rK_Y\) is Cartier, \(rK_Y^2T\) is integral for any integral Weil divisor \(T\). When \(p_g=1\) and \(rV=1\), the canonical divisor on \(Z\) is reduced and prime, \(R/(s)\) is a domain, and the product inequality applies also to its saturation \cite[Lemmas 3.12--3.13]{ref1}.

When the normalized canonical field contains a rational subfield of index two, let \(\iota\) be the corresponding involution, and put \(j|_{R_m}=(-1)^m\iota^*\). A section is invariant if its normalized function lies in the subfield, and odd if \(j\) acts by \(-1\). We use the canonical-form criterion \cite[Lemma 3.14]{ref1}: let \(C\subset R^j\) be a finitely generated graded domain with \(\operatorname{Frac}C=\operatorname{Frac}R^j\). A nonzero homogeneous rational top form of degree \(d\) on \(\operatorname{Spec}C\), extending to a resolution of a proper birational model, gives an odd section of degree \(d+1\). In particular, if \(B\subset R^j\) is polynomial with that same fraction field, a primitive invariant fraction over \(B\) having a polynomial coordinate of weight \(a\) in its denominator gives an odd section of degree \(S+1-a\), where \(S\) is the sum of generator weights. This criterion uses reflexive canonical forms and does not require a Cohen--Macaulay canonical ring.

\needspace{4\baselineskip}

\hypertarget{polar-divisors-and-degree-growth}{%
\hypertarget{section-2.5}{}\subsection{2.5. Polar divisors and degree growth}\label{polar-divisors-and-degree-growth}}

\needspace{3\baselineskip}

\hypertarget{result-2.1}{}\ResultHeading{Proposition 2.3 (multiplicity and polar divisors).} Let \(B\) be a normal
finitely generated positively graded domain of dimension \(n\ge2\),
with \(B_0=\mathbb C\). Let \(U\in\operatorname{Frac}B\setminus B\)
be homogeneous of degree \(d>0\). If its polar height-one primes are
\(P_i\), with normalized pole orders \(m_i\), then
\[
e(B[U])\ge e(B)+\frac1d\sum_i m_i e(B/P_i).              \tag{2.11}\label{eq:2.11}
\]

\ProofHeading{Proof.} Normality identifies the denominator ideal
\[
I=(B:_B U)=\bigcap_iP_i^{(m_i)}.
\]
We choose a graded complex vector-space complement \(C\) to \(I\) in
\(B\). There is an injection
\[
B\oplus\bigoplus_{k\ge1}U^kC\longrightarrow B[U].
\]
To see that this sum is direct, suppose that there is a nontrivial relation. We take the largest exponent
\(k\) and its nonzero coefficient \(c_k\). As \(c_k\notin I\), some
\(i\) has \(v_{P_i}(c_k)<m_i\). The value of \(U^kc_k\) is strictly
below \(-(k-1)m_i\), whereas every lower-exponent term has value at
least that number. Thus the term of largest exponent cannot cancel, which proves the asserted injection. Consequently
\[
H_{B[U]}(q)\ge H_B(q)+\frac{q^d}{1-q^d}H_{B/I}(q)
\]
coefficientwise. The localization \(B_{P_i}\) is a discrete valuation
ring and \(I_{P_i}\) has colength \(m_i\). The associativity formula
for multiplicities gives
\[
e(B/I)=\sum_i m_i e(B/P_i).
\]
We obtain the claimed inequality by multiplying by \((1-q)^n\) and letting \(q\to1^-\).
The complement \(C\) is a vector space, not a residue subalgebra;
no compatibility of residue lifts at distinct poles is required.
\(\square\)

\needspace{3\baselineskip}

\hypertarget{result-2.2}{}\ResultHeading{Proposition 2.4 (polynomial degree and volume).} Let
\(s,x,y,z_c\) be algebraically independent canonical sections of degrees
\(1,a,b,c\), where \(a,b,c>0\), and set
\[
t=x/s^a,\qquad u=y/s^b,\qquad z=z_c/s^c.
\]
Suppose that
\(\delta=[\mathbb C(W):\mathbb C(t,u,z)]<\infty\).
If \(H\in R_m\), with \(m>0\), has normalized function
\(h=H/s^m\in\mathbb C[t,u,z]\), write
\(d_t=\deg_t h\), \(d_u=\deg_u h\), and \(d_z=\deg_z h\).
Each positive degree gives the corresponding inequality
\[
V\ge\frac{\delta d_t}{bcm},\qquad
V\ge\frac{\delta d_u}{acm},\qquad
V\ge\frac{\delta d_z}{abm}.                              \tag{2.12}\label{eq:2.12}
\]

\ProofHeading{Proof.} Suppose first that \(d_t>0\). We consider the
polynomial subring \(\mathbb C[s,y,z_c,H]\), whose multiplicity is
\(1/(bcm)\). Since \(h\) is a nonconstant polynomial in \(t\),
\[
[\mathbb C(t,u,z):\mathbb C(h,u,z)]=d_t.
\]
Thus the field index of this subring in the canonical ring is
\(\delta d_t\), and \eqref{eq:2.10} proves the first inequality. We obtain the
other two by interchanging the roles of the variables. \(\square\)

In the applications below, we use \((a,b,c)=(4,5,6)\) and
\(\delta=2\). Thus \eqref{eq:2.12} has denominators \(30m\), \(24m\), and
\(20m\), respectively.

\needspace{4\baselineskip}
\hypertarget{section-2.6}{}\subsection{2.6. Pluriform trace}

We use the whole-member pluriform trace of \cite[Theorem 8.2]{refM}.
Its reduced-member refinement in \cite[Theorem 8.5]{refM} is stated for
\(P_1=0\); we record the version needed here, retaining the adjunction
kernel \(aR_1\).

\hypertarget{lem-pluriform-trace}{}\ResultHeading{Lemma 2.5 (Pluriform trace).}
Let \(0\ne a\in R_m\), where \(m\ge1\), and let \(F\) be a smooth
model of a coefficient-one prime component of \((a)\) on \(Y\). Put
\(c=m+1\) and \(U_{\ell,m}=P_{\ell c}-P_{\ell c-m}\). Then
\[
 U_{\ell,m}\ge P_\ell(F)\qquad(\ell\ge2).                 \tag{2.13}\label{eq:whole-member-trace}
\]
Suppose further that \(a\) defines a general member of a pencil,
\(F\) resolves a connected moving component, and the full member
\(H=(a)\) on \(Z\) is reduced. If
\(p=p_g(F)\) and
\(\rho=\operatorname{rank}(H^0(W,\Omega_W^2)\to H^0(F,K_F))\), then
\[
 U_{\ell,m}\ge P_\ell(F)+\max\{0,P_c-P_1-p+\rho\}.
                                                               \tag{2.14}\label{eq:reduced-member-trace}
\]
For \(q(F)=0\), one has \(\rho=h\).

\ProofHeading{Proof.}
The first assertion is \cite[Theorem 8.2]{refM}. For the second, full-member
adjunction and vanishing \cite[Lemma 8.1]{refM} give
\[
0\longrightarrow R_c/aR_1\longrightarrow H^0(H,\omega_H)
 \xrightarrow{\partial}H^1(Z,\omega_Z)\longrightarrow0.
\]
Ordinary trace from \(F\), followed by \(\partial\), has rank \(\rho\)
by duality; hence its liftable subspace has dimension \(p-\rho\).
Write \(Q_\ell=R_{\ell c}/aR_{\ell c-m}\), let \(T_\ell\) be the
pluriform trace image, and put \(\epsilon=\dim Q_\ell-P_\ell(F)\).
The homogeneous map
\[
R_c/aR_1\longrightarrow Q_\ell/T_\ell,\qquad [z]\longmapsto[z^\ell]
\]
is well-defined and its zero locus has dimension at least
\(P_c-P_1-\epsilon\). Compatibility with powers
\cite[Lemma 8.3]{refM} makes the residue of such a \(z\) regular on
\(F\). Reducedness makes its adjunction section zero generically on
the other components, so Cohen--Macaulayness identifies it everywhere
with ordinary trace from \(F\). Thus the zero locus has dimension at
most \(p-\rho\). Since \(\epsilon\ge0\), this proves
\eqref{eq:reduced-member-trace}. The last assertion is
\cite[Lemma 8.4]{refM}. \(\square\)

\hypertarget{reduction-to-three-fibre-types}{%
\needspace{10\baselineskip}\hypertarget{section-3}{}\section{3. Reduction to three fibre types}\label{reduction-to-three-fibre-types}}

We use the notation and results recalled in \hyperlink{section-2}{Section 2}, including fourteenth-canonical birationality \cite[Theorem 1.1]{ref1}. Suppose throughout that
\[
p_g(X)>0,\qquad V:=K_Y^3<\frac1{30}.
\]
Here \(Y\) is a terminal \(\mathbb Q\)-factorial minimal model of \(X\), \(\pi:W\to Y\) is a suitable resolution, and \(Q=\pi^*K_Y\). When coefficients are considered, \(K_Y\) is the effective divisor defined by a fixed nonzero canonical section \(s\).

\needspace{4\baselineskip}

\hypertarget{the-first-moving-system}{%
\hypertarget{section-3.1}{}\subsection{3.1. The first moving system}\label{the-first-moving-system}}

\hypertarget{prop-first-moving-system}{}\ResultHeading{Proposition 3.1 (The first moving system).} Under the assumptions of this section, \(q=0\), \(p_g=P_2=P_3=1\), \(1\le h=\chi\le2\), and \(P_4\in\{2,3\}\).

\ProofHeading{Proof.} We begin with the following bounds from the literature:
\[
\begin{array}{c|c|c}
\text{condition}&\text{volume bound}&\text{reference}\\\hline
q(X)>0&V\ge3/8&\text{\cite[Theorem 1.5]{ref16}}\\
\chi(\mathcal O_X)\le0&V\ge1/30&\text{\cite[Theorem 1.1]{ref7}}\\
P_2\ge2&V\ge1/14&\text{\cite[Proposition 4.2]{ref5}}\\
P_3\ge3&V\ge1/24&\text{\cite[Proposition 4.3]{ref5}}\\
P_4\ge4&V\ge7/160&\text{\cite[Corollary 3.3]{ref5}}.
\end{array}
\]
Since multiplication by a nonzero canonical section is injective, these bounds imply
\[
q(X)=0,\quad p_g=P_2=1,\quad \chi(\mathcal O_X)>0,\quad
P_3\le2,\quad P_4\le3.
\]
We first exclude \(P_3=2\).

The complete tricanonical system is then a pencil with connected base \(\mathbb P^1\) and moving degree one. Let \(F\) be a general fibre, \(\sigma:F\to S\) its minimal-model contraction, and
\[
D=Q|_F,\qquad A=\sigma^*K_S.
\]
Canonical restriction gives \(D\ge_{\mathbb Q}A/4\) and
\[
V\ge\frac{D^2}{3}\ge\frac{A^2}{48}.
\]
It follows that \(A^2=1\). By Debarre's inequality \cite{ref11}, the surface \(S\) is regular. Restriction of absolute two-forms, recalled in \hyperlink{section-2.3}{Section 2.3}, gives \(p_g(S)\ge2\); Noether's inequality gives the reverse inequality. Consequently \(S\) is a \((1,2)\)-surface and \(\chi=h^{2,0}=1\).

Let \(C\) be a general resolved canonical curve of \(S\). \cite[Lemma 2.5]{ref1} and
fourteenth-canonical birationality give
\[
\frac5{14}\le\xi=DC\le12V<\frac25.
\]

We next consider the restriction of the fifth-canonical system. Its moving divisor \(N\) satisfies \(AN\le60V<2\), so the restriction is composed with the canonical pencil and has rank at most two. If \(P_5\ge4\), the kernel of restriction contains another moving surface. If its pencil is the outer pencil, then \(5Q\ge_{\mathbb Q}2F\), and restriction yields \(V\ge2/49\). If it is a different pencil with the same inner curves, the curve restriction estimate gives
\[
(n+6)\xi\ge2+\lceil n\xi\rceil.
\]
At \(n=3\) this contradicts \(\xi<2/5\). Hence \(P_5\le3\).

The finite basket calculation described below leaves exactly
\[
\begin{array}{c|c|c}
P_4&B&V\\\hline
2&4(1,2)+4(1,3)+(2,7)+(3,8)&5/168\\
2&4(1,2)+3(1,3)+(2,5)+(4,13)&2/65\\
2&4(1,2)+(1,4)+(8,23)&3/92\\
3&5(1,2)+(1,3)+2(1,4)+(2,7)+(3,8)&5/168\\
3&5(1,2)+2(1,4)+(2,5)+(4,13)&2/65.
\end{array}
\]
We now impose the local-index condition at the image of the marked point \(p\). This index clears \(\xi\), but four rows contain no such fraction in \([5/14,12V]\). In the third row one must have \(\xi=9/23\) at the unique index-\(23\) point. This is an actual quotient point: a nonquotient terminal point of index exceeding four is of type \(cA/r\), and its axial weight counts identical constituents in its basket; here that weight is one. The quotient is \(\frac1{23}(1,-1,8)\). The normalized equivariant marked disc has coordinate orders congruent to \(19,4,14\) modulo \(23\). A canonical coefficient of order nine is impossible, since
\[
9\notin\langle19,4,14\rangle.
\]
We have therefore excluded \(P_3=2\), and hence \(P_3=1\).

Proposition~2.1 gives \(1\le\chi=h^{2,0}\le P_4-P_2\le2\), and excludes \(P_4=1\). Thus \(P_4=2\) or \(3\).

\needspace{4\baselineskip}

\(\square\)

\hypertarget{a-quartic-net-cannot-occur}{%
\hypertarget{section-3.2}{}\subsection{3.2. A quartic net cannot occur}\label{a-quartic-net-cannot-occur}}

\hypertarget{prop-quartic-net}{}\ResultHeading{Proposition 3.2 (The quartic pencil).} Under the assumptions of this section, \(P_4=2\) and \(h=\chi=1\).

\ProofHeading{Proof.} If \(P_4=3\) and the quartic image is a curve, the moving degree is at least two, whence \(V\ge1/18\). Suppose therefore that the quartic image is a surface. Write \(M=M_4\) and let \(S\in|M|\) be general. Put \(N=M|_S\), \(D=Q|_S\) and \(A=\sigma^*K_{S_0}\), where \(S_0\) is the minimal model. Adjunction gives \(K_S\ge_{\mathbb Q}5N/4\). Hence
\[
A^2\ge\frac54AN,\qquad A^2\le100V<\frac{10}{3}.
\]
In particular \(A^2\in\{2,3\}\). Write \(N\equiv eC\) for a connected general fibre of the restricted pencil. Since \(S_0\) is not a \((1,2)\)-surface, the covering-curve estimate gives \(AC\ge2\). Thus
\[
3\ge A^2\ge\frac54eAC.
\]
It follows that \(A^2=3\), \(e=1\) and \(AC=2\). Surface Hodge index and adjunction show that \(C_0^2=0\) and \(g(C)=2\). The Stein base of the quartic map is therefore \(\mathbb P^2\).

Let \(f:W\to\mathbb P^2\) be this connected genus-two fibration. By relative duality, the trace map embeds \(R^1f_*\omega_W\) into \(\omega_{\mathbb P^2}\), and is an isomorphism at the generic point. Using torsion-freeness and Koll'ar vanishing \cite{ref19}, we obtain
\[
H^1(W,\omega_W(M))=0.
\]
Indeed \(H^1(\mathbb P^2,f_*\omega_W(1))=0\), while the trace injection gives \(H^0(\mathbb P^2,R^1f_*\omega_W(1))=0\). Adjunction then yields
\[
p_g(S)=h^0(W,\omega_W(M))-1+h^{2,0}(W)\ge2+h^{2,0}(W).
\]
Noether's inequality gives \(p_g(S)\le3\). Hence \(\chi=1\) and \(S_0\) has invariants \((K^2,p_g,q)=(3,3,0)\).

The actual divisor induced by \(s|_S\) descends to an effective divisor
\[
B\sim K_{S_0}-C_0,\qquad K_{S_0}B=1,\qquad
0\le D_0:=\sigma_*D\le B.
\]
Let \(F_0\) be the fibre over the zero of the canonical quartic base section. We claim that
\[
2B-F_0\ge0.
\]
To prove the claim also for nonreduced fibres, we use the genus-two fibre algebras described in \cite[Lemma 4.1 and its proof]{refCP06}. In a two-connected fibre its degree-one subalgebra is a polynomial algebra. Otherwise the fibre algebra is a complete intersection with generators of degrees \(1,1,2,3\), a quadratic relation \(x_0^2-\lambda x_0x_1\), and a relation monic in the degree-three generator. Thus every nilpotent linear element has square zero. Over a DVR this says \(u^4\in t\mathscr R_4\Rightarrow u^2\in t\mathscr R_2\). Applying this implication to the canonical section of the relatively minimal genus-two fibration gives \(2B-F_0\ge0\), as required.

There is exactly one component \(\Gamma\) of \(B\) with positive canonical degree. It has coefficient one, \(K_{S_0}\Gamma=1\), and coefficient two in \(F_0\). Let \(T\subset W\) be the divisor containing its strict transform. Then
\[
\operatorname{ord}_T(s)=1,\qquad a:=\operatorname{coeff}_TQ=AD,\qquad
\frac35<a\le20V<\frac23.
\]
We explain why the first inequality is strict. Equality would give \(D\equiv A/5\): an effective representative of \(D-A/5\) is orthogonal to the nef big divisor \(A\), and its negative square would contradict nefness of \(D\). But \(T\) is exceptional over a point, so \(D\widetilde\Gamma=0\), whereas \(A\widetilde\Gamma=1\). The discrepancy \(1-a<1\) also proves that the centre is a non-Gorenstein point.

The basket budgets are
\[
\sum b=14,\qquad\sum\gamma=4,\qquad
\sum\epsilon=\frac13-V,\qquad\frac3{100}\le V<\frac1{30}.
\]
The gap bounds of the next subsection give a local denominator at most \(33\). We can therefore enumerate the baskets exactly: there are twenty before the remaining tests; monotonicity leaves nineteen, and the displayed strict coefficient interval leaves nine. For a marked quotient \(\frac1r(1,-1,b)\) with \(a=k/r\), the discrepancy-less-than-one valuation has weights
\[
(\langle-kb^{-1}\rangle_r,\langle kb^{-1}\rangle_r,r-k).
\]
The fourth-canonical fixed order along \(T\) is two, and its image is a line. There must therefore be at least two initial monomials of weight \(4k-2r\). The finite monomial test leaves only
\[
B_*=2(1,2)+(1,7)+(3,8)+(8,17),\qquad V=31/952,\qquad a=5/8.
\]
The half-fibre inequality gives \(2/5\le\xi:=QC\le5/12\). The actual divisor \(Q|_C\) is hyperelliptically invariant and bounded by a canonical divisor of degree two. A conjugate pair would have equal coefficients \(\xi/2\in[1/5,5/24]\), incompatible with the available local indices. A single Weierstrass point leaves only \(\xi=7/17\). At the unique quotient \(\frac1{17}(1,-1,8)\) the coordinate orders are \(14,3,10\), but
\[
7\notin\langle14,3,10\rangle.
\]
This contradiction excludes the quartic net and proves that \(P_4=2\). Proposition~2.1 then gives \(h=\chi=1\).

\needspace{4\baselineskip}

\(\square\)

\hypertarget{exhaustive-basket-calculations}{%
\hypertarget{section-3.3}{}\subsection{3.3. Exhaustive basket calculations}\label{exhaustive-basket-calculations}}

We give the arithmetic behind the preceding and subsequent lists. The first step is an a priori bound on every local denominator, so the enumeration requires no additional cutoff. For a normalized basket pair \((b,r)\), put
\[
\gamma=\min(r-2b,b),\quad d=b-\gamma,\quad
\epsilon=\frac b2-\frac{b^2}r-\frac\gamma6.
\]
Reid's formula, with \(P_2=1\), gives
\[
\sum b=10\chi+5-P_3,\qquad
\sum\gamma=5\chi-1+3P_3-P_4,\qquad
\sum\epsilon=\frac{\chi-2+P_4}{6}-V.
\]
The summands are nonnegative, and only \((1,2)\) and \((1,3)\) have \(\epsilon=0\). Enumerate the positive-gap pairs in a fixed order, allowing repetitions; a state with remaining integer budgets \(G,D\) is completed uniquely by \(G(1,3)+D(1,2)\). Its volume is the original gap budget minus the sum of the positive gaps. This finite recursion is exhaustive once the denominator bound is established;
all plurigenera are computed exactly by \eqref{eq:2.1}. For \(\chi=P_3=1\), \(P_4=2\),
the budgets are \((\sum\gamma,\sum d)=(5,9)\), and \cite[Lemma 4.9]{ref1} gives
\(r\le25\). \cite[Lemmas 3.12--3.13]{ref1} supply the Cartier-index and domain tests,
the latter also applying to the graded saturation of \(R/(s)\).

\hypertarget{lem-preliminary-volume}{}\ResultHeading{Lemma 3.3 (A preliminary volume bound).} Every remaining threefold satisfies \(V\ge1/48\).

\ProofHeading{Proof.} Before continuing, we show that our remaining threefolds satisfy \(V\ge1/48\). If \(V<1/50\), canonical restriction forces a \((1,2)\) quartic fibre. Fourteenth-canonical birationality gives \(5/14\le\xi\le20V\), hence \(V\ge1/56\). The recursion gives eight baskets in \([1/56,1/50)\):
\[
\begin{array}{c|c|c}
B&V&(P_5,P_6,P_7)\\\hline
3(1,2)+(3,7)+(4,9)+(4,11)&25/1386&(3,4,5)\\
7(1,2)+(1,3)+(3,7)+(3,10)&2/105&(2,4,4)\\
4(1,2)+(1,3)+3(2,5)+(3,7)&2/105&(3,3,4)\\
2(1,3)+(2,5)+(3,7)+(7,15)&2/105&(3,4,5)\\
4(1,2)+3(1,3)+(1,4)+(6,13)&1/52&(2,4,4)\\
2(1,2)+2(3,8)+(6,13)&1/52&(3,4,5)\\
(1,2)+2(1,3)+(5,12)+(6,13)&1/52&(3,4,5)\\
(2,5)+(4,11)+(8,17)&18/935&(3,4,5).
\end{array}
\]
The actual-index test excludes the first four. At index thirteen the only possible coefficient is \(5/13\), and the coordinate semigroup is \(\langle10,3,8\rangle\), which does not contain five. Thus only \((V,\xi)=(1/52,3/8)\) or \((18/935,4/11)\) remains.

In these two rows the sixth restriction has rank three and is generically finite of degree two. To justify this statement, we check the possible kernels of restriction. A larger kernel supplies another moving surface; the same outer pencil gives \(V\ge\xi/12\), a different inner pencil gives thirteen-canonical birationality, and the same inner pencil gives \((n+7)\xi\ge2+\lceil n\xi\rceil\), impossible at \(n=3\). Thirteen-canonical birationality itself contradicts \(13\xi<5\). If a rank-three restriction were composed with the canonical pencil, its moving divisor would contain two canonical curves and again give \(V\ge\xi/12\). Its generically finite degree is even by the marked genus-two semigroup, and its moving square is less than four; hence the degree is two.

The fifth generator supplies a parameter for the inner pencil. Together with the generators in degrees \(1,4,6\) it gives a polynomial subring \(B=\mathbb C[s,x_4,y_5,z_6]\) of field index two. The extra seventh section is invariant because \(7\xi<3\). If \(FW=G\) is its primitive relation and \(f=\deg F\), the Hilbert comparison gives \(V\ge(f+7)/420\). Both displayed volumes are below \(3/140\), so \(f=1\) and \(sW=x^2\) after subtracting an old section. Corollary~2.2 excludes this relation. The basket recursion has no \(\chi=1\) row in \([1/50,1/48)\).

\(\square\)

The recursion in \([1/48,1/30)\) now gives thirty-eight \(\chi=1\) baskets. Seven violate the domain inequality above, leaving thirty-one. These thirty-one necessary baskets are the common starting table for the remaining fibre arguments. They are not asserted to be realizable.

\needspace{4\baselineskip}

\hypertarget{the-possible-quartic-fibres}{%
\hypertarget{section-3.4}{}\subsection{3.4. The possible quartic fibres}\label{the-possible-quartic-fibres}}

\hypertarget{prop-quartic-fibre-types}{}\ResultHeading{Proposition 3.4 (The possible quartic fibres).} The minimal general fibre of the quartic pencil has type \((1,2,0)\), \((2,2,0)\), or \((2,3,0)\).

\ProofHeading{Proof.} Let \(F\) again denote a general quartic fibre. We have
\[
D\ge_{\mathbb Q}A/5,\qquad
V\ge D^2/4\ge AD/20\ge A^2/100.
\]
Applying Noether's and Debarre's inequalities, together with restriction of absolute two-forms, we are left with
\[
(K_S^2,p_g(S),q(S))\in
\{(1,2,0),(2,2,0),(2,3,0),(3,2,0),(3,3,0),(2,1,1),(3,1,1)\}.
\]
We remove the last four possibilities.

If \(q(F)=1\), \cite[Lemma 3.2]{ref1} gives a
nonzero two-form restricting to zero on \(F\). \cite[Lemma 3.16]{ref1} then gives
\(V\ge1/20\), a contradiction.

For a free curve system \(C\) with \(A\ge_{\mathbb Q}C\), \(AC\ge3\) and \(g(C)\ge4\), put \(\xi=DC\ge3/5\). The curve restriction inequality at \(m=12\) gives
\[
12\xi\ge6+\lceil2\xi\rceil\ge8,\qquad
V\ge\xi/20\ge1/30.
\]
For a minimal \((3,3)\)-surface, the canonical net is generically finite, has no fixed curve, and has at most one simple base point; see \cite[Proposition 5.1(2)--(3)]{ref20}. Its general resolved member has genus four and canonical degree three, so this type is excluded.

For a \((3,2)\)-surface write \(K_S=M+Z\). Hodge index and parity leave
\[
(K_SM,M^2,g)=(2,0,2),(3,1,3),(3,3,4).
\]
The last case is already excluded. A generically finite free system \(N\) on a blowup of \(S\) has \(AN\ge4\): otherwise pushing down gives \(K_SN_0=N_0^2=3\), so \(N_0\equiv K_S\); the torsion-freeness theorem for \(K_S^2=2\chi(\mathcal O_S)-3\) \cite[Lemma 1]{ref14} and \(q(S)=0\) imply \(N_0\sim K_S\), contradicting the existence of at least three independent sections. But \(A(M_6|_F)\le6AD<4\). The sixth restriction contains the canonical system by vanishing, and is therefore the canonical pencil of rank two.

The five numerical rows with \(3/100\le V<1/30\) all have \(P_6=4\). Thus the sixth image is a nondegenerate surface \(\Sigma\subset\mathbb P^3\). If \(d=\deg\Sigma\) and \(C\) is its connected curve fibre,
\[
36V\ge Q M_6^2=d\,DC,\qquad d\ge2.
\]
The genus-three case has \(DC\ge3/5\) and is impossible. In the genus-two case \(DC\ge2/5\), so \(d=2\). If \(\Sigma\) is smooth, its two rulings give \(M_6\sim F+G\), with \(G|_F=C\), and
\[
6V\ge Q^2F+Q^2G\ge\frac3{25}+\frac{DC}{4}\ge\frac{11}{50}.
\]
If \(\Sigma\) is a cone, its resolution is \(\mathbb F_2\), whose hyperplane class is \(E+2\ell\). Thus \(M_6\ge2F\) and \(6V\ge6/25\). The outer ruling has degree one because its parameter is already a ratio of sixth sections. Both cases contradict \(V<1/30\).

These exclusions leave precisely the three fibre types that we shall treat in the rest of the proof:
\[
\boxed{(1,2,0),\quad(2,2,0),\quad(2,3,0).}
\]

\needspace{4\baselineskip}

\(\square\)

\hypertarget{exclusion-of-the-230-fibres}{%
\needspace{10\baselineskip}\hypertarget{section-4}{}\section{\texorpdfstring{4. Exclusion of the \((2,3,0)\) fibres}{4. Exclusion of the (2,3,0) fibres}}\label{exclusion-of-the-230-fibres}}

\needspace{4\baselineskip}

\hypertarget{the-common-geometry-of-a-23-fibre}{%
\hypertarget{section-4.1}{}\subsection{\texorpdfstring{4.1. The common geometry of a \((2,3)\) fibre}{4.1. The common geometry of a (2,3) fibre}}\label{the-common-geometry-of-a-23-fibre}}

\hypertarget{prop-horikawa-common}{}\ResultHeading{Proposition 4.1 (Common Horikawa geometry).} In the \((2,3,0)\)-fibre case, \(1/2\le a<2/3\), \((P_5,P_6,P_7)=(3,4,5)\), the sixth evaluation bundle splits as \((1,0,0)\), \([\mathbb C(W):L_6]=2\), \(L_7=L_6\), and \(V\ge1/35\). The four baskets displayed below are the only remaining numerical possibilities.

\ProofHeading{Proof.} Let \(S\) be a \((2,3)\)-surface. We use the classical double-plane description of its canonical model \cite[Theorem 3.7(1)]{ref20}:
\[
R(S,K_S)=\mathbb C[x_0,x_1,x_2,w]/(w^2-F_8),\qquad
\deg x_i=1,\quad\deg w=4.
\]
The branch octic is reduced and has ADE singularities. We shall use actual divisors, including reducible canonical members.

We first observe that a generically finite free system \(N\) on a blowup of \(S\), with \(AN<4\), must be the complete canonical net. Indeed, the pushdown \(N_0\) has \(N_0^2\ge2\). If \(K_SN_0=2\), Hodge index gives \(N_0^2=2\) and \(N_0\equiv K_S\); algebraic simple connectedness \cite[Proposition 3.5]{ref20} and \(q(S)=0\) imply linear equivalence, and no moving square is lost. If \(K_SN_0=3\), parity gives \(N_0^2=3\). Riemann--Roch gives \(h^0(N_0)\ge4\). After at most one simple base point is resolved, its general member has genus four, while the induced linear series has degree at most three and at least three sections, contrary to Clifford's theorem.

We can apply this observation to the sixth restriction: it contains \(|K_F|\) and satisfies \(AN_6<4\), so it is the complete canonical net. For a general canonical genus-three curve, the numerical restriction inequality gives, for every sufficiently large \(m\),
\[
m a\ge4+\lceil(m-10)a\rceil,\qquad a=AD.
\]
\cite[Lemma 2.4]{ref1} therefore gives
\[
\frac12\le a\le20V<\frac23.
\]

We next determine the evaluation bundles over the base of the quartic pencil. This will allow us to turn the number of global sections into a volume estimate. If a fibre is \(\sum_j n_{pj}T_{pj}\) and \(a_{pj}=\operatorname{ord}_{T_{pj}}s\), set \(c_p=\min_j a_{pj}/n_{pj}\). By comparing orders along the components of each fibre, we find that the saturation of the line generated by \(s^m\) is
\[
\mathcal L_m=\mathcal O_{\mathbb P^1}\Bigl(\sum_p\lfloor mc_p\rfloor p\Bigr).
\]
The low plurigenera imply \(c_p<1/3\) and exactly one point \(p_0\) with \(1/4\le c_{p_0}<1/3\). An absolute two-form supplies a primitive relative canonical section at every fibre: divisibility by a base parameter would make its restriction zero on every component, contradicting the positive constant integral of \(\alpha\wedge\bar\alpha\) over fibres. Multiplying by the saturated relative canonical line gains one degree at every point with \(c_p>0\). Two such points would give a second bicanonical section. Therefore \(p_0\) is the only positive point, and
\[
\mathcal L_4=\mathcal L_6=\mathcal O(1).
\]

Let \(\mathcal V_m\) be the image of evaluation of \(R_m\) in \(f_*\omega_W^m\). It is globally generated, has \(P_m\) sections, and in degree six has rank three. Semipositivity, Leray, and the fixed-part decomposition give
\[
f_*\omega_{W/\mathbb P^1}=\mathcal O(2)\oplus\mathcal O(1)\oplus\mathcal O.
\]
Multiplication by \(s^4\) extends from \(f_*\omega_{W/\mathbb P^1}(-2)(p_0)\) to \(f_*\omega_W^5\), proving \(P_5\ge3\). Corollary~2.2 gives \(P_6>P_5\), so \(P_6\ge4\).

Suppose, for a contradiction, that \(P_6\ge5\). If the evaluation bundle has a splitting summand of degree at least two, then \(6Q\ge_{\mathbb Q}2F\), and canonical restriction gives \(V\ge1/24\). Otherwise all splitting degrees are at most one, and the splitting is \((1,1,0)\) or \((1,1,1)\). The second ruling of the resulting quadric scroll, or the Segre decomposition, gives \(M_6\ge F+G\) with \(A(G|_F)\ge2\). Resolving the vertex retains an effective exceptional term, and
\[
6V\ge D^2+Q^2G\ge1/10+1/10=1/5,
\]
a contradiction. Hence \(P_6=4\), \(\mathcal V_6=(1,0,0)\), and the sixth map has degree two onto \(\mathbb P^3\).

For the seventh restriction, \(AN_7\le7a<14/3\), so \(AN_7\le4\). A general canonical genus-three curve is hyperelliptic, and no series of degree at most four is birational. Since the sixth field has index two, this proves
\[
L_7=L_6,\qquad[\mathbb C(W):L_6]=2.
\]

The basket recursion in \([1/40,1/30)\), with \(P_5\ge3\) and \(P_6=4\), gives \((P_5,P_6,P_7)=(3,4,5)\).

Finally, \(P_7\ge5\) implies \(V\ge1/35\). Otherwise \(AN_7<4\) and the rigidity just proved gives \(N_7=A\) and evaluation rank three. The nonnegative splitting degrees sum to at least two. A degree-two summand yields \(7Q\ge2F\) and \(7V\ge1/5\). If all degrees are at most one, the second ruling of the quadric scroll or the Segre decomposition gives \(M_7\ge F+G\) with \(Q^2G\ge1/10\), and again \(7V\ge1/5\).

Applying \eqref{eq:whole-member-trace} to a general fourth-canonical
member and its coefficient-one moving component gives
\[
P_{10}-P_6\ge P_2(F)=K_S^2+\chi(\mathcal O_S)=6.
\]
Combining this inequality with the \(\chi=1\) basket recursion leaves exactly
\[
\begin{array}{c|c}
B&V\\\hline
3(1,2)+(1,3)+(5,11)+(5,13)&25/858\\
3(1,2)+(1,3)+(2,5)+(3,8)+(5,11)&41/1320\\
2(1,3)+2(2,5)+(8,17)&8/255\\
4(1,2)+(2,5)+(4,9)+(4,11)&16/495.
\end{array}
\]
We exclude these four baskets in Sections 4.2--4.4.

\needspace{4\baselineskip}

\(\square\)

\hypertarget{horikawa-fibres-with-one-absolute-two-form}{%
\hypertarget{section-4.2}{}\subsection{4.2. Horikawa fibres with one absolute two-form}\label{horikawa-fibres-with-one-absolute-two-form}}

We prove that the four baskets at the end of the preceding subsection cannot occur with a \((2,3)\) fibre and \(\chi=1\). Put \(t=x/s^4\) for the quartic parameter. The fixed-part decomposition \cite[Theorem 1.1]{ref4} identifies the unique absolute two-form \(\alpha\) with the trivial summand of
\[
f_*\omega_{W/\mathbb P^1}=\mathcal O(2)\oplus\mathcal O(1)\oplus\mathcal O.
\]
Let \(v\) be defined by \(\iota_v s=\alpha\) and use a prime for \(v\).

\hypertarget{lem-horikawa-detector}{}\ResultHeading{Lemma 4.2 (A reduced transverse detector).}
In the Horikawa case, let \(0\ne\alpha\in H^0(W,\Omega_W^2)\), let \(\iota_vs=\alpha\), and choose \(x\in R_4\setminus sR_3\). There is a prime divisor \(T\) on \(W\), nonexceptional over the terminal model, such that
\[
 \operatorname{ord}_T(s)=1,\qquad x|_T\ne0,\qquad
 \alpha|_T\ne0,\qquad
 \operatorname{mult}_T F_\infty=4,\qquad
 \operatorname{coeff}_TQ=1.
\]
Here \(F_\infty\) is the fibre of \(t=x/s^4\) over infinity. Its canonical saturation coefficient is \(1/4\). At this component the exact pole orders of \(t,t',t''\) are \(4,6,8\), and
\[
 t''-24\lambda_{\mathrm{pol}} t^2\in E_7,\qquad \lambda_{\mathrm{pol}}\ne0.
\]

\ProofHeading{Proof.}
Choose the reduced visible prime supplied by Corollary~2.2 and take its strict transform on \(W\). It survives on the canonical model and has zero discrepancy over the terminal model. Since \(\theta_\alpha^2=\lambda_{\mathrm{pol}}\bar x\) and \(\theta_\alpha\) is a unit there, \(x\) is a unit and the canonical multiplicity is one. The rational function \(t=x/s^4\) has pole four, so \(T\) occurs with multiplicity four in \(F_\infty\). For every component \(U\) of that fibre, the quartic fixed-part inequality gives \(\operatorname{mult}_U F_\infty\le4\operatorname{coeff}_UQ\); equality holds at \(T\). Thus the saturation coefficient is exactly \(1/4\). The local transverse calculation gives the displayed pole orders, and \eqref{eq:2.8} with \(m=4\) gives the last identity. \(\square\)

In particular the canonical saturation is \(\mathcal O(1)\) through degree seven. Fix this component throughout the Horikawa calculations. Thus the selected transverse multiplicity is \(e=1\), and the coefficient \(k\) of \(t^2\) in the equation for \(z'=t''\) is \(24\lambda_{\mathrm{pol}}\ne0\).

\hypertarget{lem-horikawa-field-package}{}\ResultHeading{Lemma 4.3 (The Horikawa field and branch package).} In the setting of Lemma~4.2, one may use coordinates \(t,r,z=t'\) for the quadratic sixth field. The actual fibre divisor is \(D=a\Gamma+D_{\mathrm{null}}\), with \(A\Gamma=1\) and \(AD_{\mathrm{null}}=0\). Every invariant normalized \(m\)-canonical section is polynomial in \(r,z\) over \(\mathbb C(t)\), of fibre degree at most \(\lfloor ma/2\rfloor\). The branch polynomial has degree seven and completes to a reduced ADE octic. The primitive seventh denominators are \(s,s^2,s^3,s^4,s^5,s^6,L,sL,s^2L\), where \(L=x-bs^4\); the first and third are impossible.

\ProofHeading{Proof.} Choose a fifth generator \(y\) and set
\[
r=y/s^5,\qquad z=t'=z_6/s^6.
\]
Then
\[
E_5=\langle1,t,r\rangle,\quad E_6=\langle1,t,r,z\rangle,\quad
K=L_6=\mathbb C(t,r,z),\quad[\mathbb C(W):K]=2.
\]
The three coordinates are algebraically independent because the sixth map is generically finite. The involution negates \(s\) and \(\alpha\): its restriction to a general fibre is the canonical double-plane involution, and restriction of absolute forms is injective. Thus \(v\) preserves \(K\). If \(h\) is the additional seventh function, the Euler identities give the exact system
\[
\begin{aligned}
t'&=z,\\
r'&=\lambda h+Bz+Cr+Dt+E,\\
z'&=k t^2+\mu h+Gz+Hr+It+J_0,
\end{aligned}
\]
with constant coefficients. At the selected component, \(e=1\) and \(k=24\lambda_{\mathrm{pol}}\ne0\). Also \(\operatorname{pole}_T(r)\le5e\); differentiation raises a positive pole \(n\) to exactly \(n+e+1\), and a logarithmic derivative has pole at most \(e+1\).

We now determine the divisor on the double plane that carries the poles of our normalized sections. If its distinguished line were not a branch component, an order-one horizontal canonical component would have coefficient \(\eta=a/2\in[1/4,1/3)\) in \(Q\). At a \(cA/r\) centre write the index-one equation as
\[
xy+g(z^r,u)=0,\qquad\mu_r(1,-1,b,0),\qquad0<b\le r/2.
\]
The residue volume form has character \(b\). If the inertia order is \(e_0\mid r\) and \(k_0=e_0\eta\), then the coefficient of \(s\) has value \(k_0\) and character \(-b\). Put \(v_0=\langle k_0b^{-1}\rangle_{e_0}\). The coordinate values are at least
\[
(e_0-v_0,v_0,e_0-k_0,e_0).
\]
Since \(k_0<e_0/3\), a leading canonical monomial can only be \(y^b\), with \(k_0=bv_0\). Horizontality requires a fourth-canonical numerator with nonconstant residue after division by this coefficient to the fourth power. A nonpure initial monomial must be \(y^ix\), \(y^iz\) or \(y^iu\). Its value and full character respectively give
\[
\begin{array}{c|c}
i=4b+1-e_0/v_0&i\equiv4b+1\pmod r\\
i=5b-e_0/v_0&i\equiv5b\pmod r\\
i=4b-e_0/v_0&i\equiv4b\pmod r.
\end{array}
\]
Hence \(e_0/v_0\) is a positive multiple of \(r\). Since \(e_0\le r\), it follows that \(e_0=r\), \(v_0=1\) and \(\eta=b/r\). Every constituent of the four baskets has \(b/r\ge1/3\), a contradiction. Non-\(cA\) terminal centres have index two, three or four; the first two cannot clear a coefficient in the stated interval, and the fourth would contribute an absent \((1,4)\) constituent. Thus the distinguished line is a branch component.

Let \(\Gamma\) be its reduced inverse image. The actual support statement is
\[
D=a\Gamma+D_{\mathrm{null}},\qquad A\Gamma=1,\qquad AD_{\mathrm{null}}=0.
\]
On a general canonical genus-three curve the marking is Weierstrass, with semigroup \(\langle2,7\rangle\). This marking also controls every invariant function \(f\in E_m\cap K\): it is a polynomial in \(r,z\) over \(\mathbb C(t)\) and
\[
\deg_{r,z}f\le\left\lfloor\frac{ma}{2}\right\rfloor.
\]
Indeed \(\operatorname{div}(f|_F)+mD\ge0\) follows by pulling back the actual effective divisor of the pluricanonical section. A pole along any plane curve other than the distinguished line would lift to an \(A\)-positive curve outside \(\Gamma\), which is impossible. The descended function is therefore regular on the affine plane. A polynomial of degree \(d\) has pole \(2d\) along \(\Gamma\), giving the inequality. In particular the derivative \(h'\in E_9\cap K\) has fibre degree at most two.

The normalized Jacobian
\[
J=\frac{dt\wedge dr\wedge dz}{s}
\]
is nonzero and anti-invariant, and \(\mathbb C(W)=K(J)\). On the affine canonical plane, the branch equation is a residual septic; hence
\[
P:=J^2\in\mathbb C(t)[r,z],\qquad\deg_{r,z}P=7.
\]
Its homogenization together with the distinguished line is a reduced ADE octic. In particular that octic has multiplicity at most three at each point. Closedness of \(\alpha\) gives
\[
\frac{J'}J=\operatorname{div}v,\qquad P'=2(\operatorname{div}v)P.
\]

The identity \(v(P)=\Lambda P\) says that the branch polynomial is a Darboux polynomial. We use this standard invariant-hypersurface language as in \cite[Proposition 1]{refPereira}; the weighted restrictions and the specific branch classifications required below are proved in the argument.

It follows that the primitive denominator of the extra seventh section has only the factors \(s\) and factors \(x-bs^4\): the fibrewise polynomiality just proved leaves poles only over points of the outer base. The Hilbert comparison gives
\[
V\ge\frac{7+\deg F}{420},\qquad\deg F\le6.
\]
Consequently, we need only consider the denominators
\[
s,s^2,s^3,s^4,s^5,s^6,\qquad L,sL,s^2L,\quad L=x-bs^4.
\]
An \(s\) denominator gives \(h=t^2\) after subtracting an old section, contradicting its pole \(8e>7e\). An \(s^2\) denominator gives \(h=t(r+ct)\); it is treated below. For \(s^3\), normality excludes a numerator containing no \(xz_6\) term, since completing the \(y^2\) term would make \(y/s\) integral of degree four. Thus its normalized form is
\[
h=tz+br^2+ctr+dt^2.
\]
Cancellation of the pole of \(tz\) forces \(b\ne0\) and \(\operatorname{pole}_T(r)=(9e+1)/2>4e\). The divergence is
\[
2\lambda br+(\lambda c+\mu)t+C+G.
\]
Its logarithmic-derivative bound forces \(\lambda=0\). But then \(r'\) has pole \((11e+3)/2>5e+1\), whereas its affine-linear right side has pole at most \(5e+1\), a contradiction.

\(\square\)

\needspace{3\baselineskip}

\hypertarget{result-4.1}{}\ResultHeading{Proposition 4.4 (finite modification).} In the Horikawa setup above, an additional normalized pluricanonical section
\[
h=\frac{u}{t},\qquad u=c r+d z+e_0\ne0,
\tag{4.1}\label{eq:4.1}
\]
cannot satisfy a differential system of the form
\[
\begin{aligned}
t'&=z,\\
r'&=A_0h+Bz+Cr+Dt+E,\\
z'&=k t^2+F_0h+Gz+Hr+It+J_0,
\end{aligned}
\tag{4.2}\label{eq:4.2}
\]
with all coefficients constant. A prime denotes \(v\). The assertion holds for \(h\in E_m\) in any positive degree \(m\) for which these hypotheses hold.

\ProofHeading{Proof.} Closedness of \(\alpha\) gives \(\mathcal L_v s=d\alpha=0\). Taking divergence in the independent coordinates \((t,r,z)\) therefore yields
\[
\frac{v(J)}J=\lambda_0+\frac{\lambda_1}{t},\qquad
\lambda_0=C+G,\quad\lambda_1=A_0c+F_0d.
\tag{4.3}\label{eq:4.3}
\]
We do not set \(\lambda_1\) equal to zero. Since \(P=J^2\),
\[
v(P)=2(\lambda_0+\lambda_1/t)P.
\tag{4.4}\label{eq:4.4}
\]

Expand the rational coefficients of \(P\) at \(t=0\). Let \(\nu\) be their minimum order and write
\[
P=t^\nu\bigl(Q(r,z)+tQ_1(r,z)+\cdots\bigr),
\qquad 0\ne Q\in\mathbb C[r,z].
\tag{4.5}\label{eq:4.5}
\]
The coefficient of \(t^{\nu-1}\) in \eqref{eq:4.4} is
\[
\nu zQ+u(A_0Q_r+F_0Q_z)=2\lambda_1Q.
\tag{4.6}\label{eq:4.6}
\]
If \(n=\deg Q\), its first term has degree \(n+1\) when \(\nu\ne0\). The other terms have degree at most \(n\), because \(u\) is affine linear. Hence \(\nu=0\).

Thus \(P\) is regular over the generic point of \(t=0\) and has nonzero residue there. Consider the DVR
\[
\mathcal O=\mathbb C(r,z)[t]_{(t)}
\]
in \(L_6\). The equation \(J^2=P\) defines an unramified quadratic extension at this valuation: \(P\) and \(2J\) are units. At any prime above it, \(t\) is a uniformizer, \(r,z\) remain independent residue coordinates, and
\[
s=\frac{dt\wedge dr\wedge dz}{J}
\]
has order zero. The nonzero polynomial \(u\) is also a unit at this generic point. Therefore
\[
\operatorname{ord}(s^m h)=-1.
\tag{4.7}\label{eq:4.7}
\]
This is impossible for a global pluricanonical section on a smooth projective model. Indeed, realize this divisorial valuation on a projective model of the function field and resolve the birational map to \(W\). The pullback of a holomorphic pluricanonical form is holomorphic, whereas \eqref{eq:4.7} is its intrinsic order at that prime.

This argument uses only a smooth neighbourhood of the generic divisor of the affine double cover. No canonicity claim about its other singularities is needed. \(\square\)

\needspace{3\baselineskip}

\hypertarget{result-4.2}{}\ResultHeading{Proposition 4.5.} In the Horikawa setup above, no additional seventh section can have primitive relation
\[
s^2u_7=xy,
\qquad\text{equivalently } h=u_7/s^7=tr,
\tag{4.8}\label{eq:4.8}
\]
after allowed changes of generators. This exclusion does not require a common-field statement in degree eleven, twelve or fourteen.

\ProofHeading{Proof.} Write the exact Euler system as
\[
\begin{aligned}
t'&=z,\\
r'&=a_1tr+bz+cr+dt+e_0,\\
z'&=k t^2+A_1tr+Bz+Cr+Dt+E_0,
\end{aligned}
\tag{4.9}\label{eq:4.9}
\]
with complex constant coefficients. Its divergence identity is
\[
\frac{v(J)}J=a_1t+c+B.
\tag{4.10}\label{eq:4.10}
\]
At the active prime, the logarithmic-derivative bound and the pole \(4e>e+1\) of \(t\) force \(a_1=0\).

Since \(tr\in E_7\), one has \(\operatorname{pole}_T(r)\le3e\). Thus \(r'\) has pole at most \(4e+1\), including when \(r\) is regular. The term \(bz\) would have pole \(5e+1\), strictly above every other term in its equation. Consequently \(b=0\), and
\[
r'=cr+dt+e_0.
\tag{4.11}\label{eq:4.11}
\]
We do not need to prove \(B+c=0\). In particular the mixed-Hodge finite-dimensional obstruction is not required for the new local argument.

\needspace{4\baselineskip}
\textbf{The case \(d\ne0\).}

Equation \eqref{eq:4.11} forces \(r\) to have positive pole
\[
\operatorname{pole}_T(r)=3e-1.
\tag{4.12}\label{eq:4.12}
\]
Indeed, a regular \(r\) cannot produce the pole \(4e\) of \(dt\), and comparison of the leading poles of \(r'\) and \(dt\) gives \eqref{eq:4.12}.

We have \(e=1\) and \(k\ne0\) by Lemma~4.2.

Put \(f=r\). Because \(d\ne0\),
\[
t=(f'-cf-e_0)/d,\qquad z=(f''-cf')/d.
\tag{4.13}\label{eq:4.14}
\]
Elimination from \eqref{eq:4.9} gives the exact equation
\[
\begin{aligned}
f'''={}&(c+B)f''+\frac{k}{d}(f'-cf-e_0)^2
 +A_1f(f'-cf-e_0)\\
&-Bc f'+Cd f+D(f'-cf-e_0)+dE_0.
\end{aligned}
\tag{4.14}\label{eq:4.15}
\]

For \(e=1\), the actual-time relation and \eqref{eq:4.12} give pole one for \(f\) in local time. Its dominant equation is
\[
f'''=(k/d)(f')^2.
\]
Every other term of \eqref{eq:4.15} has smaller pole. The indicial roots are \(-1,1,6\). The first-nonintegral-coefficient argument therefore makes \(f\) an integral Laurent series in \(\tau\). Formula \eqref{eq:4.14} gives the same conclusion for \(t,z\).

This contradicts the local-time obstruction.

\needspace{4\baselineskip}
\textbf{The case \(d=0\).}

The first-order equation \(r'=cr+e_0\) makes \(r\) holomorphic in local time, with an ordinary power series whose initial value can depend on the transverse coordinates. This follows either from its elementary solution or from the absence of any possible polar or fractional leading term.

Now \(t\) has time-pole two and satisfies
\[
t''=k t^2+B t'+(A_1r+D)t+Cr+E_0.
\tag{4.15}\label{eq:4.16}
\]
The coefficients involving \(r\) are ordinary power series in \(\tau\). The leading term is \(6k^{-1}\tau^{-2}\), and the indicial roots are \(-1,6\). Thus \(t\), \(z=t'\) and \(r\) all have integral Laurent expansions. The local-time obstruction gives the final contradiction. \(\square\)

\hypertarget{prop-horikawa-seventh-rank}{}\ResultHeading{Proposition 4.6 (The seventh restriction).} In the Horikawa setting, \(\rho_7=4\), \(\mathcal V_7=\mathcal O(1)\oplus\mathcal O^3\), \(a\ge4/7\), and \(AN_7=4\).

\ProofHeading{Proof.} Suppose \(\rho_7=3\). Global generation, \(h^0(\mathcal V_7)=5\), and the saturated canonical \(\mathcal O(1)\) give
\[
\mathcal V_7=\mathcal O(1)^2\oplus\mathcal O.
\tag{4.16}\label{eq:4.17}
\]
The other globally generated rank-three splitting with five sections is \(\mathcal O(2)\oplus\mathcal O^2\). It cannot contain a saturated \(\mathcal O(1)\): every map to a trivial summand is zero, and the map to \(\mathcal O(2)\) has a zero. The inclusion \(\mathcal V_6\hookrightarrow\mathcal V_7\) therefore has a quotient of length one. The quotient is a skyscraper sheaf of length one; in a local frame its new section is an old fibrewise-linear section divided by a parameter.

At a finite point \(t=b\), subtracting old sections and shifting the coordinate gives
\[
h=(c r+d z+e_0)/t.
\]
The Euler identities give exactly \eqref{eq:4.2}, so \hyperlink{result-4.1}{Proposition 4.4} excludes this possibility.

At \(t=\infty\), the extra function is a nonzero combination of \(t^2,tr,tz\) modulo \(E_6\). A nonzero \(tz\) coefficient has pole \(9e+1\), while \(tr\) has pole at most \(9e\) and \(t^2\) has pole \(8e\). Such a unique pole is incompatible with membership in \(E_7\). The \(tz\) coefficient is zero. A pure \(t^2\) is likewise impossible because \(8e>7e\). The remaining combination is a nonzero scalar multiple of
\[
tR,\qquad R=r+bt,
\]
after a change of fifth generator. It is the denominator \(s^2\) excluded in \hyperlink{result-4.2}{Proposition 4.5}. Thus \(\rho_7\ne3\).

\needspace{4\baselineskip}
\textbf{The rank and the coefficient.}

The sixth canonical net gives \(\rho_7\ge3\). The globally generated \(\mathcal V_7\) has five sections and contains the canonical \(\mathcal O(1)\), so its rank is at most four. Hence
\[
\rho_7=4,\qquad \mathcal V_7=\mathcal O(1)\oplus\mathcal O^3.
\tag{4.17}\label{eq:4.18}
\]
By small-system rigidity, \(a<4/7\) would force rank three. This proves \(a\ge4/7\).

The value \(a=4/7\) is retained. On a resolution of the seventh restriction, its moving divisor \(N_7\) satisfies \(AN_7\le7a<14/3\). The small-system rigidity and rank four force the integer intersection
\[
AN_7=4.
\tag{4.18}\label{eq:4.19}
\]
This is compatible with \(7a=4\); a strict inequality at this endpoint would be incorrect.

\needspace{3\baselineskip}

\(\square\)

\hypertarget{result-4.3}{}\ResultHeading{Lemma 4.7.} A system of the form
\[
t'=z,\quad r'=Bz+Cr+Dt+E,\quad
z'=kt^2+Gz+Hr+It+J_0
\tag{4.19}\label{eq:4.20}
\]
cannot satisfy the geometric package.

\ProofHeading{Proof.} Lemma~4.2 gives \(e=1\) and \(k\ne0\). Replace \(r\) by \(R=r-Bt\); then \(R'=cR+dt+e_0\).

If \(d\ne0\), \(R\) has a simple pole in \(\tau\), and eliminating \(t,z\) gives
\[
R'''=(k/d)(R'-cR-e_0)^2+\text{terms of smaller pole weight}.
\]
Its indicial polynomial is \((q+1)(q-1)(q-6)\). If \(d=0\), \(R\) is an ordinary analytic function of \(\tau\), and the remaining second-order equation for \(t\) has leading term \(t''=kt^2\) and indicial polynomial \((q+1)(q-6)\). In both cases the first-nonintegral-coefficient test gives integral Laurent expansions for all three coordinates. We then obtain a contradiction from the local-time obstruction. \(\square\)

\needspace{4\baselineskip}

\hypertarget{pure-denominators-in-the-horikawa-case}{%
\hypertarget{section-4.3}{}\subsection{4.3. Pure denominators in the Horikawa case}\label{pure-denominators-in-the-horikawa-case}}

\hypertarget{lem-horikawa-fourth}{}\ResultHeading{Lemma 4.8 (The pure fourth denominator).} A primitive seventh denominator \(s^4\) is impossible in the Horikawa setting.

\ProofHeading{Proof.} For \(s^4u=G\), completing the product and subtracting an old section gives
\[h=RZ+\varepsilon t^2,\]
where \((R,Z)\) are affine linear in \((r,z)\) over \(\mathbb C[t]\), with constant nonzero fibre Jacobian. Write the two Euler equations as \(R'=\lambda h+\text{affine terms}\) and \(Z'=\mu h+kt^2+\text{affine terms}\). The fibre-degree-three part of \(h'\in E_9\) is
\[(\lambda Z+\mu R)RZ.\]
The invariant fibre-degree bound makes it zero. Hence \(\lambda=\mu=0\). Reversing the affine coordinate change gives the triangular system just excluded. Thus the case \(s^4\) reduces to \hyperlink{result-4.3}{Lemma 4.7} and is excluded.

\(\square\)

\hypertarget{lem-horikawa-cusp}{}\ResultHeading{Lemma 4.9 (The fifth-denominator normal form).} If a primitive denominator \(s^5\) occurs, the normalized relation and Euler system can be put in the form \eqref{eq:4.21}--\eqref{eq:4.22}, with \(\eta a_h\ne0\). Their actual leading balance is \eqref{eq:4.23}, and the logarithmic derivative satisfies \eqref{eq:4.24}.

\ProofHeading{Proof.} Assume \(s^5u=G\). The weight-twelve terms not divisible by \(s\) are \(z_6^2\) and \(x^3\). The coefficient of \(z_6^2\) cannot vanish: at the active component the surviving \(x^3\) term would be a unit, while every other term vanishes. Complete the square in the sixth generator. If the \(x^3\) coefficient were zero, the relation would have the form
\[
s^5u=Z_6^2+s^2H_{10}.
\]
Then \((Z_6/s)^2\in R_{10}\), and normality of the canonical ring would give \(Z_6/s\in R_5\), contradicting the independent sixth generator. Thus both coefficients are nonzero.

The reduced detector has \(e=1\). The completed normalized relation is
\[
h=Z^2+\eta t^3+\delta r^2+\varepsilon tr+\zeta t^2,
\quad Z=z+\beta r+\gamma t,
\quad\eta\ne0,
\tag{4.20}\label{eq:4.21}
\]
up to an element of \(E_6\). The highest-pole term in the divergence forces the coefficient of \(h\) in \(Z'\) to be zero. If its coefficient \(a_h\) in \(r'\) is also zero, the triangular-system lemma applies. Otherwise the system has the exact shape
\[
\begin{aligned}
t'&=Z-\beta r-\gamma t,\\
Z'&=kt^2+GZ+Hr+It+J_0,\\
r'&=a_hh+BZ+Cr+Dt+E,\qquad a_h\ne0.
\end{aligned}
\tag{4.21}\label{eq:4.22}
\]
Because \(r\) has time pole at most \(5/2\), its leading terms give
\[
t=\frac6k\tau^{-2}+o(\tau^{-2}),\quad
Z=-\frac{12}k\tau^{-3}+o(\tau^{-3}),\quad
\eta=-\frac{2k}{3}.
\tag{4.22}\label{eq:4.23}
\]
Indeed, membership in \(E_7\) forces the time-pole terms of order six in \(h\) to cancel, which gives the last identity. The logarithmic-derivative bound is
\[
2\delta r+\varepsilon t=O(\tau^{-1}).
\tag{4.23}\label{eq:4.24}
\]
In estimates with a nonintegral exponent, \(O\) refers to the order of the actual Puiseux expansion.

\(\square\)

\needspace{3\baselineskip}

\hypertarget{result-4.4}{}\ResultHeading{Lemma 4.10.} In \eqref{eq:4.22}, algebraic independence of \(t,r,Z\) forces
\[
\beta=H=0.
\tag{4.24}\label{eq:4.25}
\]

\ProofHeading{Proof.} We differentiate in a transverse variable while keeping \(\tau\) fixed, and denote the resulting variation of \((t,Z,r)\) by \((T,U,V)\). The leading coefficient of \(t\) in \eqref{eq:4.23} is constant. Thus, if \(T\ne0\), its first exponent \(n\) satisfies \(n>-2\). The variational equations are
\[
\begin{aligned}
T'&=U-\beta V-\gamma T,\\
U'&=(2kt+I)T+GU+HV,\\
V'&=a_h\{2ZU+(-2kt^2+\varepsilon r+2\zeta t)T\}\\
&\quad+\{a_h(2\delta r+\varepsilon t)+C\}V+BU+DT.
\end{aligned}
\tag{4.25}\label{eq:4.26}
\]
The coefficient of \(V\) in the third equation is \(O(\tau^{-1})\) by \eqref{eq:4.24}.

First suppose \(\beta\ne0\). Substitute \(V=(U-T'-\gamma T)/\beta\) into the second equation. If \(T=T_0\tau^n+\cdots\), its leading forcing is \(12T_0\tau^{n-2}\). The case \(n=1\) would require a logarithm and is impossible. Otherwise
\[
U=\frac{12T_0}{n-1}\tau^{n-1}+\cdots,
\tag{4.26}\label{eq:4.27}
\]
apart from a possible ordinary homogeneous solution with nonzero constant term. When \(n>1\), that constant term is ruled out by the unique order-three pole \(2a_hZU\) in the third equation. When \(n<1\), the displayed term already dominates it.

The terms \(V'\) and \(O(\tau^{-1})V\) have order at least \(n-2\). In contrast, the coefficient of \(\tau^{n-4}\) from the two leading terms \(a_h(2ZU-2kt^2T)\) is
\[
-\frac{72a_hT_0}{k}\frac{n+3}{n-1}.
\tag{4.27}\label{eq:4.28}
\]
All other terms have strictly larger exponent; for example \(rT\) has order at least \(n-5/2\). Thus \eqref{eq:4.28} must vanish, forcing \(n=-3\), contrary to \(n>-2\). If \(T=0\), the first two equations give an ordinary homogeneous equation for \(U\), and the third again forces \(U=V=0\). Therefore every transverse variation vanishes, contradicting the nonzero Jacobian \(J\).

Now \(\beta=0\). If \(H\ne0\), the first two equations give
\[
U=T'+\gamma T,\qquad
HV=T''+(\gamma-G)T'-(2kt+G\gamma+I)T.
\tag{4.28}\label{eq:4.29}
\]
For a first exponent \(n>-2\) of \(T\), \(V\) has order at least \(n-2\). The coefficient of \(\tau^{n-4}\) in the third equation is now
\[
-\frac{24a_hT_0}{k}(n+3),
\]
whereas \(V'\) and the terms involving \(V\) have order at least \(n-3\). This is impossible. The possible cancellation in the leading coefficient of \eqref{eq:4.29}, at \(n=4\), only increases the order of \(V\) and does not alter the contradiction. If \(T=0\), then \(U=V=0\). Again the Jacobian would vanish. Hence \(H=0\). \(\square\)

\hypertarget{prop-horikawa-fifth}{}\ResultHeading{Proposition 4.11 (The pure fifth denominator).} A primitive seventh denominator \(s^5\) is impossible in the Horikawa setting.

\ProofHeading{Proof.} The cubic fibre part of \(h'\) is \(2a_h\delta r(Z^2+\delta r^2)\). Since \(a_h\ne0\) and \(h'\) has fibre degree at most two, \(\delta=0\). The logarithmic-derivative bound then gives \(\varepsilon=0\). By \hyperlink{result-4.4}{Lemma 4.10}, we also have \(\beta=H=0\), so the equations for \(t,Z\) are independent of \(r\). The second-order equation for \(t\) has leading term \(t''=kt^2\) and indicial roots \(-1,6\). Thus \(t,Z\) have integral Laurent expansions in actual time. The equation for \(r\) is linear, with constant coefficient of \(r\) and integral Laurent forcing. A first nonintegral exponent cannot occur in this equation either. All three coordinates are integral in actual time, contradicting quadratic reconstruction.

\(\square\)

\hypertarget{prop-horikawa-sixth}{}\ResultHeading{Proposition 4.12 (The pure sixth denominator).} A primitive seventh denominator \(s^6\) is impossible in the Horikawa setting.

\ProofHeading{Proof.} Primitivity of \(s^6u=G\) forces a nonzero \(x^2y\) term. After scaling and subtraction of old sections,
\[
h=t^2r+\gamma t^3+\beta z^2+brz+ctz+dr^2+\varepsilon tr+\zeta t^2.
\tag{4.29}\label{eq:4.30}
\]
In the Euler system, the divergence contains \(a_ht^2\) and only lower-pole terms besides it. Since \(t^2\) has pole \(8e>5e+1\), the logarithmic bound forces \(a_h=0\). If \(b_h=0\), the triangular-system lemma excludes the system. Suppose \(b_h\ne0\). The divergence now gives
\[
\beta=0,\qquad br+ct\text{ has pole at most }e+1.
\tag{4.30}\label{eq:4.31}
\]
If \(b\ne0\), set \(R=r+(c/b)t\). It has pole at most \(e+1\). In \eqref{eq:4.30} the coefficient of \(t^3\) becomes \(\gamma-c/b\). A nonzero value would give a unique pole \(12e\), too large for \(h\in E_7\). Thus this coefficient vanishes, and
\[
h=t^2R+bRz+dR^2+\varepsilon_1tR+\zeta_1t^2
\pmod{E_6}.
\tag{4.31}\label{eq:4.32}
\]
If \(b=0\), \eqref{eq:4.31} forces \(c=0\). Set \(R=r+\gamma t\) instead. Equation \eqref{eq:4.32} holds with \(b=0\). Here \(R\) initially has pole at most \(5e\).

In either case, a positive pole \(n\) of \(R\) makes \(t^2R\) the unique highest-pole term of \eqref{eq:4.32}: its pole is \(8e+n\), while \(Rz\) has pole \(5e+1+n\) and \(R^2\) has pole \(2n\). The bounds on \(n\) just stated make both strictly smaller. Hence \(R\) is regular at \(T\).

Set \(U=R+\zeta_1\). After another subtraction from \(E_6\), \eqref{eq:4.32} becomes
\[
h=U(t^2+bz+\varepsilon_1t)+dU^2.
\tag{4.32}\label{eq:4.33}
\]
If \(U\) were a unit at \(T\), this would have pole \(8e>7e\). Thus \(U\) is regular and vanishes along \(T\).

Since \(a_h=0\), the derivative of \(U\) is a constant linear combination of \(1,U,t,z\). A function regular at \(T\) has derivative with pole at most \(e\), so the \(z\) and \(t\) coefficients must vanish. Thus \(U'=c_1U+c_0\), and the finite-dimensional Hodge lemma gives \(U'=0\). In actual-time coordinates, \(U\) is independent of \(\tau\); its zero restriction to \(T\) then implies \(U=0\) in a neighbourhood. This contradicts the algebraic independence of the coordinates. The denominator \(s^6\) is excluded.

\(\square\)

We have now excluded all six pure denominators, so it remains to treat those with a finite pole on the base.

\needspace{4\baselineskip}

\hypertarget{finite-denominators-in-the-horikawa-case}{%
\hypertarget{section-4.4}{}\subsection{4.4. Finite denominators in the Horikawa case}\label{finite-denominators-in-the-horikawa-case}}

\hypertarget{lem-horikawa-finite-normal-form}{}\ResultHeading{Lemma 4.13 (The finite-pole normal form).} For a primitive seventh denominator \(L\), the branch polynomial has a genuine pole at \(L=0\), and the normalized seventh function has the completed-product form \eqref{eq:4.36}.

\ProofHeading{Proof.} For \(F=L\), translate the base to make \(L=0\) equal to \(t=0\). Dividing the numerator by \(t\) and subtracting \(E_6\) gives
\[
h=H(r,z)/t,\qquad
\deg H\le2,\quad [z^2]H=0.
\tag{4.33}\label{eq:4.34}
\]
Here the \(z^2\)-term is absent because its weight is twelve, whereas the numerator has weight eleven. Let
\[
P=t^\nu(Q(r,z)+O(t)),\quad Q\ne0,\quad
\mathcal D=a_h\partial_r+b_h\partial_z.
\]
The lowest Laurent coefficient of \(v(P)=2\operatorname{div}(v)P\) is
\[
\nu zQ+H\mathcal DQ=2(\mathcal DH)Q.
\tag{4.34}\label{eq:4.35}
\]
Here \(\nu\) is an integer and \(\deg Q\le7\).

We claim that \(\nu<0\). The point is that \(s^7h\) is a global pluricanonical section. To compare its order with \(\nu\), we work over the generic DVR \(\mathbb C(r,z)[t]_{(t)}\). If \(\nu\) is even, the quadratic extension is unramified after removing the even power, and the order of \(s^7h\) is \(-7\nu/2-1\). If \(\nu\) is odd, the extension has ramification index two, and that order is \(7(1-\nu)-2\). Neither is nonnegative for \(\nu\ge0\). This is a divisorial calculation with the rational form \(s=dt\wedge dr\wedge dz/J\), justified on a common resolution; it assumes no canonicity of a global affine model.

If the \(rz\) coefficient of \(H\) is zero, write \(H=ar^2+br+cz+d\) and \(Q=\sum_{j=0}^M q_j(r)z^j\). The coefficient of \(z^{M+1}\) in \eqref{eq:4.35} is
\[
\nu q_M+ca_hq_M'=0.
\]
No nonzero polynomial \(q_M\) satisfies this when \(\nu\ne0\). Therefore the \(rz\) coefficient is nonzero. Completing the product gives
\[
h=(RZ+\kappa)/t,\qquad t'=Z-bR+c,
\tag{4.35}\label{eq:4.36}
\]
with \(R'=a_hh+\) affine terms and \(Z'=kt^2+b_hh+\) affine terms.

The reduced detector of Lemma~4.2 applies to \eqref{eq:4.36}.

\(\square\)

\needspace{3\baselineskip}

\hypertarget{result-4.5}{}\ResultHeading{Lemma 4.14.} For \eqref{eq:4.36}, a necessary condition is
\[
2a_h\in\tfrac12+\mathbb Z,\qquad 2a_h\le\tfrac52.
\tag{4.36}\label{eq:4.37}
\]

\ProofHeading{Proof.} The leading terms of \(t,Z\) are again \(6\tau^{-2}/k\) and \(-12\tau^{-3}/k\). If \(R\) has pole \(p>2\) in time, its first equation gives \(p=2a_h\); the section bound permits only \(p=5/2\).

Otherwise \(R\) has pole at most two. Use leading pole weights \((t,R,Z)=(2,2,3)\). At a first positive nonintegral shift \(q\), the two equations for \(t,Z\) have the usual diagonal indicial factor \((q+1)(q-6)\): terms involving \(R\) are of strictly lower pole weight. The last diagonal factor, from \(R'=a_hRZ/t+\text{lower terms}\), is
\[
q-2+2a_h.
\tag{4.37}\label{eq:4.38}
\]
The off-diagonal terms do not change these factors. If there were no nonintegral shift, the local-time obstruction would apply. Since all exponents are half-integral, \eqref{eq:4.38} requires \(2a_h\) to be a nonintegral half-integer; in the present pole range it is less than two. Combining this with \(p=5/2\) proves \eqref{eq:4.37}. In particular every case with \(2a_h\in\mathbb Z\) is excluded. \(\square\)

\needspace{3\baselineskip}

\hypertarget{result-4.6}{}\ResultHeading{Corollary 4.15 (The denominator \(L\)).} The primitive denominator \(L\) is impossible.

\ProofHeading{Proof.} For the completed finite-pole normal form, the fibre-degree-three part of \(h'\) is
\[
\frac{RZ}{t^2}\bigl((\lambda-1)Z+(\mu+b)R\bigr).
\]
The invariant fibre-degree bound forces \(\lambda=1\). The preceding local-time test requires \(2\lambda\) to be a nonintegral half-integer, which is impossible. The constant \(c_0\) and the lower-degree terms do not affect this calculation. \(\square\)

For the two mixed denominators, weighted numerator enumeration, followed by subtraction of an element of \(E_6\), gives
\[
h=\gamma t^2+\eta tr+\frac{H(r,z)}t,\qquad \deg H\le2,
\tag{4.38}\label{eq:4.39}
\]
where \(\eta=0\) for \(sx\). For \(s^2x\), primitivity gives \(\eta\ne0\): the degree-thirteen numerator modulo \(s\) is a nonzero multiple of \(x^2y\).

Write \(H=H_2+H_1+H_0\) by fibre degree and \(\mathcal D_0=\lambda\partial_r+\mu\partial_z\). Since \(\rho_7=4\) whereas \(E_6\) restricts to the span of \(1,r,z\), necessarily \(H_2\ne0\).

\needspace{3\baselineskip}

\hypertarget{result-4.7}{}\ResultHeading{Proposition 4.16 (Quadratic compatibility).} For either mixed denominator in \eqref{eq:4.39},
\[
\mathcal D_0H_2=z,\qquad
\mathcal D_0H-z=A_0\in\mathbb C.
\tag{4.39}\label{eq:4.40}
\]
Moreover, for \(s^2x\),
\[
\lambda=0,\quad\mu\ne0,\quad
H_2=\frac{z^2}{2\mu}+a_2r^2.
\tag{4.40}\label{eq:4.41}
\]

\ProofHeading{Proof.} The cubic fibre part of \(h'\) is
\[
\frac{H_2(\mathcal D_0H_2-z)}{t^2}.
\tag{4.41}\label{eq:4.42}
\]
The term \(-zH_2/t^2\), from differentiating \(1/t\), is indispensable. The invariant fibre-degree bound and \(H_2\ne0\) imply \eqref{eq:4.40}. Direct differentiation of the exact Euler system gives
\[
\operatorname{div}v=\lambda\eta t+\frac zt+\frac{A_0}t+c+g.
\tag{4.42}\label{eq:4.43}
\]
This equals \(v(J)/J\). At \(T\), its term \(\lambda\eta t\) has pole \(4e\), while \(z/t\) has pole \(e+1\) and \(A_0/t\) is regular. The logarithmic-derivative bound forces \(\lambda\eta=0\). In the weight-six case \(\eta\ne0\), so \(\lambda=0\); now \eqref{eq:4.40} gives \eqref{eq:4.41}. \(\square\)

For the weight-five denominator we will need the stronger identity \(\mathcal D_0H=z\). We prove it next; the quadratic identity in \hyperlink{result-4.7}{Proposition 4.16} already suffices in weight six.

\needspace{3\baselineskip}

\hypertarget{result-4.8}{}\ResultHeading{Proposition 4.17.} In the case \(F_{\rm den}=sx\), the constant \(A_0\) in \eqref{eq:4.40} vanishes.

\ProofHeading{Proof.} Suppose \(A_0\ne0\). Put \(g_9=\Delta u\in R_9\). In the polynomial ring on variables \(s,x,y,z_6,U,G_9\) of weights \((1,4,5,6,7,9)\), consider
\[
f_1=sxU-G_{12},\qquad
f_2=sxG_9-\Delta G_{12}+sz_6U.
\tag{4.43}\label{eq:4.44}
\]
Here \(\Delta G_{12}\) is computed from the exact Euler system, replacing \(u\) by \(U\). Primitivity makes \(f_1\) irreducible; the nonzero coefficient \(sx\) of \(G_9\) makes \(f_2\) a nonzero element of its domain quotient. Thus \eqref{eq:4.44} is a regular sequence of degrees \(12,14\), with total Hilbert leading coefficient
\[
e_{\rm total}=\frac{12\cdot14}{1\cdot4\cdot5\cdot6\cdot7\cdot9}
=\frac{168}{7560}.
\tag{4.44}\label{eq:4.45}
\]
Our convention is \(\dim T_m=e(T)m^3/6+O(m^2)\).

After inverting \(sx\), both new variables can be eliminated. There is a unique main component, generically reduced, whose ring is the actual subring \(B[u,g_9]\). Other components of dimension four must lie on \(sx=0\).

On \(x=0\), \(s\ne0\), the first relation is the homogenized equation \(H=0\). The coefficient of \(U\) in the second is \(-A_0s^7\), a unit. Thus the second relation is a genuine equation on each component of \(H=0\), and there is no four-dimensional component on this open boundary.

At \(s=0\), weight counting gives
\[
G_{12}\bmod s=\alpha_0x^3+\beta_0z_6^2,\qquad \beta_0\ne0.
\tag{4.45}\label{eq:4.46}
\]
Indeed, if \(\beta_0=0\), primitivity gives \(\alpha_0\ne0\). At the active prime \(x\) is a unit, so the two sides of \(sxu=G_{12}\) would have respectively positive and zero orders.

If \(\alpha_0\ne0\), the reduced boundary is \((s,\alpha_0x^3+\beta_0z_6^2)\), of weighted degree twelve. In the first hypersurface its generic local ring is a DVR with parameter \(s\). The coefficient \(sx\) of the independent variable \(G_9\) has order one, so \(f_2\) has order at most one. This boundary contributes at most \(12/7560\).

If \(\alpha_0=0\), the reduced boundary is \((s,z_6)\), of weighted degree six. The derivative of \(f_1\) with respect to \(s\) is generically nonzero, since it contains the independent term \(xU\). Its local equation makes \(s\) have order two in the parameter \(z_6\). The coefficient of \(G_9\) then has order two, bounding the boundary multiplicity by two. Its contribution is again at most \(12/7560\). A cancellation raising that order cannot occur identically in the independent variable \(G_9\).

There is no extra four-dimensional component at \(s=x=0\), since \eqref{eq:4.46} then forces \(z_6=0\). Additivity of Hilbert multiplicities over minimal primes therefore gives
\[
e(B[u,g_9])\ge\frac{168-12}{7560}=\frac{13}{630}.
\tag{4.46}\label{eq:4.47}
\]
The graded fraction field of this subring is \(K(\xi)\), while that of the canonical ring is \(\mathbb C(W)(\xi)\), of degree two. Choosing two homogeneous elements independent over the first field gives two shifted copies of the subring inside \(R\). Comparison of cumulative Hilbert functions yields
\[
V\ge2e(B[u,g_9])\ge\frac{13}{315}>\frac1{30}.
\]
This contradiction proves \(A_0=0\). \(\square\)

In this argument we subtract the contribution of the boundary before applying the volume comparison. This is why the estimate concerns the actual subring, even though the auxiliary complete intersection may have other components.

\hypertarget{lem-horikawa-darboux-reduction}{}\ResultHeading{Lemma 4.18 (Reduction to four Darboux parameters).} For the primitive denominator \(sx\), the septic branch polynomial satisfies \eqref{eq:4.51}; its leading weighted part lies in \(\ker\mathcal L_B\), with \(B\) one of the four values in \eqref{eq:4.53}.

\ProofHeading{Proof.} Assume \(F_{\rm den}=sx\). Write
\[
h=\gamma t^2+H/t,\quad
H=a_2r^2+b_2rz+c_2z^2+d_2r+e_2z+f_2.
\]
\hyperlink{result-4.8}{Proposition 4.17} gives
\[
2a_2\lambda+b_2\mu=0,\quad
b_2\lambda+2c_2\mu=1,\quad
d_2\lambda+e_2\mu=0.
\tag{4.47}\label{eq:4.48}
\]

\needspace{4\baselineskip}
\textbf{Active balance and the exact Darboux equation.}

Since \(h\) has pole at most seven along \(T\), the leading equation in \(z'\) is \(t''=kt^2\). Thus
\[
t=(6/k)\tau^{-2}+\cdots,\quad
z=-(12/k)\tau^{-3}+\cdots,\quad r=o(\tau^{-3}).
\tag{4.48}\label{eq:4.49}
\]
The two pole-eight terms of \(h\) cancel, giving
\[
\gamma=-2kc_2/3.
\tag{4.49}\label{eq:4.50}
\]
This follows from the stated pole bounds.

By \eqref{eq:4.43} and \(A_0=0\), the divergence is \(z/t+c+g\). Hence
\[
v(J/t)=(c+g)J/t.
\]
The Hodge consequence forces \(c+g=0\). Consequently the actual septic satisfies the exact equation
\[
v(P)=\frac{2z}{t}P.
\tag{4.50}\label{eq:4.51}
\]

\needspace{4\baselineskip}
\textbf{Why only four resonance parameters can occur.}

Give \(t,r,z\) pole weights \(2,3,3\). Linearization at the balance \eqref{eq:4.49}, using \eqref{eq:4.48}--\eqref{eq:4.50}, gives
\[
(q+1)\bigl((q-3)(q-4)-6B\bigr),\qquad B=\lambda b_2.
\tag{4.51}\label{eq:4.52}
\]
The root \(-1\) corresponds to translation of the pole. The other two roots have sum seven.

We need more than the existence of one positive fractional root. In local coordinates \((\tau,y_1,y_2)\), differentiate the actual solution in the two transverse directions at fixed \(\tau\). The leading coefficients in \eqref{eq:4.49} are constants, so the two scaled variations
\[
(\tau^2\delta t,\ \tau^3\delta r,\ \tau^3\delta z)
\]
have positive order. They are independent because \(J\ne0\). They solve a regular-singular variational system
\[
\tau U'=M(\tau^{1/2})U
\]
whose constant matrix has the eigenvalues in \eqref{eq:4.52}. The dimension of its logarithm-free Puiseux solutions of positive order is bounded by the total multiplicity of positive admissible indicial roots. To see this, compare the leading exponent and eigenvector; cancel equal leading vectors between solutions and repeat at the next nonzero exponent. Equivalently, the coefficient recursion supplies free data only at indicial roots. Fewer than two such freedoms would make both transverse variations dependent.

Both nontranslation roots must therefore be positive admissible half-integers. If they are integral, a first positive nonintegral half-shift in the actual solution is impossible: its coefficient equation has invertible indicial matrix and zero forcing from previous integral terms. The local-time obstruction then applies. The remaining pairs, including the repeated pair conservatively, are
\[
(1/2,13/2),\quad(3/2,11/2),\quad
(5/2,9/2),\quad(7/2,7/2).
\]
Thus
\[
B\in\left\{\frac{35}{24},\frac58,\frac18,-\frac1{24}\right\}.
\tag{4.52}\label{eq:4.53}
\]
Every retained value differs from \(0\) and \(1\). Equations \eqref{eq:4.48} imply \(\lambda\mu\ne0\).

\needspace{4\baselineskip}
\textbf{A complete finite Darboux calculation.}

Put \(L_0=z-(\mu/\lambda)r\) and rescale \(t,L_0,z\) simultaneously by the nonzero constant \(k\). The highest weighted part of the vector field is
\[
t'=z,\qquad L_0'=t^2,\qquad
z'=\frac{2+B}{3}t^2+\frac{z^2-BL_0^2}{2t}.
\tag{4.53}\label{eq:4.54}
\]
For weights \((2,3,3)\), the highest part \(Q_0\) of a solution of \eqref{eq:4.51} lies in the kernel of
\[
\mathcal L_B=
tz\partial_t+t^3\partial_{L_0}
+\left(\frac{2+B}{3}t^3+\frac{z^2-BL_0^2}{2}\right)\partial_z-2z.
\tag{4.54}\label{eq:4.55}
\]
Coefficients of \(P\) are rational in \(t\). Their Laurent expansions at infinity have a highest weighted term because the fibre degree is finite. Polynomiality in \(t\) is not assumed.

\needspace{3\baselineskip}

\(\square\)

\hypertarget{result-4.9}{}\ResultHeading{Lemma 4.19 (complete weight bound).} A nonzero weighted homogeneous Laurent polynomial in \(\ker\mathcal L_B\), of fibre degree at most seven, has weight \(4\le W\le25\).

\ProofHeading{Proof.} Let \(N\le7\) be its largest fibre degree, with top part \(t^if_N(L_0,z)\). The degree-\(N+1\) equation is
\[
(z^2-BL_0^2)\partial_zf_N+(2i-4)zf_N=0.
\]
Since \(B\ne0\) and \(f_N\) is polynomial, the two distinct factors of \(z^2-BL_0^2\) have equal nonnegative multiplicity \(d=2-i\). Homogeneity gives
\[
f_N=C L_0^{N-2d}(z^2-BL_0^2)^d,\quad
0\le d\le\lfloor N/2\rfloor.
\]
Thus \(W=2i+3N=4+3N-2d\), proving the bounds. \(\square\)

\hypertarget{result-4.10}{}\ResultHeading{Proposition 4.20 (The mixed fifth denominator).} The primitive denominator \(sL\) is impossible for every residual Horikawa basket.

\ProofHeading{Proof.} For each value in \eqref{eq:4.53} and each integer \(4\le W\le25\), use the basis
\[
t^{(W-3(j+\ell))/2}L_0^jz^\ell,\qquad
j+\ell\le7,\quad j+\ell\equiv W\pmod2.
\]
The four possibly nonzero entries in its image are:

\begin{longtable}[]{@{}ll@{}}
\toprule\noalign{}
Target fibre monomial & Coefficient, with \(i=(W-3j-3\ell)/2\) \\
\midrule\noalign{}
\endhead
\bottomrule\noalign{}
\endlastfoot
\(L_0^jz^{\ell+1}\) & \(i+\ell/2-2\) \\
\(L_0^{j-1}z^\ell\) & \(j\) \\
\(L_0^jz^{\ell-1}\) & \(\ell(2+B)/3\) \\
\(L_0^{j+2}z^{\ell-1}\) & \(-\ell B/2\) \\
\end{longtable}

The appropriate power of \(t\) is fixed by target weight \(W+3\). Exact Gaussian elimination gives:

\begin{longtable}[]{@{}lll@{}}
\toprule\noalign{}
Values of \(B\) & Weight \(4\) & Weights \(5,\ldots,25\) \\
\midrule\noalign{}
\endhead
\bottomrule\noalign{}
\endlastfoot
all four in \eqref{eq:4.53} & kernel \(\mathbb Ct^2\) & zero kernel \\
\end{longtable}

This is a finite calculation of 88 rational matrices. The exact matrix algorithm is included in the computational appendix. \hyperlink{result-4.9}{Lemma 4.19} proves that no unexamined weight can occur.

It follows that the highest weighted part of \(P\) is \(ct^2\) with \(c\in\mathbb C\). But \(t^2\) solves the \emph{full} equation \eqref{eq:4.51}, since \(t'=z\). Subtracting \(ct^2\) would leave either zero or a solution of highest weight below four. The latter is impossible by \hyperlink{result-4.9}{Lemma 4.19}. Hence \(P=ct^2\). This makes \(J/t\) constant over \(\mathbb C\), contradicting the nontrivial quadratic extension \(\mathbb C(W)=K(J)\).

\needspace{3\baselineskip}

\(\square\)

\needspace{3\baselineskip}

\hypertarget{result-4.11}{}\ResultHeading{Proposition 4.21 (the denominator \(s^2L\) is impossible).} Assume the exact Euler differential system, the active-divisor and full-field local-time properties, the rank-four seventh restriction, and the invariant fibre-degree bound
\[
\deg_{r,z}(E_9\cap\mathbb C(t,r,z))\le2.
\]
Then a primitive seventh denominator \(F=s^2L\) cannot occur.

\ProofHeading{Proof.} Translate the quartic generator so that \(L=x\). The primitive relation is
\[
s^2xu=G_{13}(s,x,y,z_6).
\]
The only monomial of weight thirteen in variables of weights \(4,5,6\) is \(x^2y\). Since \(s\nmid G_{13}\), its coefficient is nonzero. Weighted division by \(x\), followed by subtracting an element of \(E_6\) from \(h=u/s^7\), gives
\[
h=\eta tr+\gamma t^2+\frac{H(r,z)}t,
\qquad \eta\ne0,\qquad \deg H\le2.
\tag{4.55}\label{eq:4.56}
\]
Write \(H_2\) for the homogeneous quadratic part of \(H\). The restriction rank four implies \(H_2\ne0\): otherwise \(h\) belongs to the \(\mathbb C(t)\)-span of \(1,r,z\), and the seventh restriction has rank at most three. The exact differential system is
\[
\begin{aligned}
t'&=z,\\
r'&=\lambda h+Bz+Cr+Dt+E,\\
z'&=kt^2+\mu h+Gz+Ir+Jt+K.
\end{aligned}
\tag{4.56}\label{eq:4.57}
\]
Its coefficients are complex constants.

The fibre-degree-three part of \(h'\in E_9\) equals
\[
\frac{H_2\bigl((\lambda\partial_r+\mu\partial_z)H_2-z\bigr)}{t^2}.
\]
The degree bound and \(H_2\ne0\) therefore give
\[
(\lambda\partial_r+\mu\partial_z)H_2=z.
\tag{4.57}\label{eq:4.58}
\]
Meanwhile the divergence of \eqref{eq:4.57} is
\[
\eta\lambda t+\frac{\lambda H_r+\mu H_z}{t}+C+G.
\tag{4.58}\label{eq:4.59}
\]
At an active prime \(T\), the first term has pole \(4e\) if nonzero, whereas the second has pole at most \(e+1\). The divergence is the logarithmic derivative \(J'/J\), whose pole is at most \(e+1\). Hence \(\lambda=0\). Equation \eqref{eq:4.58} then implies
\[
\mu\ne0,\qquad
H_2=\frac{z^2}{2\mu}+a r^2.
\tag{4.59}\label{eq:4.60}
\]
No assertion about the linear part of \(H\) is needed.

Set \(R=r-Bt\). Then
\[
R'=CR+dt+E,\qquad d=D+CB.
\tag{4.60}\label{eq:4.61}
\]
If \(d\ne0\), comparison of poles gives \(\operatorname{pole}_T(R)=3e-1\); if \(d=0\), \(R\) is regular. In particular \(\operatorname{pole}_T(R)\le3e-1\). Subtracting another element of \(E_6\) from \(h\) gives
\[
h=\eta tR+\gamma_1t^2+
\frac{z^2/(2\mu)+aR^2+bR+cz+d_0}{t}.
\tag{4.61}\label{eq:4.62}
\]
To check the normalization, expand \(H(R+Bt,z)/t\); its difference from the displayed residue is a linear combination of \(1,t,R,z\). The coefficient \(\gamma_1\) includes the contribution \(\eta B\). The equation for \(z'\) retains the coefficient \(k\) of \(t^2\), with only its lower linear terms changed.

At our reduced detector, \(e=1\) and \(k\ne0\) and \(R\) has pole at most two. The term \(\eta tR\) has pole at most six. Since \(h\) has pole at most seven, its two pole-eight terms in \eqref{eq:4.62} cancel. From \(t\sim(6/k)\tau^{-2}\) and \(z\sim-(12/k)\tau^{-3}\), this gives
\[
\gamma_1=-\frac{k}{3\mu}.
\]
The principal equations are
\[
t'=z,\qquad
z'=\frac{2k}{3}t^2+\frac{z^2}{2t},\qquad R'=dt,
\]
with upper pole weights \((t,R,z)=(2,1,3)\). The indicial polynomial is
\[
(q+1)(q-1)(q-4).
\]
Again there is no positive nonintegral half-shift. The actual functions have no logarithmic terms, so the same local-time obstruction applies. This proves the proposition. \(\square\)

The rank-three case and all nine primitive denominator types have now been excluded. Therefore a quartic fibre of type \((2,3,0)\) cannot occur with \(\chi(\mathcal O_X)=1\) and \(V<1/30\).

\needspace{4\baselineskip}

\hypertarget{exclusion-of-the-220-fibres}{%
\needspace{10\baselineskip}\hypertarget{section-5}{}\section{\texorpdfstring{5. Exclusion of the \((2,2,0)\) fibres}{5. Exclusion of the (2,2,0) fibres}}\label{exclusion-of-the-220-fibres}}

In this section we exclude the remaining case in which the quartic pencil has a \((2,2)\) minimal fibre. By the reductions above, we may assume that
\[
p_g=h^{2,0}=\chi=1,\quad q=0,\quad P_2=P_3=1,\quad P_4=2,
\quad V<\frac1{30},
\]
and that the four-canonical pencil has general minimal fibre
\(S\) with \((K_S^2,p_g(S),q(S))=(2,2,0)\). We use
\(D=Q|_F\), \(A=\sigma^*K_S\), \(a=AD\), and \(d=D^2\).
The preceding restriction results give
\[
D\ge_{\mathbb Q}A/5,\qquad
\frac25\le a,\qquad \frac a5\le d\le4V,\qquad a<\frac23.       \tag{5.1}\label{eq:5.1}
\]
Our argument first bounds the small moving systems on the fibre. We then use their actual canonical markings to exclude the resulting numerical list.

\needspace{4\baselineskip}

\hypertarget{small-systems-and-the-numerical-list}{%
\hypertarget{section-5.1}{}\subsection{5.1. Small systems and the numerical list}\label{small-systems-and-the-numerical-list}}

\needspace{3\baselineskip}

\hypertarget{result-5.1}{}\ResultHeading{Lemma 5.1.} A nonconstant moving system on \(S\) has canonical
degree at least two. A moving system of canonical degree less than
four has at most three sections. If it has a surface image, its
moving divisor \(N\) satisfies
\[
K_SN=N^2=3,\qquad h^0(S,N)=3.                              \tag{5.2}\label{eq:5.2}
\]
Its resolved map has degree two or three onto \(\mathbb P^2\).

\ProofHeading{Proof.} We first show that \(\operatorname{Pic}(S)\) has no nontrivial torsion.
A torsion bundle of order \(n\ge2\) would give a connected étale cover
with \(K^2=2n,\ \chi=3n,\ p_g\ge3n-1\), contrary to Noether's
inequality. Since \(q(S)=0\), numerically equivalent line bundles
are therefore linearly equivalent.

We now consider a moving divisor. It is nef, and Hodge index together with adjunction parity excludes
canonical degree one: its square would have to be zero, with odd
canonical intersection. If the system has a surface image and canonical degree two, these inequalities force square two and numerical, hence
linear, equivalence to \(K_S\), contradicting \(p_g(S)=2\).
For canonical degree three, its square is one or three. Square one
would give a birational plane, impossible for a surface of general
type. Hence its square is three. Its resolved square is two or
three, since square one again gives a birational plane. In particular, there is at most one simple base point. A general resolved member has genus
four, and the restricted moving system has degree at most three.
Clifford's theorem bounds its dimension by two, so the original
system has at most three sections.

A system composed with a pencil and having three sections contains
at least twice a connected moving member. Its base is rational
since \(q(S)=0\), and that member has canonical degree at least two.
Thus its canonical degree is at least four. \(\square\)

\hypertarget{prop-small-net-volumes}{}\ResultHeading{Proposition 5.2 (Small-net volume estimates).} In the \((2,2,0)\)-fibre case, \(V>1/40\). If \((P_5,P_6,P_7)=(3,4,5)\), then \(V\ge1/35\). The only remaining baskets are those in \eqref{eq:5.6}.

\ProofHeading{Proof.} Let \(\mathcal V_m\) be the image of global evaluation along the
outer pencil and \(\rho_m\) its rank. This is a globally generated
bundle on \(\mathbb P^1\), with exactly \(P_m\) sections:
\[
\mathcal V_m=\bigoplus_{i=1}^{\rho_m}\mathcal O(a_{mi}),
\qquad a_{mi}\ge0,\qquad\sum_i a_{mi}=P_m-\rho_m.            \tag{5.3}\label{eq:5.3}
\]
Multiplication by \(s\) makes the ranks nondecreasing. A summand of
degree \(k\) gives \(Q\ge_{\mathbb Q}(k/m)F\), and hence, writing
\(\mu=k/m\),
\[
V\ge\frac{2\mu^3}{(1+\mu)^2}.                              \tag{5.4}\label{eq:5.4}
\]
If two summands are positive, sections vanishing at one prescribed
outer point give a moving divisor \(G\) with
\(mQ\ge_{\mathbb Q}F+G\) and \(A(G|_F)\ge2\). Resolve its base
points, so that \(G\) is nef. Then
\[
mV\ge d+Q^2G\ge\frac a5+\frac1{10}.                        \tag{5.5}\label{eq:5.5}
\]
Indeed \(4Q-F\) is effective and
\(4Q^2G\ge QFG\ge A(G|_F)/5\).

We next combine these inequalities with the basket list to obtain \(V>1/40\). Since the endpoint will matter, we give the argument explicitly. Suppose that
\(V\le1/40\), then \(a\le1/2\) and the sixth restriction has rank
two or three. Rank three forces \(a=1/2,d=1/10,V=1/40\),
the impossible \(K_S^2=2\) endpoint of \cite[Proposition 4.3]{ref1}. In rank two, \(P_6=4\) gives splitting \((2,0)\) or
\((1,1)\), implying respectively \(V\ge1/24\) or \(V\ge3/100\).
Every row with \(P_6=3\) has \(P_5=3\), contrary to Corollary~2.2. Thus \(V>1/40\).

If \((P_5,P_6,P_7)=(3,4,5)\), one furthermore has \(V\ge1/35\).
Otherwise \(\rho_6=\rho_7=3\), \(a\ge1/2\), and the seventh
splitting is \((2,0,0)\) or \((1,1,0)\). The latter gives
\(7V\ge1/5\). The former gives \(\mu=2/7\), \(D\ge_{\mathbb Q}2A/9\);
the sixth net has \(AN=3\), so
\(d\ge DN/6\ge1/9\) and \(V\ge2/63\).

Thus \eqref{eq:2.1} and the packing enumeration give the following necessary list. We keep rows of the same volume distinct, since their higher plurigenera need not agree.
\[
\begin{array}{c|l|c}
V&B&(P_5,P_6,P_7,P_8,P_9,P_{10})\\ \hline
2/77&8(1,2)+(2,7)+(4,11)&(2,4,4,6,6,9)\\
25/858&3(1,2)+(1,3)+(5,11)+(5,13)&(3,4,5,7,8,10)\\
7/240&3(1,2)+2(1,3)+(2,5)+(7,16)&(3,4,5,6,8,9)\\
5/168&8(1,2)+(1,3)+(2,7)+(3,8)&(2,4,4,6,6,9)\\
5/168&4(1,2)+(1,3)+2(3,7)+(3,8)&(3,4,5,6,7,9)\\
19/630&3(1,2)+2(1,3)+(2,5)+(3,7)+(4,9)&(3,4,5,6,8,9)\\
2/65&8(1,2)+(2,5)+(4,13)&(2,4,4,6,7,10)\\
41/1320&3(1,2)+(1,3)+(2,5)+(3,8)+(5,11)&(3,4,5,7,8,10)\\
8/255&2(1,3)+2(2,5)+(8,17)&(3,4,5,7,9,11)\\
16/495&4(1,2)+(2,5)+(4,9)+(4,11)&(3,4,5,7,8,10).
\end{array}                                                     \tag{5.6}\label{eq:5.6}
\]

\needspace{3\baselineskip}

\(\square\)

\hypertarget{result-5.2}{}\ResultHeading{Lemma 5.3.} The sixth map is generically finite of degree two
or three onto \(\mathbb P^3\).

\ProofHeading{Proof.} By \hyperlink{result-5.1}{Lemma 5.1}, we have \(\rho_6=2\) or \(3\). We first exclude rank two. In that case,
\eqref{eq:5.4} excludes splitting \((2,0)\), so the splitting is
\(\mathcal O(1)^{\oplus2}\) and the image is a smooth quadric.
We claim that the inner canonical pencil is free of genus two. If instead it had genus three, the curve estimate would give
\(\xi=DC\ge1/2\) by \cite[Lemma 2.4]{ref1} applied to
\(m\xi\ge4+\lceil(m-10)\xi\rceil\).
Then the two rulings would give
\(6V\ge d+\xi/4\ge1/10+1/8>1/5\).

We may therefore work in the genus-two case and choose ruling coordinates
\(E_6=\langle1,t,u,tu\rangle\). Write
\(t'=A(t)+B(t)u\) with affine \(A,B\); Hodge independence gives
\(B\ne0\). Since \(u', (tu)'\in E_8\), the identity
\[
(tu)'-tu'-A(t)u=B(t)u^2
\]
puts \(1,u,u^2\) in \(\mathbb C(t)E_8\). On the general surface
fibre this subsystem has moving divisor \(2C\), so \(8a\ge4\).
The ruling estimate now gives
\(6V\ge a/5+\xi/4\ge1/10+1/10=1/5\), since the curve
restriction inequality gives \(\xi\ge2/5\).
This contradicts \(V<1/30\). Thus \(\rho_6=3\); \hyperlink{result-5.1}{Lemma 5.1}
gives the asserted degree. \(\square\)

\needspace{4\baselineskip}

\needspace{3\baselineskip}
\hypertarget{prop-five-over-168}{}\ResultHeading{Proposition 5.4.}
The basket
\[
 B_*=4(1,2)+(1,3)+2(3,7)+(3,8),\qquad V=5/168,
\]
does not occur in the \((2,2,0)\)-fibre case.

\ProofHeading{Proof.}
Here \((P_1,P_5,P_6,P_{10})=(1,3,4,9)\), and \(P_2(F)=5\).
Since \(q(F)=0\), we have \(\rho=h=1\). If the full general quartic
member on \(Z\) were reduced, \eqref{eq:reduced-member-trace} would give
\(5=P_{10}-P_6\ge5+(3-1-2+1)=6\).
It therefore has a fixed prime of multiplicity at least two. Its strict
transform has positive canonical degree, and \(168K_Y\) is Cartier;
\cite[Lemma 3.12]{ref1} gives \(K_Y^2Z_4\ge2/168\) for the quartic
fixed part \(Z_4\).

Choose the special pencil model of \cite[Proposition 5.6]{refM}, and put
\(L=\mu^*K_Y|_F\), \(A=\sigma^*K_S\), and
\(d=L^2=K_Y^2F_Y\). The comparison there gives
\[
 d\le4V-\frac2{168}=\frac3{28},\qquad
 LA\le d+\frac{d^2}{V}\le\frac{69}{140}<\frac12.
\]
On a common refinement resolving the sixth system, let \(N\) be its
moving restriction. Then \(6\widehat L\ge_{\mathbb Q}N\), and
\hyperlink{result-5.1}{Lemma 5.1} together with
\hyperlink{result-5.2}{Lemma 5.3} gives \(\widehat A N=3\).
Thus \(LA\ge1/2\), a contradiction. \(\square\)

\hypertarget{the-actual-small-net-marking}{%
\hypertarget{section-5.2}{}\subsection{5.2. The actual small-net marking}\label{the-actual-small-net-marking}}

\needspace{3\baselineskip}

\hypertarget{result-5.3}{}\ResultHeading{Lemma 5.5.} For the sixth net \(N\) there are unique effective
divisors \(E,U\) with
\[
N\sim K_S+E,\qquad K_S\sim E+U,\qquad
K_SE=K_SU=1,\quad E^2=U^2=-1,\quad EU=2,\quad NE=0.
                                                               \tag{5.7}\label{eq:5.7}
\]
If \(A_0\) is the actual canonical divisor induced by \(s\), and
\(D_S=\sigma_*D\), then the actual constant sixth member is
\(N_0=A_0+E\), and
\[
A_0=E+U,\qquad 0\le D_S\le A_0,\qquad 6D_S\ge2E+U.        \tag{5.8}\label{eq:5.8}
\]

\ProofHeading{Proof.} We set \(E=N-K_S\) and \(U=2K_S-N\). The displayed
intersection numbers follow from \eqref{eq:5.2}. Neither has two sections,
since a moving divisor cannot have canonical degree one.
Riemann--Roch gives
\(h^0(E)+h^0(U)=2+h^1(E)\), so both have exactly one section.
Nefness of \(N\) and \(NE=0\) make each component of \(E\)
\(N\)-orthogonal.

To identify the actual divisors, we write \(s=dt\wedge\beta\) and
\(\alpha|_F=t'\beta\). These are independent canonical forms.
Multiplication by the section of \(E\) embeds them into the
three-dimensional net. Let \(L_\alpha,L_\beta\) be their independent
linear coordinates, \(x_0\) the constant sixth member, and \(L_r\)
the member representing \(t'\in E_6\). Dominance onto the plane
gives \(L_\alpha x_0=L_\beta L_r\); coprimeness implies
\(x_0\) is proportional to \(L_\beta\). Hence \(N_0=A_0+E\).
Actual moving-part comparison gives \(6D_S\ge A_0+E\), and
discrepancies give \(D_S\le A_0\).

It remains to show that \(E\le A_0\). We give a cycle argument so that nonreduced components are included.
It has a unique positive canonical-degree component \(L\), of
degree and multiplicity one. The inequality \(E\le5A_0\)
puts \(L\) in \(A_0\). Starting with \(L\), add a \(K_S\)-null
component in \(\operatorname{Supp}(E)\) whenever its intersection
with the current cycle is positive. The process stays below \(E\), because \(E\) has zero
intersection with each such component, and below \(A_0\), because
a canonical divisor has zero intersection with it. It terminates
at a cycle \(\Gamma\). If \(\Gamma<E\), the difference is supported
on a negative definite configuration of \((-2)\)-curves and
\(\Gamma(E-\Gamma)\le0\); hence \(E^2<\Gamma^2\).
But \(K_S\Gamma=1\), Hodge index and parity give
\(\Gamma^2\le-1\), contradicting \(E^2=-1\).
Thus \(E\le A_0\). Uniqueness identifies \(A_0-E\) with \(U\).
This proves \eqref{eq:5.8}. \(\square\)

The positive canonical-degree components of \(E,U\) are either
distinct with canonical multiplicity one, or coincide with
canonical multiplicity two. Their coefficients in \(D\) are
strictly between zero and one, by \eqref{eq:5.1} and \eqref{eq:5.8}. Each containing
horizontal divisor therefore lies over one non-Gorenstein point
of \(Y\): a nonexceptional divisor, a divisor with curve centre,
or one over a Gorenstein point would have integral discrepancy.
Thus these curves are \(D\)-null and their coefficients are
cleared by one actual local index.

\needspace{3\baselineskip}

\hypertarget{result-5.4}{}\ResultHeading{Lemma 5.6 (adjunction at a marked curve).} Let \(F\) be a smooth projective surface, let \(D\) be a nef and big \(\mathbb Q\)-divisor, and let \(A\) be a divisor such that \(K_F-A\ge0\). Suppose that \(\lambda D-K_F\ge0\). Let \(C\) be an integral \(D\)-null curve with \(AC=b>0\), and assume that \(C\) is not a component of \(K_F-A\). Write
\(c=\operatorname{coeff}_C D\),
\(m=\operatorname{coeff}_C K_F\). Then
\[
\lambda c-m\ge\frac b{b+2}.                               \tag{5.9}\label{eq:5.9}
\]

\ProofHeading{Proof.} We put \(k=-C^2\) and \(h=K_FC\). Hodge index gives
\(k>0\), and \(K_F-A\ge0\), with no component \(C\), gives
\(h\ge b\). Removing the known coefficient
\(\gamma=\lambda c-m\) of \(C\) from the effective divisor
\(\lambda D-K_F\) gives \(0\le-h+\gamma k\).
Adjunction yields \(k\le h+2\). Thus
\(\gamma\ge h/k\ge h/(h+2)\ge b/(b+2)\). \(\square\)

\needspace{3\baselineskip}

\hypertarget{result-5.5}{}\ResultHeading{Corollary 5.7.} Under the assumptions of this section, no rational seventh section satisfies \(sW=x^2\).

\ProofHeading{Proof.} This is Corollary~2.2. \(\square\)

\hypertarget{two-plurigenus-plateaux}{%
\hypertarget{section-5.3}{}\subsection{5.3. Two plurigenus plateaux}\label{two-plurigenus-plateaux}}

We now isolate the two arguments for \((P_5,P_6,P_7)=(2,4,4)\). Their proofs use only the differential identities and the prime degree of the sixth map. We can therefore apply them again to the \((1,2,0)\) fibres in \hyperlink{section-6}{Section 6}.

\needspace{3\baselineskip}

\hypertarget{result-5.6}{}\ResultHeading{Lemma 5.8 (the two plateaux).} Assume the hypotheses of the
canonical-quotient proposition, \(t'=v(x/s^4)\ne0\), and
\[
(P_1,\ldots,P_7)=(1,1,1,2,2,4,4).
\]
Suppose \(L_6\) is rational of degree two or three in
\(\mathbb C(W)\). Neither of the following configurations occurs:

\begin{enumerate}
\def\labelenumi{(\alph{enumi})}
\item
  \(P_8=P_9=6\);
\item
  the field \(\mathbb C(t,t')\) is the base field of a connected fibration in curves of genus two or three, \(V=2/65\), and
  \((P_8,P_9,P_{10},P_{11},P_{12})=(6,7,10,10,14)\).
\end{enumerate}

\ProofHeading{Proof of (a).} Apply \cite[Lemma 3.19]{ref1}.

\ProofHeading{Proof of (b).} We begin with the coordinates supplied by \cite[Proposition 3.17]{ref1} and \cite[Corollary 3.18]{ref1}, write \(\theta^2=\lambda\bar x\ne0\), \(k=24\lambda\), \(r=t'\) and \(z=t''-kt^2\in E_6\). At a reduced transverse component where \(x\) is a unit, the actual time is \(\tau=\rho^2\), and
\[t=(6/k)\tau^{-2}+\cdots,\qquad t^2\notin E_6.\tag{5.10}\label{eq:5.10}\]
The later plateau
\(R_{11}=sR_{10}\) gives the exact sections
\[
\mathcal E=r^2-\frac{2k}{3}t^3\in E_{10},\qquad
xR_9/s\subset R_{12}.                                   \tag{5.11}\label{eq:5.11}
\]
Indeed \(\bar r_6^2=(2k/3)\bar x^3\), so
\[
r_6^2-(2k/3)x^3\in sR_{11}=s^2R_{10}.
\] For \(g_9\in R_9\), its derivative belongs to
\(R_{11}=sR_{10}\), hence \(\theta\bar g_9=0\); multiplying by
\(\theta\) gives \(\bar x\bar g_9=0\).

First suppose \(1,t,r,z\) are independent. Set
\[
w=z',\quad g=z''-2ktz,\quad
\mathcal P=\langle1,t,r,z,t^2,tr,tz,\mathcal E\rangle .
\]
The space \(\mathcal P\) has dimension eight, \(w\in E_8\),
\(g\in E_9\), and \(g'\in E_{10}\).
If \(\mathcal P,w,g\) are independent, they form \(E_{10}\).
The equation for \(g'\) has fifth-order indicial polynomial
\((q+1)(q-1)(q-6)^2(q-8)\), so reconstruction gives a contradiction.

We next show that \(w\notin\mathcal P\). Suppose to the contrary that
\[
w=c\mathcal E+dtr+etz+L,\qquad
L\in\langle1,t,r,z,t^2\rangle.
\]
For \(c=e=0\), the only pole-five term forces \(d=0\),
leaving the integral-time third-order equation already treated.
For \((c,e)\ne(0,0)\), the highest-pole part has coefficient
\[
\frac{q+1}{k}\bigl(6e(q-6)-24c\bigr)
\]
at the next shift. It has at most one positive root. Every
coefficient of the existing Puiseux series is consequently
determined by constants and at most one free transverse
coefficient. The negative root only moves the pole, whose
location is fixed by \(\tau=0\). Two independent transverse
variations are therefore impossible, whereas
\(dt\wedge dr\wedge dz\ne0\).

It follows that the remaining dependence has the form
\[
g=fw+c\mathcal E+dtr+etz+L.                              \tag{5.12}\label{eq:5.12}
\]
We claim \(w\in K=\mathbb C(t,r,z)\).
To prove this claim, we apply the Euler identity to \eqref{eq:5.12}, using the division in \eqref{eq:5.11}, and then differentiate in the usual way. We obtain, respectively,
\[
R=(4c+5e)rz+5dr^2\in E_{10},                             \tag{5.13}\label{eq:5.13}
\]
\[
H=(2c+e)rz+etw+\frac{5kd}{3}t^3\in E_{10}
       \pmod{\mathcal P+\mathbb Cw}.                    \tag{5.14}\label{eq:5.14}
\]
For a direct check of the coefficients, we may differentiate
\(\mathcal E\): \(\mathcal E'=2rz\); the final Euler division
uses precisely \(R_{11}=sR_{10}\).

Suppose \(w\notin K\). If \(R\notin\mathcal P\), it completes
\(\mathcal P+\mathbb Cw\) to \(E_{10}\). Comparing the coefficient
of \(w\) in \eqref{eq:5.14} forces \(e=0\).
If \(c\ne0\), the section \(tg\in E_{12}\), supplied by
\eqref{eq:5.11}, gives
\(T=ct\mathcal E+dt^2r\) modulo
\(B_{12}^0+\langle w,tw\rangle\), where
\(B=\mathbb C[s,x,r_6,z_6]\).
Its \(t^4\) term makes it independent of the eleven-dimensional
\(B_{12}^0\). Since \(P_{12}=14\), these spaces exhaust \(E_{12}\).
But \(R'\) has coefficient \(4cr\) along \(w\), not in
\(\langle1,t\rangle\), a contradiction.
If \(c=0\), the pole-five term \(dtr\) in \eqref{eq:5.12} forces \(d=0\):
\(g'\in E_{10}\) implies that a pole of \(g\) has order at most
four. This again contradicts \(R\notin\mathcal P\).

If \(R\in\mathcal P\), polynomial independence in \eqref{eq:5.13} gives
\(d=0,c=-5e/4\). For \(e=0\), \eqref{eq:5.12} is the integral-time
fourth-order equation in the proof of \cite[Lemma 3.19]{ref1}. For \(e\ne0\),
\eqref{eq:5.14} gives \(tw-\tfrac32rz\in E_{10}\).
The section \(tg\) supplies
\(T=e(t^2z-\tfrac54t\mathcal E)\), again independent of
\(B_{12}^0\), and the same dimension exhaustion holds.
The derivative of \(tw-\tfrac32rz\) has coefficient
\(ft-r/2\) along \(w\), again outside \(\langle1,t\rangle\).
This proves \(w\in K\).

If \([L:K]=3\), write its actual eighth section primitively as
\(G/F\) over \(B\), with denominator degree \(d_F\).
The hypersurface subring gives
\[
\frac2{65}\ge\frac{3(8+d_F)}{1152},\qquad d_F\le3.
\]
The denominator is a power of \(s\). Odd degrees would make
numerator and denominator divisible by \(s\), so only degrees
zero and two occur. Their normalized monomials lie in
\(\langle1,t,r,z,t^2,tr,tz\rangle\subset\mathcal P\), a
contradiction. Thus the sixth degree is two in this case.

If instead \(z\in\langle1,t,r\rangle\), choose \(u\) completing
\(E_6\), and set \(w=u'\), \(g=u''-2ktu\),
\(\mathcal P=\langle1,t,r,u,t^2,tr,tu,\mathcal E\rangle\).
If \(\mathcal P,w,g\) are independent, the equation for \(g'\)
is linear in \(u\), with principal exponents \(-3,2,4\).
If \(w\in\mathcal P\), its first-order equation excludes
fractional powers, whether or not its coefficient of \(u\)
has a double pole. These alternatives contradict reconstruction.
The remaining dependence is
\[
u'=w,\qquad w'=fw+(a_1t+a_0)u+B(t,r).                        \tag{5.15}\label{eq:5.15}
\]
The connected curve over \(\mathbb C(t,r)\) is the canonical
curve, of genus two or three. The function \(u\) has degree two
or three on it. Its pole divisor bounds that of \(w\):
at a point where differentiation increases a pole by \(d>0\),
the successive poles of \(u,w,w'\) would be \(n,n+d,n+2d\),
incompatible with \eqref{eq:5.15}. A pole of \(w\) at a regular point of
\(u\) is excluded in the same way. This computation is made
after completing and extending the residue field; the derivation
extends uniquely over the finite separable coefficient field.
Clifford's theorem and Riemann--Roch give
\(h^0(\operatorname{div}_\infty u)=2\), so
\(w\in\mathbb C(t,r)+\mathbb C(t,r)u\subset K\).
If the sixth degree were three, the same denominator bound would
put \(w\in\mathcal P\), already excluded. Thus the degree is two
also in this case.

We can now treat the two cases simultaneously. Write
\(B=\mathbb C[s,x,a_6,b_6]\), and let \(w_8=G/F\) be the
actual eighth derivative, rational over \(B\), with
\(w\notin\mathcal P\). The field-index bound gives
\[
d_F\le9.
\]
Parity and primitivity exclude odd \(d_F\); degrees zero and two
put \(w\) in \(\mathcal P\). Thus \(d_F=4,6,8\).
If \(F\) involves a sixth coordinate, then
\(F=s^{d_F-6}L_6+h(s,x)\), where \(L_6\ne0\).
Replacing one sixth coordinate by \(L_6\), use the polynomial
ring
\[
B'=\mathbb C[s,x,b_6,F_{d_F}].
\]
After clearing powers of \(s\), the given fraction supplies an
actual section with equation \(FU=g(s,x,b_6)\ne0\).
The hypersurface \(FU=g\) is normal canonical; successive
general slices through its singular locus have equation
\(FU=z^q\), and inversion of adjunction applies.
Its homogeneous canonical form has weight
\(1+4+6=11\). The canonical-cone pullback construction therefore
gives a nonzero odd canonical section of degree twelve.
This is the proof of \cite[Lemma 3.14(c)]{ref1} applied to \(B'\); its degree is independent of the
power of \(s\) used in the coordinate change.

Let \(A=B[w_8]\). For \(d_F=6,8\), its degree-twelve part has
dimension thirteen, while \(P_{12}=14\). Hence the invariant
part of \(R_{12}\) is exactly \(A_{12}\), and the odd part has
dimension one. Every ninth section is invariant: an odd \(g_9\)
would give the two independent odd twelfth sections
\(s^3g_9,xg_9/s\), by \eqref{eq:5.11}.
Choose \(g_9\notin sR_8\). Represent these two twelfth sections
by polynomials \(p_{12},q_{12}\) in variables of weights
\((1,4,6,6,8)\). Their relation is
\[
xp_{12}-s^4q_{12}
 =c\,s^{8-d_F}(F_{d_F}W-G_{8+d_F}).
\]
The coefficient of \(W\) on the left lies in
\(\langle x^2,s^4x,s^8\rangle\), whereas the right contains
\(cs^2L_6\). Thus \(c=0\). Coprimeness gives \(s^4\mid p_{12}\),
forcing \(g_9\in sR_8\), a contradiction.
Consequently \(F\in\mathbb C[s,x]\), and
\[
w\in\mathbb C(t)[a,b],\qquad \deg_{a,b}w\le2.             \tag{5.16}\label{eq:5.16}
\]

In the independent case take \(a=r,b=z\), so
\(v=r\partial_t+(kt^2+z)\partial_r+w\partial_z\).
Write
\[
w=A(t)z^2+B(t)rz+C(t)r^2+D(t)z+E(t)r+F(t).
\]
Equation \eqref{eq:5.12} becomes
\[
v(w)=fw+c\mathcal E+dtr+(2k+e)tz+L.
\]
Comparing \(z^3,z^2,r^3,rz,r^2\) gives
\[
A=B=0,\quad \dot C=0,\quad \dot D=-2C,\quad
\dot E=fC+c-CD,
\]
where dots denote differentiation in \(t\).
The \(z\)-coefficient says
\(E+D^2=fD+(2k+e)t+\text{constant}\).
Its \(t^2\)-coefficient is \(5C^2=0\).
Thus \(D\) is constant, \(E\) affine, and the \(r\)-coefficient
makes \(F\) quadratic. This puts \(w\in\mathcal P\), impossible.

In the dependent case
\(t'=r,\ r'=kt^2+c_0r+bt+a\). Write
\[
w=A(t)u^2+B(t)ru+C(t)r^2+D(t)u+E(t)r+F(t).
\]
Equation \eqref{eq:5.15} gives \(A=0\) and
\[
\dot B+B^2=0,\qquad
\dot D+2BD=(f-c_0)B.
\]
Their rational solutions are
\[
(B,D)=(0,\text{constant}),\quad
\left(\frac1{t-a_*},\frac{f-c_0}{2}+
                    \frac{d_0}{(t-a_*)^2}\right).
\]
The scalar equation for \(t\) gives integral actual-time
expansions for \(t,r\) at \eqref{eq:5.10}. The resulting first-order
equation for \(u\) has integral Laurent forcing and coefficient
of residue zero or \(-2\): in the second case
\(r/(t-a_*)=-2/\tau+\) an integral holomorphic series.
Its homogeneous exponent is integral, so a first fractional
exponent of \(u\) is impossible. Reconstruction gives the final
contradiction. This proves (b). \(\square\)

Thus \hyperlink{result-5.6}{Lemma 5.8} removes all three rows with \((P_5,P_6,P_7)=(2,4,4)\)
from \eqref{eq:5.6}.

\needspace{4\baselineskip}

\hypertarget{the-seventh-field-and-the-denominator-reductions}{%
\hypertarget{section-5.4}{}\subsection{5.4. The seventh field and the denominator reductions}\label{the-seventh-field-and-the-denominator-reductions}}

We turn to the six remaining rows, all of which satisfy
\((P_5,P_6,P_7)=(3,4,5)\).
The fifth restriction has rank two, since rank one would give
a degree-two outer pencil in degree five, contradicting \eqref{eq:5.4}.
Choose
\[
B=\mathbb C[s,x_4,y_5,z_6],\qquad
K=\mathbb C(t,u,z),\qquad
t=x/s^4,\ u=y/s^5,\ z=z_6/s^6.
\]
Then \(B\) is polynomial, and the field index is two or three.

We first record the saturation argument that will control the seventh field. If \(B\) is a unique
factorization domain, \(W\notin\operatorname{Frac}B\), and
coprime \(A,C\in B\) satisfy \(Av,Cv\in B+BW\), comparison
of the coefficients of \(1,W\) gives \(v\in B+BW\).
For a new seventh section \(W\notin B_7\), the dimensions of
\((B+BW)_m\), if \(W\) is irrational, are
\[
\begin{array}{c|rrrr}
m&8&9&12&13\\ \hline
\dim(B+BW)_m&6&7&14&16.
\end{array}                                                    \tag{5.17}\label{eq:5.17}
\]
For the rows \(25/858,41/1320,16/495\), an extra eighth section
and \(P_{12}=14\), multiplied by \(s^4,x\), contradict this
saturation statement. For \(7/240,19/630\), use an extra ninth
section and \(P_{13}=16\). Thus \(L_7=K\) in those five rows.

In the five rows just considered, write a rational new seventh section primitively as \(W=G/F\),
with \(d_F=\deg F>0\). The hypersurface subring gives
\[
V\ge\frac{\delta(7+d_F)}{840},
\qquad \delta=[\mathbb C(W):K]\in\{2,3\}.                 \tag{5.18}\label{eq:5.18}
\]
For \(\delta=3\), this forces \(d_F\le2\).
Degree one is impossible by \hyperlink{result-5.5}{Corollary 5.7}, and degree two by
\hyperlink{result-6.2}{Lemma 6.3}. Thus the five rows now under consideration have
quadratic sixth field. Its involution negates \(s,\alpha\),
so \(v(K)\subset K\), by \cite[Corollary 3.9]{ref1}.

We collect the numerical data needed by the denominator lemmas below. Stating them together will allow us to use the same algebraic argument for both fibre types.
\[
\begin{array}{c|c|c|c}
V&P_8&P_{12}&\text{additional data}\\ \hline
25/858&7&14&P_9=8,\ P_{10}=10,\ P_{13}=17,\ P_{15}=23\\
7/240&6&14&P_{10}=9,\ P_{11}=11,\ P_{13}=16\\
19/630&6&14&P_{10}=9,\ P_{11}=11,\ P_{13}=16\\
41/1320&7&14&P_9=8,\ P_{10}=10,\ P_{13}=17,\ P_{15}=24\\
16/495&7&14&P_9=8,\ P_{10}=10,\ P_{13}=17,\ P_{15}=24 .
\end{array}                                                    \tag{5.19}\label{eq:5.19}
\]
All have \(d_F\le6\) by \eqref{eq:5.18}.
Degree one is excluded by \hyperlink{result-5.5}{Corollary 5.7}.
Degree two is excluded by \hyperlink{result-6.2}{Lemma 6.3}.
For degree five, use \hyperlink{result-6.4}{Lemma 6.5} when \(P_8=6\), and
\hyperlink{result-6.3}{Lemma 6.4} when \(P_8=7\). Degree six gives
\(V\ge13/420\) by \eqref{eq:5.18}, which already excludes the two
rows with \(P_8=6\); for the others use \hyperlink{result-6.5}{Lemma 6.6}.
When \(P_8=7\) and \(d_F\ge2\),
the required eighth-field equality follows immediately:
the invariant subring has at least thirteen sections in degree
twelve; an odd eighth section multiplied by the two-dimensional
\(R_4\) would give \(P_{12}\ge15\).
For \(d_F=3,4\) and \(V\le8/255\), use \hyperlink{result-7.1}{Theorem~7.1}(b).
At \(16/495\), use Lemmas 6.8 and 6.9 with the last line of \eqref{eq:5.19}.
Thus these five rows are excluded. The only remaining numerical
row is \(V=8/255\); no equality \(L_7=K\) is assumed for it.

\needspace{4\baselineskip}

\hypertarget{the-final-basket-and-its-actual-markings}{%
\hypertarget{section-5.5}{}\subsection{5.5. The final basket and its actual markings}\label{the-final-basket-and-its-actual-markings}}

Let
\[
B=2(1,3)+2(2,5)+(8,17),\qquad V=\frac8{255}.
\]
The actual non-Gorenstein indices are among \(3,5,17\).
This follows from the terminal classification: a local basket
has the actual local index, except for the \(cAx/4\) case,
whose basket includes index four \cite[Section 2]{ref17}.
The marking of \hyperlink{result-5.3}{Lemma 5.5}, valid independently of \(L_7\),
and \eqref{eq:5.1} give
\[
\frac12\le a\le\frac{32}{51}.
\]
If the positive components of \(E,U\) are distinct, their
coefficients \(b,c\) satisfy \(b\ge1/3\) by \eqref{eq:5.8} and
\(c\ge4/15\) by \eqref{eq:5.9} with \(\lambda=5,m=1\).
Each is cleared by one of the actual indices \(3,5,17\).
The only possibility with \(b+c\le32/51\) is
\[
(b,c)=(1/3,5/17).
\]
If they coincide, their coefficient is \(a\), and the only
possibilities in the stated interval are
\[
a=9/17,\qquad 10/17,\qquad 3/5.                          \tag{5.20}\label{eq:5.20}
\]

To treat the index-seventeen markings, we use the following valuation calculation. It avoids a classification of the marked divisorial valuation.

\needspace{3\baselineskip}

\hypertarget{result-5.8}{}\ResultHeading{Lemma 5.9 (valuation characters).} Consider a smooth cyclic quotient of index \(r\ge2\), with nonzero coordinate characters
\(w_1,w_2,w_3\). Put \(b=w_1+w_2+w_3\), and assume that \(\gcd(b,r)=1\).
Suppose an invariant canonical section \(s=f\nu\) has pullback
coefficient \(k/r\), where \(0<k<r\) and \(\gcd(k,r)=1\), at a divisorial valuation centred
at the quotient point. Set \(j=k(-b)^{-1}\bmod r\).
For its normalized extension \(\widetilde v\) to the cover,
\[
\widetilde v(f)=k,\qquad
\widetilde v(X_i)\ge\langle jw_i\rangle_r.
\]
If \(\eta_m=g\nu^m\) is invariant, then
\[
\widetilde v(g)\equiv mk\pmod r,\qquad
v(\eta_m/s^m)=\frac{\widetilde v(g)-mk}{r}.                \tag{5.21}\label{eq:5.21}
\]
For nonzero character of \(g\), its value is at least the
minimum of the three positive coordinate bounds.

\ProofHeading{Proof.} If the ramification index is \(e\), then
\(\widetilde v(f)=ek/r\) is integral. Since \(\gcd(k,r)=1\) and the cover has degree \(r\), we have \(e=r\); the residue extension is
trivial and the entire cyclic group is inertia.
Its character \(h\) on the cotangent line of the valuation
ring is faithful. Indeed a finite-order automorphism in
characteristic zero acting trivially both on the residue field
and that line is the identity, as follows successively from its
first nonzero higher term.
A semi-invariant of character \(c\) and value \(n\) therefore
has \(hn=c\bmod r\). Apply this to \(f\), of character \(-b\),
and then to the coordinates and to \(g\).
A nonzero-character coefficient belongs to the maximal ideal,
giving the last bound. \(\square\)

There is exactly one index-seventeen basket constituent.
The terminal point is \(cA/17\), whose basket has its axial
weight many identical index-seventeen constituents \cite[Remark 2.1]{refCH11}. Axial weight one makes its index-one cover smooth:
in \(XY+g(Z^{17},U)=0\) the coefficient of \(U\) is nonzero.
Thus the actual quotient is
\(\frac1{17}(1,16,8)\), and the volume form has character eight.
For coefficient \(5/17\), \hyperlink{result-5.8}{Lemma 5.9} gives \(j=10\) and
coordinate values at least \(10,7,12\), whereas \(f\) would
have value five. Since \(f\) lies in the maximal ideal this
is impossible. The distinct-component marking is excluded.

For either shared index-seventeen marking, write
\(s=dt\wedge\beta\), \(r=t'=\alpha|_F/\beta\).
A positive canonical-degree curve cannot be fixed in \(|K_S|\):
otherwise its moving part would have canonical degree one,
contrary to \hyperlink{result-5.1}{Lemma 5.1}. Since \(\beta\) vanishes twice along
the common positive component, \(r\) has exact pole two there.
Every function in \(E_5\) has pole at most
\(\lfloor5a\rfloor=2\).
The constant sixth member \(A_0+E\) has multiplicity three
there, and the sixth net has no fixed curve; hence its unique
new direction modulo \(E_5\) has exact pole three.
Thus
\[
r\in E_5,\qquad E_5=\langle1,t,r\rangle .
\]
This contradicts Corollary~2.2.

It remains to treat \(a=3/5\). Here we use the fixed divisor of the outer pencil.
If the actual divisors of \(s\) and \(x\) on the canonical model
had no common prime, multiplication by \(x\) would inject
\[
R_n/sR_{n-1}\longrightarrow R_{n+4}/sR_{n+3}.
\]
Indeed \(xa=sb\) would force \(a/s\) to have no divisorial pole,
hence to be a section by reflexivity. But this basket has
\[
P_{12}-P_{11}=3>2=P_{16}-P_{15}.
\]
Thus the quartic system on the canonical model has a fixed
prime. Compute intersections on a crepant \(\mathbb Q\)-factorial
terminal model \(Y\), and write \(4K_Y\sim T+Z_4\), where
\(T\) is a general pencil member. The strict transform of that
fixed prime has positive \(K_Y^2\)-intersection by ampleness on
the canonical model; exceptional divisors have zero such intersection.
Consequently \(\zeta=K_Y^2Z_4>0\).

We need two intersection estimates. First, if \(rK_Y\)
is Cartier and \(B\) is an integral Weil divisor, then
\(rK_Y^2B\) is integral. Take general smooth surfaces in
\(|rK_Y+nH|\) and \(|nH|\), avoiding the isolated singularities;
the difference of their integral intersections of \(K_Y,B\)
is precisely this number. Here \(r=255\), so
\[
d=K_Y^2T=4V-\zeta\le\frac{31}{255}.                     \tag{5.22}\label{eq:5.23}
\]
Second, the original special-pencil model gives
\[
a\le d+\frac{d^2}{V}.                                   \tag{5.23}\label{eq:5.24}
\]
To see this, restrict
\(\mu^*(K_Y+T)=K_{W_{\mathrm{sp}}}+F+G\), \(G\ge0\),
to a general fibre and intersect with \(D\):
\[
a\le DK_F\le d+K_YT^2.
\]
The threefold Hodge index inequality gives
\((K_Y^2T)^2\ge K_Y^3(K_YT^2)\).
It follows by applying surface Hodge index on a smooth general
member of an ample multiple of \(K_Y+\varepsilon H\) and
letting \(\varepsilon\to0\); no nefness of \(T\) is assumed.
This proves \eqref{eq:5.24}.

Combining \eqref{eq:5.23}--\eqref{eq:5.24} yields
\[
\frac35=a\le\frac{31}{255}+\frac{(31/255)^2}{8/255}
=\frac{403}{680}<\frac35,
\]
a contradiction. We have therefore excluded every marking in \eqref{eq:5.20}.

\needspace{3\baselineskip}

\hypertarget{result-5.9}{}\ResultHeading{Theorem 5.10.} Under \(p_g(X)>0\) and \(V<1/30\), the
four-canonical pencil cannot have minimal fibre of type
\((2,2,0)\). \(\square\)

\needspace{4\baselineskip}

\hypertarget{quartic-pencils-of-12-surfaces}{%
\needspace{10\baselineskip}\hypertarget{section-6}{}\section{\texorpdfstring{6. Quartic pencils of \((1,2)\) surfaces}{6. Quartic pencils of (1,2) surfaces}}\label{quartic-pencils-of-12-surfaces}}

We now treat the \((1,2)\) fibres. Retaining the notation of the preliminaries, we write \(D=Q|_F\), \(A=\sigma^*K_{F_0}\), \(a=DA\), \(d=D^2\), and \(V=Q^3\). The earlier reduction gives
\[
p_g=P_2=P_3=\chi=h^{2,0}=1,\quad q=0,\quad P_4=2,\quad
1/48\le V<1/30.
\]
The canonical restriction and actual marking satisfy
\[
Q\ge_{\mathbb Q}F/4,\qquad D\ge_{\mathbb Q}A/5,\qquad
5/14\le\xi\le a\le20V<2/3.
\tag{6.1}\label{eq:6.1}
\]
On the original special pencil model the selected divisors satisfy \(0\le D\le K_F\le5D\). This actual comparison is used only there, or after pulling back both sides.

\needspace{4\baselineskip}

\hypertarget{the-finite-numerical-list}{%
\hypertarget{section-6.1}{}\subsection{6.1. The finite numerical list}\label{the-finite-numerical-list}}

We apply the recursion of \hyperlink{section-3.3}{Section 3.3} with the budgets \((\sum\gamma,\sum d)=(5,9)\) and
\(r\le25\) from \cite[Lemma 4.9]{ref1}. This gives 38 necessary possibilities. Of these, seven have
\(r_YV=1\) and violate the product inequality of \cite[Lemma 3.13]{ref1} in the
saturated canonical quotient. The following table records a violating pair of degrees for each of these seven baskets.
Here \(h_m=\dim\mathcal S_m=P_m-P_{m-1}\) except that
\(h_2=P_2-P_1+1\), by \cite[Proposition 3.3]{ref1}. The saturation is a domain
because it is a subring of the fraction field of \(R/(s)\).

\begin{longtable}[]{@{}lrl@{}}
\toprule\noalign{}
Discarded basket & \(V\) & \((i,j;h_i,h_j,h_{i+j})\) \\
\midrule\noalign{}
\endhead
\bottomrule\noalign{}
\endlastfoot
\(3(1,2)+(1,3)+(3,8)+(7,16)\) & \(1/48\) & \((2,15;1,2,1)\) \\
\(5(1,2)+(9,23)\) & \(1/46\) & \((2,4;1,1,0)\) \\
\(7(1,2)+2(1,3)+(2,7)+(3,7)\) & \(1/42\) & \((2,6;1,2,1)\) \\
\(3(1,2)+2(1,3)+3(3,7)\) & \(1/42\) & \((2,6;1,1,0)\) \\
\(5(1,2)+3(2,5)+(3,8)\) & \(1/40\) & \((2,4;1,1,0)\) \\
\(6(1,2)+3(1,3)+(1,4)+(4,9)\) & \(1/36\) & \((2,8;1,2,1)\) \\
\(3(1,2)+2(1,3)+(4,9)+(5,12)\) & \(1/36\) & \((2,11;1,2,1)\) \\
\end{longtable}

Thus each displayed row satisfies \(r_YV=1\) but violates \(h_{i+j}\ge h_i+h_j-1\). We are left with the following 31 necessary possibilities.

\begin{longtable}[]{@{}
  >{\raggedleft\arraybackslash}p{(\columnwidth - 8\tabcolsep) * \real{0.0450}}
  >{\raggedright\arraybackslash}p{(\columnwidth - 8\tabcolsep) * \real{0.4300}}
  >{\raggedleft\arraybackslash}p{(\columnwidth - 8\tabcolsep) * \real{0.1000}}
  >{\raggedright\arraybackslash}p{(\columnwidth - 8\tabcolsep) * \real{0.1550}}
  >{\raggedright\arraybackslash}p{(\columnwidth - 8\tabcolsep) * \real{0.2700}}@{}}
\toprule\noalign{}
\begin{minipage}[b]{\linewidth}\raggedleft
No.
\end{minipage} & \begin{minipage}[b]{\linewidth}\raggedright
Basket
\end{minipage} & \begin{minipage}[b]{\linewidth}\raggedleft
\(V\)
\end{minipage} & \begin{minipage}[b]{\linewidth}\raggedright
\((P_5,P_6,P_7)\)
\end{minipage} & \begin{minipage}[b]{\linewidth}\raggedright
\((P_8,P_9,P_{10},P_{12},P_{13})\)
\end{minipage} \\
\midrule\noalign{}
\endhead
\bottomrule\noalign{}
\endlastfoot
1 & \((1,3)+(5,13)+(8,17)\) & \(14/663\) & \((3,4,5)\) & \((7,8,10,14,17)\) \\
2 & \((1,2)+(2,5)+(4,11)+(7,15)\) & \(7/330\) & \((3,4,5)\) & \((7,8,10,13,16)\) \\
3 & \((1,2)+2(1,3)+(2,5)+(3,7)+(6,13)\) & \(59/2730\) & \((3,4,5)\) & \((6,8,9,13,15)\) \\
4 & \(3(1,2)+(1,3)+(3,7)+(3,8)+(4,9)\) & \(11/504\) & \((3,4,5)\) & \((6,7,8,12,14)\) \\
5 & \(5(1,2)+(2,5)+(7,18)\) & \(1/45\) & \((3,3,4)\) & \((6,6,9,12,14)\) \\
6 & \(2(1,2)+2(1,3)+(2,5)+2(4,9)\) & \(1/45\) & \((3,4,5)\) & \((6,8,8,12,14)\) \\
7 & \(5(1,2)+3(1,3)+(1,4)+(5,11)\) & \(1/44\) & \((2,4,4)\) & \((6,7,8,12,13)\) \\
8 & \(3(1,2)+2(3,8)+(5,11)\) & \(1/44\) & \((3,4,5)\) & \((7,7,9,12,15)\) \\
9 & \(2(1,2)+2(1,3)+(5,11)+(5,12)\) & \(1/44\) & \((3,4,5)\) & \((6,8,9,13,14)\) \\
10 & \((1,3)+(2,5)+(3,8)+(8,17)\) & \(47/2040\) & \((3,4,5)\) & \((7,8,10,14,17)\) \\
11 & \(5(1,2)+2(2,5)+(5,13)\) & \(3/130\) & \((3,3,4)\) & \((6,6,9,12,14)\) \\
12 & \((1,2)+(1,3)+(5,13)+(7,15)\) & \(3/130\) & \((3,4,5)\) & \((7,8,10,14,17)\) \\
13 & \(2(1,2)+(2,5)+(4,11)+(6,13)\) & \(17/715\) & \((3,4,5)\) & \((7,8,10,13,16)\) \\
14 & \(8(1,2)+(3,8)+(3,10)\) & \(1/40\) & \((2,4,4)\) & \((6,6,9,12,14)\) \\
15 & \((1,2)+(1,3)+(2,5)+(3,8)+(7,15)\) & \(1/40\) & \((3,4,5)\) & \((7,8,10,14,17)\) \\
16 & \(2(1,2)+2(1,3)+(2,5)+(3,7)+(5,11)\) & \(29/1155\) & \((3,4,5)\) & \((6,8,9,13,15)\) \\
17 & \(2(1,2)+(1,3)+(5,13)+(6,13)\) & \(1/39\) & \((3,4,5)\) & \((7,8,10,14,17)\) \\
18 & \(8(1,2)+(2,7)+(4,11)\) & \(2/77\) & \((2,4,4)\) & \((6,6,9,12,14)\) \\
19 & \(4(1,2)+2(3,7)+(4,11)\) & \(2/77\) & \((3,4,5)\) & \((6,7,9,13,15)\) \\
20 & \(3(1,2)+(2,5)+(4,11)+(5,11)\) & \(3/110\) & \((3,4,5)\) & \((7,8,10,13,16)\) \\
21 & \(2(1,2)+(1,3)+(2,5)+(3,8)+(6,13)\) & \(43/1560\) & \((3,4,5)\) & \((7,8,10,14,17)\) \\
22 & \(4(1,2)+2(3,8)+(4,9)\) & \(1/36\) & \((3,4,5)\) & \((7,7,9,13,16)\) \\
23 & \(3(1,2)+(1,3)+(5,11)+(5,13)\) & \(25/858\) & \((3,4,5)\) & \((7,8,10,14,17)\) \\
24 & \(3(1,2)+2(1,3)+(2,5)+(7,16)\) & \(7/240\) & \((3,4,5)\) & \((6,8,9,14,16)\) \\
25 & \(8(1,2)+(1,3)+(2,7)+(3,8)\) & \(5/168\) & \((2,4,4)\) & \((6,6,9,13,15)\) \\
26 & \(4(1,2)+(1,3)+2(3,7)+(3,8)\) & \(5/168\) & \((3,4,5)\) & \((6,7,9,14,16)\) \\
27 & \(3(1,2)+2(1,3)+(2,5)+(3,7)+(4,9)\) & \(19/630\) & \((3,4,5)\) & \((6,8,9,14,16)\) \\
28 & \(8(1,2)+(2,5)+(4,13)\) & \(2/65\) & \((2,4,4)\) & \((6,7,10,14,16)\) \\
29 & \(3(1,2)+(1,3)+(2,5)+(3,8)+(5,11)\) & \(41/1320\) & \((3,4,5)\) & \((7,8,10,14,17)\) \\
30 & \(2(1,3)+2(2,5)+(8,17)\) & \(8/255\) & \((3,4,5)\) & \((7,9,11,16,19)\) \\
31 & \(4(1,2)+(2,5)+(4,9)+(4,11)\) & \(16/495\) & \((3,4,5)\) & \((7,8,10,14,17)\) \\
\end{longtable}

Table 6.1. The necessary baskets in the \((1,2)\)-fibre case.

We next use the differences \(A_n=P_n-P_{n-1}\) to detect a fixed divisor. If \(s\) and \(x\in R_4\) have no common prime on the normal canonical model, multiplication by \(x\) induces an injection
\[R_n/sR_{n-1}\longrightarrow R_{n+4}/sR_{n+3}.\]
Indeed, \(xb=sa\) makes \(b/s\) regular at every prime divisor, retaining the full multiplicities of \(\operatorname{div}(s)\); reflexivity then extends it globally. Hence a decrease \(A_n>A_{n+4}\) detects a fixed prime of the quartic system on the canonical model.

Write \(4K_Y\sim T+Z_4\), where \(T\) is the image of a general connected fibre, and put \(\zeta=K_Y^2Z_4\). A fixed prime surviving on the canonical model has positive canonical intersection because \(K_Z\) is ample. Cartier-index integrality therefore gives \(\zeta\ge1/r_Y\) and
\[d=K_Y^2T=4V-\zeta\le4V-1/r_Y.\]
On the special pencil model, restricting
\(\pi^*(K_Y+T)=K_{W_{\rm sp}}+F+G\), \(G\ge0\), and intersecting with \(D\) gives
\[a\le DK_F\le d+K_YT^2\le d+d^2/V.\]
The last inequality is threefold Hodge index with nef and big polarization \(K_Y\), obtained by ample perturbation; \(T\) need not be nef. The right side increases for \(d\ge0\), so substitution yields
\[
a\le U(V,r_Y)=20V-9/r_Y+1/(r_Y^2V).
\tag{6.2}\label{eq:6.2}
\]
Applying \eqref{eq:6.2} to the following decreases in the Hilbert function removes five further rows.

\begin{longtable}[]{@{}rrlr@{}}
\toprule\noalign{}
Row & \(r_Y\) & \((n,A_n,A_{n+4})\) & \(U\) \\
\midrule\noalign{}
\endhead
\bottomrule\noalign{}
\endlastfoot
5 & 90 & \((5,1,0)\) & \(7/20\) \\
6 & 90 & \((6,1,0)\) & \(7/20\) \\
7 & 132 & \((6,2,1)\) & \(7/18\) \\
8 & 88 & \((5,1,0)\) & \(63/176\) \\
9 & 132 & \((9,2,1)\) & \(7/18\) \\
\end{longtable}

We now explain these exclusions. Every upper bound in the last column is less than \(2/5\), so the base-point comparison gives \(a=\xi\); the denominator is consequently cleared by one actual local index. Rows 5 and 6 contradict \(a\ge5/14\). In row 8 no fraction cleared by one of the actual indices \(2,8,11\) lies in \([5/14,63/176]\). Rows 7 and 9 require \(a=4/11\) at the singleton quotient \((5,11)\). Its index-one-cover canonical coefficient has character \(-5\). A coefficient of value four would require a monomial of weight four for positive coordinate valuations congruent to \((8,3,7)\) modulo eleven. There is no such monomial.

Row 11 has \(P_5=P_6=3\), contrary to Corollary~2.2. In particular, every remaining row has \(P_6=4\).

\needspace{4\baselineskip}

\hypertarget{surface-images-and-the-quadratic-sixth-field}{%
\hypertarget{section-6.2}{}\subsection{6.2. Surface images and the quadratic sixth field}\label{surface-images-and-the-quadratic-sixth-field}}

\needspace{3\baselineskip}

\hypertarget{result-6.1}{}\ResultHeading{Lemma 6.1.} Apart from the possible surface alternative in row 31, the sixth map is generically finite of degree two. In that exceptional alternative \(a=\xi=4/11\).

\ProofHeading{Proof.} We write \(K_0=\mathbb C(t,t')\) for the field of the connected canonical-curve fibration and first consider the case of a surface image. A surface sixth image has function field \(K_0\) and restriction rank two or three. Rank three gives \(6D\ge_{\mathbb Q}2C\), hence
\[
V\ge d/4\ge\xi/12.
\tag{6.3}\label{eq:6.3}
\]
At rank two, the two-dimensional restriction kernel supplies a moving subsystem \(S_1\) with \(M_6\ge F+S_1\). If it is the outer pencil, \(6Q\ge_{\mathbb Q}2F\) and canonical restriction again gives \eqref{eq:6.3}. Otherwise its restriction is composed with the canonical-curve pencil. The numerical clause of \cite[Proposition 3.6(2.1)]{ref10}, with \(m_1=6,n=3\), gives
\[
10\xi\ge\lceil3\xi\rceil+2\ge4,\qquad \xi\ge2/5.
\]
The separation hypotheses hold in degree six because \(t,t'\) are present.

The surface-image estimate \(36V\ge\deg(\Sigma_6)\xi\) makes the image a quadric. The globally generated rank-two evaluation bundle has four sections; a summand of degree at least two would give the previous outer-pencil alternative. Hence it is \(\mathcal O(1)^{\oplus2}\). On a common resolution the sixth moving part is \(F+G\), where \(G\) is the nef second ruling and \(G|_F=C\). Therefore
\[
6V\ge d+\xi/4.
\]
Choose \(u_0\) with \(E_6=\langle1,t,u_0,tu_0\rangle\), and write \(t'=A(t)+B(t)u_0\) for affine \(A,B\), with \(B\ne0\). The identity
\[
(tu_0)'-tu_0'-A(t)u_0=B(t)u_0^2
\]
places \(1,u_0,u_0^2\) in \(\mathbb C(t)E_8\). Restricting to a general outer fibre gives an actual subsystem with moving part \(2C\), so \(d\ge\xi/4\). Thus \(V\ge\xi/12\ge1/30\), impossible.

We are thus left with the alternatives satisfying \eqref{eq:6.3}. Among the rows of Table 6.1, a fraction in \([5/14,12V]\) cleared by one actual basket index exists only in rows 29 and 31; it is \(4/11\) in both. Since it is below \(2/5\), the base-point comparison gives \(a=\xi\). Row 29 is excluded by the singleton \((5,11)\) character test. This proves the surface assertion.

For a generically finite sixth map write \(E_6=\langle1,t,t',h\rangle\). On the generic canonical curve, \(h\) is nonconstant and its only pole is the Weierstrass point, of order at most \(\lfloor6\xi\rfloor\le3\). It generates the hyperelliptic field. Connectedness therefore gives \([\mathbb C(X):L_6]=2\). \(\square\)

\hypertarget{lem-canonical-subring-data}{}\ResultHeading{Lemma 6.2 (The sixth subring and eighth functions).} For each remaining generically finite row with \((P_5,P_6,P_7)=(3,4,5)\), the sixth subring is polynomial of weights \((1,4,5,6)\); its new seventh section is rational with primitive denominator degree at most six, and every eighth function lies in the sixth field.

\ProofHeading{Proof.} For \((P_5,P_6,P_7)=(3,4,5)\), the fifth map is not an outer pencil: that would give \(5Q\ge_{\mathbb Q}2F\), and hence \(V\ge8/245\), larger than every volume in the table. Thus
\[
B=\mathbb C[s,x,y,z],\qquad\deg(s,x,y,z)=(1,4,5,6)
\]
is polynomial. Write \((t,u,z_0)=(x/s^4,y/s^5,z/s^6)\) and \(K=L_6\). Since \(7\xi<5\), a new seventh section is rational over \(B\). Its primitive relation \(FW=G\) satisfies
\[
f=\deg F\ge1,\qquad V\ge(7+f)/420,\qquad f\le6.
\tag{6.4}\label{eq:6.4}
\]
All eighth functions are in \(K\) for \(V<1/32\), because \(8\xi\le160V<5\). Only rows 30 and 31 exceed this cutoff. In row 30 the fixed-part witness \((12,3,2)\) gives \(a\le403/680<5/8\). In row 31, if \(f\ge2\), then \(h_{12}(B[W])\ge13\); an eighth section outside the fraction field would add two independent twelfth sections after multiplication by \(B_4\), contradicting \(P_{12}=14\). The alternative \(f=1\) would give \(sW=x^2\), excluded by Corollary~2.2. Consequently every eighth function belongs to \(K\) in each remaining case.

\needspace{4\baselineskip}

\(\square\)

\hypertarget{uniform-algebraic-exclusions}{%
\hypertarget{section-6.3}{}\subsection{6.3. Uniform algebraic exclusions}\label{uniform-algebraic-exclusions}}

We formulate the next denominator lemmas without specifying the fibre type, so that they apply both here and in \hyperlink{section-5}{Section 5}. In each lemma, \((P_1,\ldots,P_7)=(1,1,1,2,3,4,5)\), \(q=0\), \(h^{2,0}=1\), \(V<1/30\), and the Hodge and Euler identities of the preliminaries hold. Moreover, \(B\) is the polynomial sixth subring of weights \((1,4,5,6)\), and \(W\) is a new seventh section with primitive denominator \(F\).

\needspace{3\baselineskip}

\hypertarget{result-6.2}{}\ResultHeading{Lemma 6.3.} If \((P_4,P_5,P_6,P_7)=(2,3,4,5)\) and the sixth field has prime degree, then \(\deg F=2\) is impossible.

\ProofHeading{Proof.} After subtracting an old seventh section, we reduce the relation to \(s^2W=xy\), so that \(w=tu\). By Corollary~2.2, \(t'\notin E_5\); use \(z_0=t'\). The second residue identity gives constant-coefficient equations
\[
u'=\alpha tu+Az_0+Bu+Ct+D,\qquad
z_0'=kt^2+Ftu+Gz_0+Hu+It+J.
\]
If \(\alpha\ne0\), a component detecting the \(y\) residue has \(u\)-pole \(5\ell\), whereas \(xy=s^2W\) bounds the \(t\)-pole by \(2\ell\). The first equation forces the latter pole to be \(\ell+1\). In the second equation, \(tu\) first forces \(F=0\), then \(Hu\) is either absent or less polar. Every remaining term has pole below the left-hand pole \(3\ell+3\), a contradiction. Thus \(\alpha=0\). Setting \(U=u-At\) gives
\[
U'=BU+Mt+D,\qquad
t''=at^2+ftU+bt'+cU+dt+e.
\tag{6.5}\label{eq:6.5}
\]
If \(a\ne0\), a pole of \(t\) has order two and indicial roots \(-1,6\), together with the integral first-order root. If \(a=0\), a pole requires \(fM\ne0\) and has order three; elimination gives leading equation \(U'''=fUU'\) at a double pole, with roots \(-1,4,6\). At a bounded point the system is polynomial. The actual expansions are integral, contrary to the Jacobian lemma.

\(\square\)

\needspace{3\baselineskip}

\hypertarget{result-6.3}{}\ResultHeading{Lemma 6.4.} In the quadratic case, if \(P_7=5\), \(P_8=7\), \(P_{12}=14\) and \(E_8\subset K\), then \(\deg F=5\) is impossible.

\ProofHeading{Proof.} Put \(\mathcal A=B[W]\); its eighth and twelfth dimensions are six and thirteen. First
\[
s^4Q,xQ\in\mathcal A_{12}\quad\Longrightarrow\quad Q\in\mathcal A_8.
\tag{6.6}\label{eq:6.6}
\]
Indeed write \(s^4Q=a+bW\), \(xQ=c+dW\). Primitiveness gives
\[
xb-s^4d=Fh,\qquad xa-s^4c=-Gh,\qquad h=\alpha s^4+\beta x.
\]
Coprimeness gives \(b-\beta F=s^4k_1\), \(a+\beta G=s^4k_8\), and therefore \(Q=k_8+k_1W\).

If \(\mathcal A\) is canonical, its residue form supplies an odd section in degree twelve. The invariant twelfth space is then \(\mathcal A_{12}\), and the new invariant eighth section contradicts \eqref{eq:6.6}. To complete the argument, we must examine the cases in which canonicity is not immediate. In the exceptional cases we will use the normality of the full canonical ring.

A nonzero fifth-coordinate coefficient in \(F\) gives, after changing that coordinate, \(yW=g(s,x,z)\), which is canonical. Otherwise \(F=s^5\) or \(sx\). For \(F=s^5\), the reduced numerator is
\[
G=az^2+dx^3+bsyz+cs^2y^2+es^2xz+fs^3xy+gs^4x^2,\qquad(a,d)\ne(0,0).
\]
If \(a\ne0\), complete the square in \(z\). If \(d=0\), \(z'/s\) is integral of degree five, contradicting the new sixth generator. If \(d\ne0\), the slice \(y=0\) is
\[
z'^2+dx^3+g's^4x^2-s^5W=0.
\]
Its weight-\((1,2,3,1)\) blowup is crepant. The exceptional surface \(z'^2+dx^3-s^5W=0\) avoids the two ambient quotient coordinate points and has only an \(E_8\) point. It is normal and klt; inversion of adjunction proves canonicity of the slice and then of \(\mathcal A\).

If \(a=0,b\ne0\), legal changes give \(s^5W=x^3+syz+cs^4x^2\). The slice \(W=0\) has a crepant ordinary blowup whose exceptional cubic \(x^3+syz=0\) is normal with three \(A_2\) points. The same argument applies. If \(a=b=0\), write \(x^3=s^2H\). Then \((x^2/s)^3=sH^2\), so \(x^2/s\) is an actual seventh section. It cannot be in \(B_7+\mathbb CW\), by reduction modulo \(s\), a contradiction.

For \(F=sx\), the reduced numerator is
\[
G=az^2+dx^3+bsyz+cs^2y^2+es^6z+fs^7y+gs^{12}.
\]
If \(a\ne0\), complete the square; the successive slices \(y=0\) and \(W=\lambda s+\mu x\) give a \(D_4\) surface for general \(\lambda,\mu\), even when \(d=0\). If \(a=0,b\ne0\), the equation reduces to \(sxW=dx^3+syz+gs^{12}\), with \(d\ne0\), and the preceding cubic slice proves canonicity. Finally if \(a=b=0\), write \(G=dx^3+s^2H\). Then \(T=x^2/s\) satisfies \(dT^2-WT+xH=0\); normality again contradicts primitiveness. This treats every possible denominator of degree five. \(\square\)

\needspace{3\baselineskip}

\hypertarget{result-6.4}{}\ResultHeading{Lemma 6.5.} In the quadratic sixth-field setting, if \(P_7=5\), \(P_8=6\) and \(P_{10}=9\), then \(\deg F=5\) is impossible.

\ProofHeading{Proof.} Here \(R_8=\mathcal A_8\) and \(R_{10}=\mathcal A_{10}\) for \(\mathcal A=B[W]\), whose first relation has degree twelve. Write \(D_0\) for the sixth-coordinate coefficient in \(t'\) and \(e\) for the \(w\) coefficient in \(u'\). The degree-eleven Euler commutator, reduced modulo \(s\) in \(\mathcal A\), reads
\[
5D_0yz-4exW=0.
\]
There is no relation of that degree, so \(D_0=e=0\). In particular \(t'\in E_5\), contradicting Corollary~2.2. \(\square\)

\needspace{3\baselineskip}

\hypertarget{result-6.5}{}\ResultHeading{Lemma 6.6.} In the quadratic case, if \(P_7=5\), \(P_8=7\), \(P_{12}=14\) and \(E_8\subset K\), then \(\deg F=6\) is impossible.

\ProofHeading{Proof.} The subring \(\mathcal A=B[W]\) has \(h_{12}=14\), so \(R_{12}=\mathcal A_{12}\). For a homogeneous degree-eight fraction \(q\) with \(s^4q,xq\in\mathcal A_{12}\), comparison of the coefficient of \(W\) gives
\[
xa-s^4c+(xb-s^4d)W=\lambda s^3(FW-G).
\]
If \(\lambda=0\), then \(q\in\mathcal A_8\). Otherwise \(F=sF_5\), \(G=x^2y+sH_{12}\) after rescaling, and
\[
\{q:s^4q,xq\in\mathcal A_{12}\}=\mathcal A_8+\mathbb C\,xy/s.
\tag{6.7}\label{eq:6.8}
\]
A new invariant eighth section exists, so \(Q=xy/s\) is actual and \(R_8=\mathcal A_8+\mathbb CQ\).

If \(F_5\) contains a fifth coordinate \(Y\), promote \(sY\) to a weight-six coordinate, clear powers of \(s\) in the numerator, and divide by that coordinate. The resulting fixed hypersurface is \(TU=g(s,x,z)\). Its residue form gives an odd section in degree \(1+(1+4+6)=12\), impossible. Thus \(F=s^6\) or \(s^2x\).

Write \(\Delta y=eW\bmod sR_6\) and \(\Delta z=kx^2+hQ\bmod sR_7\). The fifth/sixth commutator lies in \(sR_{12}=s\mathcal A_{12}\). The term \(Qy\) is not in \(\mathcal A_{13}\), so \(h=0\). Reducing modulo \(s\), the sole degree-thirteen relation is \(x^2y=0\), which forces \(e=0\). Consequently
\[
u'\in E_6,\qquad z_0'\in E_7+\mathbb Ct^2.
\tag{6.8}\label{eq:6.9}
\]

Our next step is to show that the coefficient of \(sz^2\) in \(G\) is zero. Completing the square and making triangular changes gives, respectively,
\[
x^2y+sz^2+cs^3y^2-s^6W,
\]
or
\[
x^2y+sz^2+cs^3y^2+ds^8y+es^{13}-s^2xW.
\]
When \(c\ne0\), replace \(y\) by \(Y-ds^5/(2c)\), and adjust \(W\), before testing Newton faces. The remaining terms in \(s,y\) are then \(cs^3Y^2+(e-d^2/(4c))s^{13}\). Each Newton face is nondegenerate: successively the \(W,z,x\) derivatives handle a face containing their corresponding mandatory monomials, and the remaining one- or two-term faces are smooth on the torus. The hypersurfaces are normal since their singular sets have codimension at least two. If \(m\) is the weight of the equation for a nonzero nonnegative valuation vector, the mandatory monomials give
\[
\sum\nu_i\ge\frac{13}{12}m
\quad\text{or}\quad
\sum\nu_i\ge\frac98m,
\]
For the first bound use coefficients \((1/2,1/2,1/12)\) on
\(x^2y,sz^2,s^6W\); for the second use \((3/8,1/2,1/4)\) on
\(x^2y,sz^2,s^2xW\). The Newton discrepancy criterion \cite{ref15} proves canonicity: for a nondegenerate hypersurface, a regular subdivision of the Newton fan resolves it, and the discrepancy of a primitive integral weight \(w\) is \(\sum w_i-1-w(F)\). The displayed strict weight bounds make this nonnegative. The residue form then gives an odd section in degree eleven, contradicting \(R_{12}=\mathcal A_{12}\). Thus that coefficient vanishes.

By Corollary~2.2, \(t'\notin E_5\). Take \(z_0=t'\), and replace \(u\) by \(U=u-at\) so that
\[
U'=bU+ct+d,\qquad
t''=24\lambda t^2+A_0w+B_0t'+C_0U+D_0t+E_0.
\tag{6.9}\label{eq:6.10}
\]
Choose the reduced \(x\)-unit component of Corollary~2.2, so \(\ell=1\). At this component, \(t,t'\) have \(\rho\)-poles \(4,6\); \(U\) is bounded if \(c=0\), and otherwise has pole \(3\ell-1\). For \(F=s^6\), the complete expression for \(w\), modulo affine terms, is
\[
t^2U+\beta t^3+\eta Ut'+\kappa U^2+rt t'+\sigma tU+\gamma t^2.
\]
Its seventh-section bound first forces \(\beta=0\), then \(c\ne0\) and \(\ell=1\). If \(A_0=0\), the indicial roots of \eqref{eq:6.10} are \(-1,6,1\). If \(A_0\ne0\), the highest-family first-shift coefficient, after eliminating the first-order variation, is a nonzero multiple of \((q-2)(q+1)/(q-1)\). All exceptional rational shifts are integral.

For \(F=s^2x\) the displayed expression is divided by \(t\). At \(\ell=1\) the same leading second-order system applies. In every case the actual expansions are integral, contrary to the Jacobian lemma. Thus the alternative \(t'\notin E_5\), which is forced by Corollary~2.2, is impossible. \(\square\)

\needspace{3\baselineskip}

\hypertarget{result-6.6}{}\ResultHeading{Lemma 6.7.} Suppose the sixth field is quadratic, \(P_7=5\), \(P_8=7\), \(P_9\ge8\), \(P_{16}\le29\), \(E_8\subset K\), and \(V\le8/255\). Then \(\deg F=5\) or \(6\) is impossible.

\ProofHeading{Proof.} We use the simultaneous polar-divisor inequality from the preliminaries. In a normal four-dimensional graded ring generated in weights \(w_1\le\cdots\le w_N\), a new rational degree-\(m\) element increases multiplicity by at least
\[
\frac1{m w_{N-2}w_{N-1}w_N}.
\tag{6.10}\label{eq:6.11}
\]
Indeed the images of three generators contain a residue transcendence basis at a polar prime. If the prime does not contain \(s\), one may select \(s\) as one of them.

First let \(\deg F=6\). The seventh subring has multiplicity \(13/840\). If normal, a new invariant eighth section gives
\[
V\ge2\left(\frac{13}{840}+\frac1{8\cdot5\cdot6\cdot7}\right)=\frac9{280}>\frac8{255}.
\]
After legal changes, \(F\) is \(z,sy,s^2x\) or \(s^6\). The first two hypersurfaces are normal: at every generic boundary divisor a derivative involving the free \(W\) coordinate or a primitive numerator is a unit. In the latter two cases, primitivity gives \(G=x^2y+sH_{12}\); the only possible nonnormal prime is \((s,x)\). A nonzero \(z^2\) coefficient in \(H_{12}\bmod s\) makes it regular as well. Otherwise shift \(y\) to write \(G=x^2y+s^2H_{11}\).

Set \(L=xW-H_{11}\) for \(F=s^2x\), and \(L=s^4W-H_{11}\) for \(F=s^6\). Then \(x^2y=s^2L\) and \(Q_0=xy/s\) is integral, since \(Q_0^2=yL\). The ring \(\mathcal A'=\mathcal A[Q_0]\) is the determinantal graph
\[
sQ_0=xy,\qquad xQ_0=sL,\qquad Q_0^2=yL.
\]
It has codimension two, is Cohen--Macaulay, and is a domain: localization at \(s\) is the graph, while every component on \(s=0\) has dimension at most three. At the possibly new boundary, write \(H_{11}\bmod(s,x)=a yz\). If \(a\ne0\), the generic equation \(Q_0^2+a y^2z=0\) is separable and regular. If \(a=0\), write \(H_{11}=xH_7+sH_{10}\) and \(H_{10}\bmod(s,x)=by^2\). For \(b\ne0\) the local equation after substituting \(x=sQ_0/y\) has a nonzero linear \(s\) coefficient. Thus \(\mathcal A'\) is normal.

If \(a=b=0\), write \(H_{10}=xH_6+sH_9\). The fraction \(Z=xy/s^2\) is integral: in the two cases its monic equation is
\[
Z^2-Z(W-H_7)+sZH_6+yH_9=0,
\]
or
\[
Z^2+ZH_7+sZH_6+yH_9-s^2yW=0.
\]
This gives a new seventh fraction of denominator degree two, contradicting primitiveness. Hence \(\mathcal A'\) is always normal in a surviving case.

The extension is finite birational, so \(e(\mathcal A')=13/840\). Its ninth piece has dimension seven. We claim that every ninth section is invariant. Indeed, an odd ninth section would give five independent odd sixteenth multiples by \(\mathcal A_7\), in addition to the 25 invariant sections of \(\mathcal A_{16}\), contradicting \(P_{16}\le29\). A new rational ninth section and \eqref{eq:6.11} give
\[
V\ge2\left(\frac{13}{840}+\frac1{9\cdot6\cdot7\cdot8}\right)
=\frac{239}{7560}>\frac8{255}.
\]

Now let \(\deg F=5\). The normal forms and integral-dependence exceptions in the proof of \hyperlink{result-6.3}{Lemma 6.4} still apply. In every canonical case \(e(\mathcal A)=1/70\). A finite polar prime of a new invariant eighth section has residue transcendence-basis weight product at most \(42\), so it would give \(V\ge29/840\). Thus all poles lie on \(s=0\).

If \(F=y\), write \(yW=g(s,x,z)\). If \(g\bmod s=az^2+dx^3\ne0\), the boundary is irreducible with quotient multiplicity \(1/70\); one pole already gives \(V\ge9/280\). Otherwise \(g=s^eH\), with \(e\ge2\). The two boundary primes \((s,y),(s,W)\) have residue products \(168,120\). A pole of order two at either is too costly, so \(sQ\in\mathcal A_9\). Modulo old sections, necessarily \(Q=xy/s\).

The graph \(\mathcal A'=\mathcal A[Q]\) is defined by the minors of
\[
\begin{pmatrix}y&s^{e-1}H&Q\\s&W&x\end{pmatrix}.
\]
It is Cohen--Macaulay of codimension two and has no extra component on \(s=0\). Its boundary primes are \((s,y,W),(s,y,Q),(s,x,W)\); at their generic points the equations solve all but a uniformizer and the generic units. Hence it is normal. The actual graph has multiplicity
\[
e(\mathcal A')=\frac{9\cdot12}{6720}-\frac1{1344}=\frac{103}{6720}.
\]
The subtracted component of the degree-\((9,12)\) complete intersection is \((s,y)\) of generic length one. The three boundary quotient products are \(192,168,240\).

Its ninth dimension is still seven. Moreover \(\dim\mathcal A'_{16}\ge25\): \(\dim\mathcal A_{16}=24\), and \(x^2Q\) is new, as seen at \((s,W)\). Thus all ninth sections are invariant by the same five-multiple argument. A new ninth section has a finite-prime residue product at most \(56\), or a boundary product at most \(240\). Therefore
\[
V\ge\frac{103}{3360}+\frac1{1080}
=\frac{191}{6048}>\frac8{255}.
\]

For \(F=s^5\), the square case in \hyperlink{result-6.3}{Lemma 6.4} has irreducible boundary \(z'^2+dx^3=0\), quotient multiplicity \(1/70\), and is excluded by \(9/280\). The remaining canonical case is
\[
s^5W=x^3+syz+cs^4x^2.
\]
Its boundary has valuations \((3,1,0,0,0)\) on \((s,x,y,z,W)\) and residue product \(210\). A pole of order three is too costly; hence \(sQ\in\mathcal A_9\) and \(Q=xy/s\). For \(T=x^2/s\),
\[
T^2=s^3xW-Qz-cs^2x^3.
\]
Normality contradicts the primitive fifth denominator.

For \(F=sx\), the square-completed equation is
\[
sxW=z'^2+dx^3+s^2(cy^2+fs^5y+gs^{10}).
\]
If \(d\ne0\), use the same irreducible-boundary cost. If \(d=0\), the unique boundary has \(\nu(s)=2,\nu(z')=1\) and residue product \(140\). A pole of order two is too costly, so \(sQ\in\mathcal A_9\). Its only possible new quotient \(xy/s\) has pole order two and is inadmissible. Finally, the nonsquare canonical case is
\[
sxW=dx^3+syz+gs^{12},\qquad d\ne0.
\]
The preceding product-\(210\) bound again gives \(Q=xy/s\), and
\[
dT^2-WT+Qz+gs^{10}x=0,\qquad T=x^2/s,
\]
contradicts primitiveness. The noncanonical exceptions were already eliminated by integral dependence in \hyperlink{result-6.3}{Lemma 6.4}. This proves both degree exclusions. \(\square\)

\needspace{4\baselineskip}

\hypertarget{the-larger-volume-twelfth-space-alternatives}{%
\hypertarget{section-6.4}{}\subsection{6.4. The larger-volume twelfth-space alternatives}\label{the-larger-volume-twelfth-space-alternatives}}

The last numerical row has more room than the uniform volume estimates of \hyperlink{section-7}{Section 7}. Its small twelfth, thirteenth and fifteenth pieces supply the following replacements, which are again independent of the fibre type.

\needspace{3\baselineskip}

\hypertarget{result-6.7}{}\ResultHeading{Lemma 6.8.} Suppose the sixth field is quadratic, \(E_8\subset K\), and \((P_7,P_8,P_{12},P_{13})=(5,7,14,17)\). If \(V<11/336\), a primitive cubic seventh denominator is impossible. If also \(P_9=8\), a pure quartic denominator is impossible.

\ProofHeading{Proof.} The cubic denominator is excluded by the seventh-degree argument of \hyperlink{section-7.3}{Section 7.2}; that argument uses neither the bound on \(P_8\) nor the volume cutoff in Theorem~7.1(b).

For a pure quartic denominator write
\[
s^4W=yz+cs^3x^2,\qquad w=uz+ct^2.
\]
The odd thirteenth residue section gives \(\dim I_5\ge2\), and \(\dim I_4\ge1\) as before. Finite poles are impossible, so \(s^4Q\in\mathcal A_{12}\). At the two boundary primes the values of \(s^5,sx,y\) are
\[
(5,5),\quad(1,1),\quad(\nu,0),\qquad
\nu=3\text{ if }c\ne0,\quad\nu=4\text{ if }c=0.
\]
Hence either both poles are at most one, giving \(Q=xy/s\), or the second is zero and the first is at most \(\nu\). Inspection of \(\mathcal A_{11}/s^3\) or \(\mathcal A_{12}/s^4\) then gives \(Q=xz/s^2\); the only extra possible term \(z^2/s^4\) would make \(z/s^2\) integral of degree four. Generic residue independence excludes hidden cancellation.

If \(q=tz\), then \(E_9=\langle1,t,u,z,w,t^2,tz,tu\rangle\). Write
\[
t'=L_0,\qquad u'=L_1+ew,\qquad
z'=L_2+jw+kt^2+htz,
\]
where the \(L_i\) are affine. The cubic part of \(w'\) is
\[
(ez+ju)(uz+ct^2)+kut^2+htuz.
\]
Its coefficients force \(e=j=k=h=0\), making \(E_6\) stable.

If \(q=tu\), then \(\mathcal A[Q]\) has dimensions \(14,16\) in degrees \(12,13\). Its invariant thirteenth piece fills that of \(R\), because the odd thirteenth section supplies the remaining dimension. Every ninth section is invariant, since an odd one would have two independent odd thirteenth multiples. A new ninth function \(g\) and \(tg\) must lie in
\[
\begin{aligned}
M={}&\langle1,t,u,z,t^2,tu,tz,u^2,uz,t^3,z^2,t^2u\rangle\\
&+\langle tuz,u^2z,uz^2+ct^2z,tu^2\rangle.
\end{aligned}
\]
Monomial comparison gives
\[
M\cap t^{-1}M=\langle1,t,u,z,t^2,tu,uz,u^2\rangle.
\]
Thus the new ninth direction is \(u^2\). Replacing \(htz\) by \(htu\) in the derivative equations, the same cubic comparison makes \(E_6\) stable. This contradiction proves the pure-quartic assertion. \(\square\)

\needspace{3\baselineskip}

\hypertarget{result-6.8}{}\ResultHeading{Lemma 6.9.} Suppose the sixth field is quadratic, \(E_8\subset K\), and
\[
(P_7,P_8,P_9,P_{10},P_{12},P_{13},P_{15})
=(5,7,8,10,14,17,24).
\]
Then a finite quartic seventh denominator is impossible.

\ProofHeading{Proof.} Write \(xW=g(s,y,z)\), where
\[
g=ayz+bsy^2+cs^5z+ds^6y+es^{11}.
\]
This normal canonical hypersurface supplies an odd section \(\eta_{13}\). For a new invariant eighth section \(Q\), its annihilator dimensions satisfy \(\dim I_4\ge1\), \(\dim I_5\ge2\). A fourth annihilator is \(s^4\) or \(x-\lambda s^4\). If \(\lambda\ne0\), an independent fifth annihilator forces every pole onto \(s=x=0\). The possible positive \(s\)-orders of \(g\) are \(1,5,6,11\), none equal to four; hence both \(s^4\) and \(x\) clear these poles. We may therefore assume \(xQ\) or \(s^4Q\) is old.

If \(xQ\) is old, a fifth annihilator containing \(y\) gives \(g=yh_6\), \(Q=h_6H_6/x\), with \(y\nmid H_6\). Otherwise \(s^5\) clears the poles and gives \(g=s^5h_6\), \(Q=W(asy+bz)/s^5\). If \(s^4Q\) is old, write \(s^4Q=P_{12}+cyW\). The other fifth annihilator is \(sx\); divisibility forces \(s^3\mid P_{12}\), and \(c\ne0\) requires \(s^5\mid g\). Since \(B_9/sB_8\) is spanned by \(xy\), the complete list modulo old sections is
\[
\begin{array}{c|c}
g&Q\\ \hline
yh_6&h_6H_6/x,\quad y\nmid H_6\\
s^5h_6&W(asy+bz)/s^5\\
\text{arbitrary}&xy/s\\
s^5h_6&axy/s+cyW/s^4,\quad c\ne0.
\end{array}
\tag{6.11}\label{eq:6.13}
\]
Here \(h_6,H_6\) are linear in \(z,sy,s^6\).

In the first case the residue form
\[
\Omega_0=\frac{ds\wedge dx\wedge dy\wedge dz}{x}
\]
has a second regular quotient form \(\Omega_1=(h_6/x)\Omega_0\). To verify this, let \(a_f=\nu(f)\) for any prime divisorial valuation centred on the graph. Regularity of \(W,Q\) gives
\[
a_x\le a_y+a_h,\qquad a_x\le a_H+a_h.
\]
If \(a_x\le a_h\), use regularity of \(\Omega_0\). Otherwise \(a_y,a_H>0\). If \(H_6\) contains \(z\), use coordinates \(s,x,y,H_6\) and the differential estimate
\(\nu(df_1\wedge\cdots\wedge df_4)\ge\sum\nu(f_i)-1\) to get
\(\nu(\Omega_1)\ge a_s+a_H-1\ge0\).
If it does not contain \(z\), the condition \(y\nmid H_6\) gives a nonzero \(s^6\) coefficient, which forces \(a_s>0\); the original coordinates give
\(\nu(\Omega_1)\ge a_s+a_z-1\ge0\).
Thus \(\Omega_1\) supplies an odd fifteenth section independent of \(s^2\eta_{13}\).

There are also three new invariant fifteenth sections \(s^7Q,szQ,WQ\) over the old twenty-dimensional space. Multiply a proposed relation by \(x^2\) and reduce modulo \(x\). Its \(WQ\) coefficient would force \(H_6\) to be a multiple of \(sy\), so vanishes; the remaining coefficients vanish by coprimeness with \(y\). Consequently \(P_{15}\ge23+2=25\), impossible.

For the second case with \(b\ne0\), use \(H_6=asy+bz\) in the same form calculation. If \(a_x>a_h\), regularity of \(W\) gives \(a_s>0\), and \(\nu(\Omega_1)\ge a_s+a_y-1\ge0\). The three new invariant fifteenth sections are \(s^2yQ,szQ,WQ\). Their independence follows by reduction modulo \(x\) and \(\gcd(s,H_6)=1\). Again \(P_{15}\ge25\). When \(b=0\), the section belongs to the final row of \eqref{eq:6.13}.

The remaining nonpure possibilities are
\[
w=h(u,z)/t,\quad q=tu,\quad h\in\langle1,u,z,u^2,uz\rangle,
\]
or
\[
w=h(u,z)/t,\quad q=uw+ktu,\quad k\ne0,\quad
h\in\langle1,u,z\rangle.
\]
In both cases the actual twelfth space is
\[
E_{12}=M=B_{12}^0+\langle w,uw,t^2u\rangle,\qquad
B_{12}^0=\langle1,t,u,z,t^2,tu,tz,u^2,uz,z^2,t^3\rangle.
\tag{6.12}\label{eq:6.14}
\]
The fifth/sixth commutator reduced modulo this space has negative-\(t\) numerator \(h(6ez-5bu^2)\), where \(e\) is the \(w\) coefficient in \(u'\) and \(b\) the \(q\) coefficient in \(z'\). It forces \(e=b=0\). Hence \(u'\in E_6\) and \(z'\in E_7+\mathbb Ct^2\).

Corollary~2.2 gives \(t'\notin E_5\). We take \(z=t'\) and make a fifth-coordinate change \(U=u-at\) so that \(U'=bU+ct+d\). At the reduced \(x\)-unit component, \(t,t'\) have \(\rho\)-poles four and six, while \(U\) has pole at most two. Thus all terms of \(w\) have \(\rho\)-pole at most six. If \(c\ne0\), the dominant system is \(t''=24\lambda t^2\), \(U'=ct\), with indicial roots \(-1,6,1\). If \(c=0\), Hodge stability makes \(U\) time-constant, and the scalar indicial roots are \(-1,6\). In both cases the first-nonintegral-coefficient test gives integral actual-time expansions, contradicting the Jacobian lemma.

Finally the pure remaining case is \(q=uw\), with \(w=h(u,z)/t\) and \(h\) linear. The low dimensions give
\[
E_9=\langle1,t,u,z,t^2,tu,w,uw\rangle,
\]
\[
E_{10}=\langle1,t,u,z,t^2,tu,tz,u^2,w,uw\rangle.
\]
In \(q'=u'w+uw'\in E_{10}\) the sole double finite pole removes the \(w\) coefficient of \(u'\). The simple-pole numerator then removes both the \(z\) coefficient of \(u'\) and the \(uw\) coefficient of \(w'\). Thus \(u'\in\langle1,t,u\rangle\). The fourth/fifth commutator \(5t'u-4tu'\in E_{10}\) removes the \(z\) coefficient of \(t'\), because \(uz\) is absent from this space. Hence \(\langle1,t,u\rangle\) is stable, a contradiction. We have now exhausted the alternatives in \eqref{eq:6.13}. \(\square\)

\needspace{4\baselineskip}

\hypertarget{completion-of-the-numerical-reduction}{%
\hypertarget{section-6.5}{}\subsection{6.5. Completion of the numerical reduction}\label{completion-of-the-numerical-reduction}}

The four rows with \((P_5,P_6,P_7)=(2,4,4)\) are excluded by \hyperlink{result-5.6}{Lemma 5.8}. Rows 14,18,25 satisfy its part (a), since \(P_8=P_9=6\). Row 28 satisfies part (b), since \((P_8,P_9,P_{10},P_{11},P_{12})=(6,7,10,10,14)\) and \(V=2/65\). Its connected-curve hypothesis follows directly from the canonical pencil: on a general outer fibre, \(t'=v(t)\) is the ratio of the independent canonical forms \(\alpha|_F\) and \((s/dt)|_F\), hence is a degree-one coordinate on the canonical-pencil base. Thus \(\mathbb C(t,t')\) is the base field of the connected genus-two fibration. \hyperlink{result-6.1}{Lemma 6.1} verifies the quadratic-field hypothesis for all four rows.

For the remaining rows, every primitive seventh denominator has degree at most six. Degree two is excluded by \hyperlink{result-6.2}{Lemma 6.3}. Degree five occurs only in rows 23,24,26,27,29,30,31: use \hyperlink{result-6.4}{Lemma 6.5} for rows 24,26,27, \hyperlink{result-6.3}{Lemma 6.4} for rows 23,29,31, and \hyperlink{result-6.6}{Lemma 6.7} for row 30. Degree six occurs only in rows 29,30,31: use \hyperlink{result-6.5}{Lemma 6.6} for rows 29,31 and \hyperlink{result-6.6}{Lemma 6.7} for row 30. The cubic and quartic alternatives in row 31 are excluded by Lemmas 6.8--6.9.

Corollary~2.2 excludes degree one. There remain two cases: a denominator of degree three or four at volume at most \(8/255\), and a surface image for the sixth map in row 31. We exclude them in \hyperlink{section-7}{Section 7}. Its cubic and quartic theorem applies since every eighth dimension in Table 6.1 is at most seven. Consequently the entire \((1,2)\)-fibre branch is impossible below \(1/30\).

\needspace{4\baselineskip}

\hypertarget{the-remaining-canonical-rings}{%
\needspace{10\baselineskip}\hypertarget{section-7}{}\section{7. The remaining canonical rings}\label{the-remaining-canonical-rings}}

We complete the exclusions required in Sections 5 and 6. We state the arguments in a common form, since several intermediate cases have the same canonical-ring data. Throughout, \(q=0\), \(h^{2,0}=1\), \(R\) is the canonical ring, and \(s\) is a nonzero generator of \(R_1\), and \(E_m=s^{-m}R_m\). A prime denotes the Hodge derivation. We use the Euler identities, the transverse-component criterion, the Jacobian lemma, and the graded multiplicity inequalities recalled in \hyperlink{section-2}{Section 2}.

In the generically finite case put
\[
B=\mathbb C[s,x,y,z_6],\qquad
\deg(s,x,y,z_6)=(1,4,5,6),\qquad
(t,u,z)=(x/s^4,y/s^5,z_6/s^6).
\]
These four homogeneous generators are algebraically independent. The field \(K=\mathbb C(t,u,z)\) has index two in \(\mathbb C(X)\). The deck involution negates both holomorphic forms \(s\) and \(\alpha\): an invariant form would descend holomorphically to a smooth rational model of \(K\), contrary to its vanishing Hodge numbers. Hence the Hodge derivation preserves \(K\). The common low-dimensional data are
\[
(P_1,P_2,P_3,P_4,P_5,P_6,P_7)=(1,1,1,2,3,4,5).
\]
The new seventh section \(W\) is rational over \(B\). Write its primitive relation as \(FW=G\), with \(F,G\in B\) coprime, and call \(\deg F\) its denominator degree. This terminology refers to the homogeneous section, not to a chosen rational coordinate on the quartic-pencil base.

\needspace{3\baselineskip}

\hypertarget{result-7.1}{}\ResultHeading{Theorem 7.1.} Under the preceding assumptions:

\begin{enumerate}
\def\labelenumi{(\alph{enumi})}
\item
  A denominator of degree one is impossible.
\item
  A denominator of degree three or four is impossible if \(P_8\le7\) and \(V\le8/255\).
\end{enumerate}

For the first assertion, the relation is \(sW=x^2\), excluded by Corollary~2.2 without a volume or eighth-plurigenus assumption. In the second, the seventh subring already has a six-dimensional eighth piece, so \(P_8=6\) or \(7\).

The remaining sixth-surface possibility in the last basket is excluded in \hyperlink{section-7.6}{Section 7.5}. In particular, the arguments below complete both cases left by \hyperlink{section-6}{Section 6}.

The local arguments use \cite[Lemma 3.10]{ref1}, while polynomial degree growth
is controlled by \hyperlink{result-2.2}{Proposition 2.4}. These apply to actual
sections and require no finite generation of a differential subring.

\needspace{4\baselineskip}

\hypertarget{normal-forms-in-denominator-degrees-three-and-four}{%
\hypertarget{section-7.2}{}\subsection{7.1. Normal forms in denominator degrees three and four}\label{normal-forms-in-denominator-degrees-three-and-four}}

\hypertarget{lem-seventh-normal-forms}{}\ResultHeading{Lemma 7.2 (Normal forms).} Under the common quadratic-field and low-degree hypotheses, the primitive cubic and quartic seventh relations have the three forms displayed below, and their rings are normal. If \(P_8\le7\), the eighth invariant space adds at most one direction to the old eighth space.

\ProofHeading{Proof.} For the normal forms and normality, no additional bound on \(P_8\) or the volume is needed. After homogeneous triangular changes of generators and subtracting an old seventh section, the possible relations are
\[
\begin{array}{ll}
\text{degree three:}&s^3W=xz_6+by^2,\\
\text{pure degree four:}&s^4W=yz_6+cs^3x^2,\\
\text{finite degree four:}&xW=ayz_6+bsy^2+cs^5z_6+ds^6y+es^{11}\ne0.
\end{array}
\]
Here and below constants in different displays need not be the same.

We verify that this list is exhaustive. A denominator of degree three must be \(s^3\). After removing numerator terms divisible by \(s^3\), a degree-ten numerator is a combination of \(xz_6,y^2,sxy,s^2x^2\). If the \(xz_6\) coefficient vanished, primitivity would make the \(y^2\) coefficient nonzero; its monic quadratic equation would make \(y/s\) an actual fourth section, a contradiction. A sixth-coordinate shift gives the displayed form. A degree-four denominator is \(s^4\) or, after a fourth-coordinate shift, \(x\). In the pure case the primitive degree-eleven numerator has a nonzero \(yz_6\) coefficient. Fifth- and sixth-coordinate shifts remove \(sy^2,sxz_6,s^2xy\), leaving only \(cs^3x^2\) modulo \(s^4B_7\). In the finite case division by \(x\) leaves exactly the five displayed monomials.

We will use the normality of these hypersurfaces. Since they are Cohen--Macaulay, it is enough to check regularity in codimension one. Direct differentiation shows that each singular locus has codimension at least two. For the finite form this also follows from the presentation \(xW-g(s,y,z_6)\): its singular locus is contained in \(x=W=g=0\), together with the derivatives of \(g\). Their multiplicities are respectively \(1/84,11/840,11/840\). In every case the old eighth piece has dimension six. For the remainder of Theorem~7.1(b), choose a new invariant eighth section \(Q\) if one exists, and write \(q=Q/s^8\); otherwise put \(q=0\). Since \(P_8\le7\) and the old eighth piece has dimension six, the invariant eighth space is the old space plus at most this one direction. A possible anti-invariant eighth section plays no role, because \(v(K)\subset K\).

\needspace{4\baselineskip}

\(\square\)

\hypertarget{denominator-degree-three}{%
\hypertarget{section-7.3}{}\subsection{7.2. Denominator degree three}\label{denominator-degree-three}}

\hypertarget{lem-cubic-denominator}{}\ResultHeading{Lemma 7.3 (The cubic denominator).} Under the common quadratic-field and low-degree hypotheses of this section, a primitive seventh denominator of degree three is impossible. This assertion does not require a bound on \(P_8\) or the additional volume cutoff.

\ProofHeading{Proof.} Put
\[
 \mathcal A=\mathbb C[s,x,y,z_6,W]/(s^3W-xz_6-by^2),
 \qquad e(\mathcal A)=1/84.
\]
Corollary~2.2 gives \(t'\notin E_5\), so choose \(z=t'\). The seventh relation becomes
\[
 w=tz+bu^2+ctu+dt^2+\operatorname{affine}(t,u,z).
\]
The first and second Euler identities give
\[
 u'=A_0w+B_0z+C_0u+D_0t+E_0,\qquad
 z'=24\lambda t^2+F_0w+G_0z+H_0u+I_0t+J_0.
\]
At the reduced \(x\)-unit component of Corollary~2.2, \(\ell=1\) and \(t,z\) have exact \(\rho\) poles \(4,6\). If \(b=0\), the term \(tz\) has pole \(9\ell+1\), greater than every other term and than the seventh-section bound \(7\ell\), a contradiction. If \(b\ne0\), cancellation forces \(u\) to have pole \((9\ell+1)/2\). The divergence is
\[
 J'/J=2A_0bu+(A_0c+F_0)t+C_0+G_0.
\]
The two nonconstant terms have distinct actual-time poles greater than one, so \(A_0=F_0=0\). The pole \((11\ell+3)/2\) of \(u'\) then exceeds every remaining term on its right. This proves \(t'\in E_5\). This contradicts Corollary~2.2, so the cubic denominator is impossible. \(\square\)

\hypertarget{a-pure-denominator-of-degree-four}{%
\hypertarget{section-7.4}{}\subsection{7.3. A pure denominator of degree four}\label{a-pure-denominator-of-degree-four}}

\hypertarget{prop-pure-quartic-denominator}{}\ResultHeading{Proposition 7.4 (The pure quartic denominator).} Under the assumptions of Theorem~7.1(b), the primitive denominator \(s^4\) is impossible.

\ProofHeading{Proof.} Write
\[
\mathcal A=
\mathbb C[s,x,y,z_6,W]/(s^4W-yz_6-cs^3x^2),
\qquad w=uz+ct^2,
\qquad e(\mathcal A)=\frac{11}{840}.
\tag{7.1}\label{eq:7.14}
\]

\needspace{4\baselineskip}
\textbf{Eighth sections.}

If \(q\ne0\), a finite pole of its eighth section \(Q\) has residue multiplicity at least \(1/30\). It would give
\[
V\ge2\left(\frac{11}{840}+\frac1{240}\right)
=\frac{29}{840}>\frac8{255}.
\]
Thus only the boundary primes \(T_y=(s,y)\) and \(T_z=(s,z_6)\) can be polar. Their residue multiplicities are \(1/168\) and \(1/140\). The simultaneous bound is
\[
\frac{m_y}{168}+\frac{m_z}{140}\le\frac{37}{1785},
\qquad 5m_y+6m_z\le17.
\tag{7.2}\label{eq:7.15}
\]
In particular \(m_y\le3\), \(m_z\le2\), so \(s^3Q\in\mathcal A_{11}\). The exact degree-eleven enumeration gives
\[
Q=a\frac{xy}{s}+b\frac{y^2}{s^2}
+d\frac{xz_6}{s^2}\pmod{\mathcal A_8}.
\tag{7.3}\label{eq:7.16}
\]
This also holds for \(c=0\), when the boundary values of \(y,z_6\) are four rather than three. The possible term \(xW/s^3\) has pole pair \((3,3)\) and is excluded by \eqref{eq:7.15}; the term \(yz_6/s^3=sW-cx^2\) is old.

If both \(b,d\) are nonzero, their two distinct polar contributions give the pair \((2,2)\), contradicting \eqref{eq:7.15}. Thus \(bd=0\). If \(b\ne0\), normality applied to the resulting monic quadratic makes \(y/s\) an actual fourth section, impossible. Consequently
\[
q=a tu+d tz,
\tag{7.4}\label{eq:7.17}
\]
When no new invariant eighth direction exists, set \(a=d=0\). We retain both this case and the case in which both terms in \eqref{eq:7.17} are present.

\needspace{4\baselineskip}
\textbf{Polynomial degree growth.}

The exact low spaces give
\[
t'=L_0,\qquad u'=L_1+e w,\qquad
z'=L_2+jw+k t^2+h q,
\tag{7.5}\label{eq:7.18}
\]
where the \(L_i\) are affine in \(1,t,u,z\). The ordinary-degree-raising part of the derivation fixes \(t\) and satisfies
\[
D_2u=e(uz+ct^2),\qquad
D_2z=j(uz+ct^2)+kt^2+h(atu+dtz).
\]
Modulo \(t\), this is \(euz\partial_u+juz\partial_z\).

If \(e\ne0\), the highest homogeneous part of \(v^r(w)\) contains \(e^ru z^{r+1}\). Its coefficient cannot cancel: to retain \(u\)-exponent one and \(t\)-exponent zero, every degree-raising step must use \(euz\partial_u\). \hyperlink{result-2.2}{Proposition 2.4}, omitting \(z\), gives \(V\ge1/20\). Similarly \(j\ne0\) gives the monomial \(j^ru^{r+1}z\) and \(V\ge1/24\). Thus \(e=j=0\).

If \(hd\ne0\), the highest homogeneous part of \(v^r(w)\) contains \((hd)^rt^ruz\). \hyperlink{result-2.2}{Proposition 2.4}, now omitting \(t\), gives \(V\ge1/30\). Hence \(hd=0\). Consequently we are left with
\[
t'=L_0,\qquad u'=L_1,\qquad
z'=K t^2+Htu+L_2.
\tag{7.6}\label{eq:7.19}
\]
The quadratic pair \((K,H)\) is nonzero, since otherwise the sixth space would be stable, contrary to the Hodge lemma.

\needspace{4\baselineskip}
\textbf{The remaining differential equations.}

Corollary~2.2 gives \(t'\notin E_5\). Thus \(t'\) is new in degree six. Choose \(z=t'\) and replace \(u\) by \(U=u-at\) to remove the sixth coefficient from \(U'\). Then
\[
U'=bU+ct+d,\qquad
 t''=Kt^2+HtU+L(1,t,t',U).
\tag{7.7}\label{eq:7.21}
\]
At the reduced \(x\)-unit detector of Corollary~2.2, \(\ell=1\), and \(t,t'\) have exact \(\rho\) poles \(4,6\). If \(c\ne0\), the pole of \(U\) is \(3\ell-1\). If \(c=0\), Hodge stability of \(\langle1,U\rangle\) gives \(U'=0\), so \(U\) is independent of actual time.

For \(\ell=1\), balance requires \(K\ne0\). When \(c\ne0\), the pole orders of \((t,U)\) are \((2,1)\) and the resonance indices are \(-1,6,1\), all integral. When \(c=0\), \(U\) is already time-constant and the second-order equation has the integral indices \(-1,6\).

No positive nonintegral rational exponent is possible. Thus \(t,U,t'\) have integral actual-time expansions and generate the sixth field. The Jacobian lemma gives a contradiction; its Jacobian clause is essential, since coordinate integrality alone would allow ramification index two. This excludes the pure quartic denominator.

\needspace{4\baselineskip}

\(\square\)

\hypertarget{a-finite-denominator-of-degree-four}{%
\hypertarget{section-7.5}{}\subsection{7.4. A finite denominator of degree four}\label{a-finite-denominator-of-degree-four}}

The primitive finite denominator can be shifted to \(x\). Its normal form is
\[
\mathcal A=
\mathbb C[s,x,y,z_6,W]/(xW-g(s,y,z_6)),
\qquad e(\mathcal A)=\frac{11}{840},
\]
\[
g=a yz_6+bsy^2+cs^5z_6+ds^6y+es^{11}\ne0.
\tag{7.8}\label{eq:7.23}
\]
As noted in \hyperlink{section-7.2}{Section 7.1}, this ring is normal. We will bound the poles of the eighth sections and of the actual ninth derivative, then compare their differential identities.

\hypertarget{lem-eighth-ninth-poles}{}\ResultHeading{Lemma 7.5 (Eighth and ninth pole filtrations).} For the finite quartic normal form \eqref{eq:7.23}, an invariant eighth section and the ninth section \(s^9w'\) have no finite poles. Their boundary poles satisfy \eqref{eq:7.27}, and their complete spaces are those in \eqref{eq:7.26}, \eqref{eq:7.28}, and \eqref{eq:7.29}.

\ProofHeading{Proof.} At a finite prime where \(x\) is a unit, the polynomial localization gives a residue polynomial subring of weight product at most \(30\). At a finite prime containing \(x\), the relation imposes an irreducible factor of \(g\) in the variables \(s,y,z_6\), while \(W\) remains free. Since \(s\) is a unit, the residue contains algebraically independent images of \(s,W\) and either \(y\) or \(z_6\), with weight product at most \(42\).

Consequently no invariant eighth section has a finite pole, since
\[
2\left(\frac{11}{840}+\frac1{8\cdot42}\right)
=\frac9{280}>\frac8{255}.
\tag{7.9}\label{eq:7.24}
\]
Since the derivation preserves \(K\), the function \(w'\) is invariant. The same finite-pole bound holds for the particular ninth section \(D_9=s^9w'\):
\[
2\left(\frac{11}{840}+\frac1{9\cdot42}\right)
=\frac{17}{540}>\frac8{255},
\qquad \frac{17}{540}-\frac8{255}=\frac1{9180}.
\tag{7.10}\label{eq:7.25}
\]
No assertion that all ninth sections are invariant is needed.

If there is no new invariant eighth direction, put \(q=0\) throughout the following classification. If \(a\ne0\), the boundary \(s=0\) is irreducible, of residue multiplicity \(11/840\). Its eighth and ninth pole bounds are both one; the ninth bound before taking integers is \(333/187<2\). Consequently, if \(q\ne0\) then \(Q=xy/s\) modulo \(\mathcal A_8\); in all cases
\[
w'\in N_1:=\langle1,t,u,z,w,t^2,tu,tz,u^2\rangle.
\tag{7.11}\label{eq:7.26}
\]

If \(a=0\), write \(g=s^eh\) with \(s\nmid h\) and \(e\in\{1,5,6,11\}\). The boundary primes \((s,x)\) and \((s,W)\) have residue multiplicities \(1/210\) and \(1/120\). Their valuations are
\[
\nu_x(s,x,y,z_6,W)=(1,e,0,0,0),\qquad
\nu_W(s,x,y,z_6,W)=(1,0,0,0,e).
\]
The eighth and ninth simultaneous bounds become
\[
4m_x(Q)+7m_W(Q)\le17,
\qquad
4m_x(D_9)+7m_W(D_9)\le19.
\tag{7.12}\label{eq:7.27}
\]
The unrounded right sides are \(296/17\) and \(333/17\).

An exact degree-eight enumeration gives either
\[
q\in\langle tu,tz\rangle\quad\text{or}\quad q=uw.
\tag{7.13}\label{eq:7.28}
\]
The latter requires \(e\ge5\), and cannot be mixed with a nonzero polynomial term: its pole pair \((4,0)\) together with those of \(tu,tz\) exceeds the eighth bound. For \(e=1\), \(tz\) is also excluded, leaving \(q=tu\) or \(q=0\).

The ninth enumeration gives
\[
w'\in N_1\quad(e=1),\qquad
w'\in N_{\ge5}:=N_1+\langle uw,zw\rangle\quad(e\ge5).
\tag{7.14}\label{eq:7.29}
\]
For completeness, we verify these enumerations directly. We use normal monomials \(s^\alpha x^iy^jz_6^kW^l\) with \(il=0\), and set
\[
\alpha=m-4i-5j-6k-7l,
\quad m_x=\max(0,-\alpha-ei),
\quad m_W=\max(0,-\alpha-el).
\]
For \(m=8\), applying the first bound in \eqref{eq:7.27}, and removing old sections, gives \eqref{eq:7.28}. For \(m=9\), the extra monomials are exactly \(y^2/s,xz_6/s\), and, when \(e\ge5\), \(yW/s^3,z_6W/s^4\). They normalize to \(u^2,tz,uw,zw\). At \((s,x)\), powers of the free residue coordinate \(W\) separate leading terms, and at \((s,W)\) powers of \(x\) do so. Within one power, the remaining residue variables are independent. The relation has already removed all simultaneous \(xW\) factors, so boundary cancellation introduces no extra possibility.

\(\square\)

\hypertarget{prop-finite-quartic-denominator}{}\ResultHeading{Proposition 7.6 (The finite quartic denominator).} Under the assumptions of Theorem~7.1(b), the primitive denominator \(x\) is impossible.

\ProofHeading{Proof.} Suppose \(a\ne0\). Legal fifth- and sixth-coordinate shifts and rescaling give
\[
tw=uz+\varepsilon.
\tag{7.15}\label{eq:7.30}
\]
By \eqref{eq:7.26}, \(q\) is either \(tu\) or zero, and \(w'\in N_1\). Differentiating \eqref{eq:7.30} and comparing rational functions in the independent variables \(t,u,z\) gives the exact system
\[
\begin{aligned}
t'&=l t+B u+A z,\\
u'&=c t+d u+A w,\\
z'&=e t+f z+B w+Kt^2+Htu,\\
w'&=c z+e u+(d+f-l)w+Ktu+Hu^2,
\end{aligned}
\tag{7.16}\label{eq:7.31}
\]
with \((d+f)\varepsilon=0\). The comparison uses the precise space \(N_1\): after substituting \(w=(uz+\varepsilon)/t\), that space has no double pole in \(t\). It is stronger than mere regularity on \(\mathcal A_s\). Direct differentiation verifies that \eqref{eq:7.31} preserves \eqref{eq:7.30}.

If \(A=0\), the span of \(1,t,u\) is stable, a contradiction. Hence \(A\ne0\). Choose the reduced \(x\)-unit detector of Corollary~2.2. Put \(r=u/t\). Then
\[
r'=F(r)+A\varepsilon/t^2,
\qquad F(r)=c+(d-l)r-Br^2,
\]
\[
t''=A(K+Hr)t^2+(l+f+2Br)t'+G(r)t+2AB\varepsilon/t,
\tag{7.17}\label{eq:7.32}
\]
where
\[
G(r)=Bc+Ae-fl+B(d-f-2l)r-2B^2r^2.
\]
The pair \((K,H)\) is nonzero, since otherwise \(1,t,u,z,w\) spans a stable space.

At this detector, \(t\) has actual-time pole \(p=2\), while a pole of \(r\) has order less than one. Such a positive pole is impossible in the Riccati equation: its derivative would dominate its square and every other term. Thus \(r\) is bounded, with value \(r_0\).

Balance in the second equation forces \(K+Hr_0\ne0\): otherwise its right side has pole less than four, while \(t''\) has pole four.

For \(p=2\), a first nonintegral relative exponent \(q>0\) of \(t\) has indicial polynomial
\[
(q-2)(q-3)-12=(q+1)(q-6).
\]
The Riccati equation introduces no earlier fractional correction: a relative \(q\) term of \(t\) first induces a fractional term of \(r\) at absolute exponent \(2p+q+1\), too late to affect these indicial equations. Thus \(t,r,t'\), and hence \(t,u,z\), have integral actual-time expansions. The Jacobian lemma, including its Jacobian clause for \(\ell=1\), gives the contradiction.

\needspace{4\baselineskip}
\textbf{Boundary order one.}

Suppose the \(yz_6\) coefficient vanishes and the \(sy^2\) coefficient is nonzero. Complete the square in the fifth coordinate. The two possible forms are
\[
tw=u^2+\varepsilon\qquad\text{or}\qquad tw=u^2+z.
\tag{7.18}\label{eq:7.34}
\]
In the second form, a legal sixth shift absorbs the constant and rescales its nonzero coefficient. In both cases \(q\) is either \(tu\) or zero, and \(w'\in N_1\).

Write \(u'=Ew+L_1\), and, for the second form, \(z'=Jw+Kt^2+Htu+L_2\). The absence of a double pole in \(w'\) gives, respectively,
\[
t'=lt+2Eu\qquad\text{or}\qquad t'=lt+2Eu+J.
\tag{7.19}\label{eq:7.35}
\]
Both alternatives place \(t'\) in \(E_5\), contrary to Corollary~2.2.

\needspace{4\baselineskip}
\textbf{Higher boundary order and a polynomial eighth function.}

The remaining numerator is linear:
\[
tw=h(u,z),\qquad h=az+bu+c\ne0.
\tag{7.20}\label{eq:7.36}
\]
First suppose \(q\in\langle tu,tz\rangle\), including \(q=0\). Write
\[
u'=Ew+L_1,\qquad z'=Jw+Kt^2+Hq+L_2.
\]
Differentiating \eqref{eq:7.36} gives \(t'w+tw'=az'+bu'\). The allowed space \(N_{\ge5}\) has no double pole in \(t\). Its coefficient forces
\[
t'|_{t=0}=aJ+bE.
\]
Thus \(t'\) is affine in \(t\), contradicting independence of \(t,t'\).

\needspace{4\baselineskip}
\textbf{Higher boundary order and the eighth function \(uw\).}

The only remaining eighth function is \(q=uw\). Double-pole comparison now gives
\[
t'=lt+aHu+aJ+bE.
\]
This also places \(t'\) in \(E_5\), contrary to Corollary~2.2.

\(\square\)

This excludes every finite quartic denominator. Together with Sections 7.2 and 7.3, it proves \hyperlink{result-7.1}{Theorem 7.1}(b).

\needspace{4\baselineskip}

\hypertarget{the-last-sixth-surface-image}{%
\hypertarget{section-7.6}{}\subsection{7.5. The last sixth-surface image}\label{the-last-sixth-surface-image}}

The remaining basket has
\[
\mathcal B=4(1,2)+(2,5)+(4,9)+(4,11),\qquad
V=\frac{16}{495},\qquad(P_5,P_6,P_7)=(3,4,5).
\]
It remains to show that the sixth canonical image of this basket cannot be a surface. The smooth-quadric alternative with its specified outer ruling was excluded in \hyperlink{section-6}{Section 6}. The remaining surface alternatives satisfy \(V\ge\xi/12\), by \hyperlink{section-6}{Section 6}. The one-point marking condition, together with \(\xi\ge5/14\) and \(\xi<2/5\), therefore gives
\[
 a=\xi=\frac4{11}.
\tag{7.21}\label{eq:7.40}
\]
The marking lies over the singleton \((4,11)\) constituent. Use the actual canonical divisors from the geometric reduction: \(Q=\pi^*K_Y\), \(D=Q|_F\), \(A=\sigma^*K_S\), and \(\Gamma\) is the unique component of the selected actual \(K_F\) with positive \(A\)-degree. In particular
\[
 A\Gamma=1,\quad \operatorname{coeff}_\Gamma K_F=1,
 \quad\operatorname{coeff}_\Gamma D=4/11.
\]

\needspace{3\baselineskip}

\hypertarget{result-7.2}{}\ResultHeading{Lemma 7.7.} In the quartic-pencil setting, if \(V<1/30\), then
\[
 E_6\cap\mathbb C(t)=\langle1,t\rangle_{\mathbb C},
 \qquad t=x_4/s^4.
\tag{7.22}\label{eq:7.41}
\]

\ProofHeading{Proof.} Suppose an additional rational function \(h(t)\) occurs, so that \(H_6=s^6h(t)\) is an actual sixth section.

First suppose \(h\) has a finite pole of order \(m\ge1\) at \(t=\lambda\). On a common resolution of the outer pencil let \(F_\infty\) and \(F_\lambda\) be its full fibre divisors. Use the effective canonical representative \(Q=\pi^*\operatorname{div}_Y(s)\). The actual section \(x_4=s^4t\) gives coefficientwise \(Q\ge F_\infty/4\) on the support of \(F_\infty\). The section \(H_6\) gives \(Q\ge mF_\lambda/6\) on the support of \(F_\lambda\). These supports are disjoint. Consequently
\[
 Q\ge\frac14F_\infty+\frac m6F_\lambda
 \ge\frac14F_\infty+\frac16F_\lambda.
\tag{7.23}\label{eq:7.42}
\]
In particular \(12Q\ge3F_\infty+2F_\lambda\sim5F\). The canonical restriction inequality gives
\[
 D\ge_{\mathbb Q}\frac5{17}A.
\]
As \(D\) and \(A\) are nef and \(A^2=1\), we obtain
\[
 V\ge\frac5{12}D^2
 \ge\frac5{12}\left(\frac5{17}\right)^2
 =\frac{125}{3468}>\frac1{30}.
\tag{7.24}\label{eq:7.43}
\]
Here \(125/3468-1/30=47/17340\). This argument uses the effective pullback \(Q\), not a new canonical divisor after a blowup.

It remains to consider a polynomial \(h(t)=\sum_{i=0}^d c_it^i\) of degree \(d\ge2\), scaled so that \(c_d=1\). Put \(T=x_4/s\), a homogeneous rational element of degree three. The actual section \(H_6\) gives
\[
 s^{3d-6}H_6
 =\sum_{i=0}^d c_iT^i s^{3(d-i)}.
\tag{7.25}\label{eq:7.44}
\]
This is a monic equation for \(T\) over the canonical ring, since \(3d-6\ge0\). Normality gives \(T\in R_3\). But \(P_3=1\), so \(x_4\) would be a multiple of \(s^4\), contrary to \(P_4=2\). This proves \eqref{eq:7.41}. \(\square\)

\hypertarget{prop-sixth-surface-exclusion}{}\ResultHeading{Proposition 7.8 (The remaining sixth-surface case).} The sixth-surface possibility in the final basket of Section 6 cannot occur.

\ProofHeading{Proof.} The rank-two alternative with the same outer pencil from \hyperlink{section-6}{Section 6} supplies a third outer function. Indeed, at a general \(F_\lambda\), its two-dimensional restriction kernel contains \(s^2(x-\lambda s^4)\). Choose an independent kernel section \(H\). In the same-outer-pencil alternative their ratio belongs to \(\mathbb C(t)\), hence \(H/s^6\in\mathbb C(t)\). It cannot lie in \(\langle1,t\rangle\): the kernel of evaluation at \(t=\lambda\) in that space is only one-dimensional. Thus this alternative contradicts \eqref{eq:7.41}. The other rank-two alternative is the already excluded smooth quadric. Thus it suffices to treat the sixth restriction of rank three.

Write \(u=y_5/s^5\) for the new fifth function. By \hyperlink{result-7.2}{Lemma 7.7} it is not an outer function. On the generic outer fibre, the canonical pole filtration and \(\lfloor5a\rfloor=1\) show that \(u\) is a degree-one coordinate on the canonical-pencil base, with a simple pole at its distinguished canonical member. A new sixth function \(z=z_6/s^6\) has an exact double pole, since the sixth restriction has rank three and \(\lfloor6a\rfloor=2\).

Here is the pole-filtration justification. Along the horizontal divisor \(E\) containing \(\Gamma\), the selected canonical divisor has multiplicity one. Since \(\alpha\) is holomorphic, the coefficients of \(v=\alpha/s\) have pole order at most one in smooth local coordinates. The function \(t\) is regular at the generic point of this horizontal divisor, so \(t'\) has pole order at most one there. The only new sixth direction has exact pole order two. Since \(t'\in E_6\), its coefficient in that direction vanishes, and \(t'\in E_5\). This contradicts Corollary~2.2 and excludes the remaining sixth-surface case. \(\square\)

\hypertarget{proof-of-the-volume-theorem-and-its-equality-case}{%
\needspace{10\baselineskip}\hypertarget{section-8}{}\section{8. Proof of the volume theorem and its equality
case}\label{proof-of-the-volume-theorem-and-its-equality-case}}

\needspace{4\baselineskip}

\hypertarget{completion-of-the-volume-bound}{%
\hypertarget{section-8.1}{}\subsection{8.1. Completion of the volume
bound}\label{completion-of-the-volume-bound}}

\ProofHeading{Proof of \hyperlink{result-1.1}{Theorem 1.1}, except for strictness.} We pass to a terminal
minimal model and suppose, for a contradiction, that \(V<1/30\). \hyperlink{section-3}{Section 3} gives
\(p_g=P_2=P_3=P_4-1=1\), \(q=0\), and a connected quartic pencil with
rational base. Its general minimal fibre has type \((1,2,0)\),
\((2,2,0)\), or \((2,3,0)\). \hyperlink{section-4}{Section 4} excludes both Euler
characteristics in the last type. \hyperlink{section-5}{Section 5} excludes type \((2,2,0)\).
For type \((1,2,0)\), \hyperlink{section-6}{Section 6} gives the complete necessary basket list
and reduces it to the denominator and surface cases excluded in \hyperlink{section-7}{Section
7}. Thus \(V\ge1/30\). \(\square\)

We now exclude \(V=1/30\) when \(\chi>0\). We retain the notation
and the Hodge, residue, and numerical tools of \hyperlink{section-2}{Section 2}. In particular
all divisor coefficients refer to the actual section \(s\). No
classification of equality cases is assumed.

\needspace{4\baselineskip}

\hypertarget{equality-in-surface-restriction}{%
\hypertarget{section-8.2}{}\subsection{8.2. Equality in surface
restriction}\label{equality-in-surface-restriction}}

We first record the consequences of equality in canonical restriction.
For a coefficient-one pencil \(mK_Y\sim F_Y+Z_Y\), we use the special model
of \cite[Proposition 2.2]{ref1}, and put \(D=\mu^*K_Y|_F\), \(A=\sigma^*K_{F_0}\).
Its effective adjunction formula is \[
(m+1)D\sim_{\mathbb Q}K_F+H,\qquad
H=G|_F+\mu^*Z_Y|_F\ge0.
\tag{8.1}\label{eq:8.1}
\] We apply this formula on the special model itself; after a further
blowup we pull back both sides.

\needspace{3\baselineskip}

\hypertarget{result-8.1}{}\ResultHeading{Lemma 8.1 (maximal restriction).} If \(A^2=m(m+1)^2V\), then
\(Z_Y=0\), the general member \(F_Y\) has Du Val singularities, and \[
\operatorname{ind}_P(K_Y)\mid m+1\quad(P\in F_Y),\qquad
\operatorname{ind}_P(K_Y)\mid m\quad(P\notin F_Y).
\tag{8.2}\label{eq:8.2}
\] In particular, the global terminal index divides \(m(m+1)\).

\ProofHeading{Proof.} We have equality throughout \[
mV\ge D^2\ge DA/(m+1)\ge A^2/(m+1)^2.
\] Hence the Hodge index theorem gives \(D\equiv A/(m+1)\). The effective divisor
\(K_F-A+H\) is therefore numerically trivial, so it vanishes. Thus \(F\) is minimal
and \(H=0\). In particular \(\mu^*Z_Y|_F=0\). On a general very ample
smooth surface in \(Y\), the restrictions of \(F_Y\) and \(Z_Y\) are
orthogonal, and their sum is nef and big. Their squares are nonnegative;
Hodge index forces \(Z_Y=0\). This is the argument of \cite[Lemma 3.5]{ref9}, applied to \(mK_Y\).

Since \((F,G|_F)\) is klt, \cite[Lemma 3.4(1)]{ref9} gives a plt pair
\((Y,F_Y)\) and \(K_F=\mu_F^*K_{F_Y}\). Thus \(F_Y\) is Du Val.
Adjunction gives \(K_{F_Y}\sim(m+1)K_Y|_{F_Y}\), and restriction to the
integral boundary of a plt pair preserves the local Cartier index \cite[Appendix A, Theorem A.1]{int28}. This proves the first assertion in \eqref{eq:8.2};
off \(F_Y\), the relation \(mK_Y\sim F_Y\) proves the second.
\(\square\)

\needspace{3\baselineskip}

\hypertarget{result-8.2}{}\ResultHeading{Lemma 8.2 (the square-three endpoint).} Suppose a connected
quartic pencil, possibly a general subpencil of a net, satisfies \[
A^2=3,\qquad AD=2/3,\qquad D^2=2/15.
\tag{8.3}\label{eq:8.3}
\] Then \(Z_Y=0\). The general member \(F_Y\) is normal klt, has exactly
one non-Gorenstein singularity, of type \(\frac13(1,1)\), and otherwise
has Du Val singularities. Moreover \[
\operatorname{ind}_P(K_Y)\mid15\quad(P\in F_Y),\qquad
\operatorname{ind}_P(K_Y)\mid4\quad(P\notin F_Y).
\tag{8.4}\label{eq:8.4}
\]

\ProofHeading{Proof.} We intersect \eqref{eq:8.1} with \(D\) to obtain \(D(K_F-A+H)=0\). Every
\(\sigma\)-exceptional curve occurs in \(K_F-A\), so nefness makes \(D\)
descend to a nef big divisor \(D_0\) on \(F_0\). With \(T=\sigma_*H\),
\[
5D_0\sim_{\mathbb Q}K_{F_0}+T,\qquad D_0T=0,\qquad K_{F_0}T=1/3.
\] The support of \(T\) is negative definite. If a component \(C\) has
canonical degree \(k>0\) and coefficient \(t_C\), adjunction and the
nonnegative intersections with other components give \[
0=(K_{F_0}+T)C\ge k-t_C(k+2),\qquad
kt_C\ge k^2/(k+2)\ge1/3.
\] Thus all these inequalities are equalities. We conclude that \(T=\Gamma/3\), where
\(\Gamma\) is a smooth rational \((-3)\)-curve of canonical degree one,
disjoint from every other component. Any remaining components have
degree zero and form negative definite \((-2)\)-configurations; their
homogeneous intersection equations force their coefficients to vanish.

We next compare the two surface models. The exceptional divisor
\(K_F+H-\sigma^*(K_{F_0}+\Gamma/3)\) is
numerically trivial and therefore zero. A first blowup of
\((F_0,\Gamma/3)\) has negative crepant boundary coefficient,
incompatible with the effective pushforward of \(H\) on each
intermediate surface. Thus \(F=F_0\) and \(H=\Gamma/3\).

A strict fixed component meeting \(F\), or an exceptional divisor with
curve centre meeting \(F\), would contribute an integral positive
coefficient to \(H\). Therefore every component of \(H\) maps to a point
on \(Y\). The effective divisor \(\mu^*Z_Y|_F\) is supported there and
has degree zero on all contracted curves; negative definiteness gives
\(\mu^*Z_Y|_F=0\). The ample-surface argument in \hyperlink{result-8.1}{Lemma 8.1} gives
\(Z_Y=0\). Now \cite[Lemma 3.4(1)]{ref9} gives a plt pair and \[
K_F+\Gamma/3=\mu_F^*K_{F_Y}.
\] The other exceptional curves are \((-2)\)-curves disjoint from
\(\Gamma\). Their contractions are Du Val, while contracting \(\Gamma\)
gives \(\frac13(1,1)\). Thus \(3K_{F_Y}\) is Cartier; equivalently one
may use the descent criterion \cite[Theorem (12.1)(i)]{int41}. Since
\(K_{F_Y}\sim5K_Y|_{F_Y}\), index restriction gives \eqref{eq:8.4}. \(\square\)

\needspace{4\baselineskip}

\hypertarget{the-endpoint-numerical-reduction}{%
\hypertarget{section-8.3}{}\subsection{8.3. The endpoint numerical
reduction}\label{the-endpoint-numerical-reduction}}

\hypertarget{prop-positive-euler-endpoint}{}\ResultHeading{Proposition 8.3 (The positive-Euler endpoint).} If \(p_g>0\), \(V=1/30\), and \(\chi>0\), then \(q=0\), \(p_g=P_2=1\), \(h=\chi=1\), \(P_4=3\), and the basket is one of the ten listed in Section 8.4.

\ProofHeading{Proof.} We assume throughout this subsection that \(p_g>0\), \(V=1/30\),
and \(h=\chi>0\). The published bounds in
\hyperlink{section-3.1}{Section 3.1} and \(P_5\ge5\Rightarrow V\ge1/25\) \cite[Theorem 2.5(ii)--(iii)]{ref7} give \[
q=0,\quad p_g=P_2=1,\quad P_3\le2,\quad P_4\le3,\quad P_5\le4,
\qquad h=h^{2,0}.
\tag{8.5}\label{eq:8.5}
\] By Proposition~2.1, \(h\le P_4-P_2\le2\). When \(P_3=2\), the first part of \hyperlink{section-3.1}{Section 3.1} remains valid at equality: the tricanonical pencil has moving degree one and regular minimal fibre \((1,2,0)\). The Hodge restriction bound \(h\le p_g(F)-1=1\) therefore gives \(h=1\).

The geometric part of \hyperlink{section-3.2}{Section 3.2} also remains valid with \(V\le1/30\). If \(P_4=3\), a pencil would have moving degree at least two and give \(V\ge1/18\). Otherwise the fourth map is a connected genus-two fibration over \(\mathbb P^2\). Its general fourth surface has minimal type \((3,3,0)\): the inequalities
\(\frac54AN\le A^2\le10/3\) still force \(A^2=3\), \(N\equiv C\), and \(AC=2\). The trace and adjunction argument gives \(p_g(S)\ge2+h\), whereas Noether's inequality gives \(p_g(S)\le3\); hence again \(h=1\). If \(h=2\), these two reductions would force \(P_3=1\), \(P_4\le2\), contradicting \(h\le P_4-P_2\). Thus \(h=1\) in every case.

We apply the recursion of Appendix B at the single volume \(1/30\). The filters are \(P_5\le4\), monotonicity, the saturated-domain inequality of \cite[Lemma 3.13]{ref1}, and the marked-curve test below. Their exact counts are
\[
\begin{array}{c|r|r|r|r|r}
(P_3,P_4)&\text{index bound}&\text{raw}&\text{low degrees}&\text{domain}&\text{marked curve}\\\hline
(1,2)&33&4&3&0&0\\
(1,3)&32&11&9&3&3\\
(2,2)&34&6&5&2&0\\
(2,3)&33&24&12&7&7.
\end{array}
\]
The case \((P_3,P_4)=(1,1)\) is already excluded by Proposition~2.1 (and has negative gap in \eqref{eq:B.1}). For the domain filter, when \(r_Y=30\), \cite[Lemma 3.13]{ref1} makes both \(\mathcal A\) and \(\mathcal S\) domains. We apply the product inequality to
\[
 a_m=\dim\mathcal S_m=P_m-P_{m-1}+\mathbf1_{m=2}.
\]

The two remaining \(P_3=P_4=2\) baskets are \[
4(1,2)+(1,3)+(1,4)+(7,20),\qquad (1,2)+(5,12)+(7,20).
\] On the tricanonical \((1,2)\) fibre, \cite[Lemmas 2.5--2.6 and Theorem 1.1]{ref1} give \(5/14\le\xi\le12V=2/5\). The local index must clear the denominator of \(\xi\); this leaves only \(\xi=2/5\) at the unique \((7,20)\) constituent. This is an
actual quotient point: a nonquotient point of this index is \(cA/20\),
whose axial weight counts equal constituents \cite[Remark 2.1]{refCH11}. On
the smooth index-one cover, the canonical coefficient has character
\(13\) modulo \(20\) and multiplicity at most two. No monomial of degree
at most two in characters \((1,19,7)\) has character \(13\). Thus neither row can occur.

The ten remaining baskets all have \(h=1\) and \(P_4=3\); they are listed in \hyperlink{section-8.4}{Section 8.4}. Each has global index divisible by thirty and greater than thirty. Appendix B specifies the finite computation. Passing these necessary tests does not assert geometric realizability. \(\square\)

\hypertarget{euler-characteristic-one}{%
\hypertarget{section-8.4}{}\subsection{8.4. Euler characteristic
one}\label{euler-characteristic-one}}

\hypertarget{prop-positive-euler-exclusion}{}\ResultHeading{Proposition 8.4 (Strictness for positive Euler characteristic).} None of the ten baskets of Proposition~8.3 is geometrically realizable. Consequently \(\chi>0\) implies \(V>1/30\).

\ProofHeading{Proof.} By Proposition~8.3, \(h=1\) and \(P_4=3\). The fourth map has the
geometry described in \hyperlink{section-8.3}{Section
8.3}. The ten remaining baskets are as follows; the triple column is
\((P_5,P_6,P_7)\).

\begin{longtable}[]{@{}
  >{\raggedright\arraybackslash}p{(\columnwidth - 6\tabcolsep) * \real{0.0550}}
  >{\raggedright\arraybackslash}p{(\columnwidth - 6\tabcolsep) * \real{0.6100}}
  >{\raggedleft\arraybackslash}p{(\columnwidth - 6\tabcolsep) * \real{0.0850}}
  >{\raggedright\arraybackslash}p{(\columnwidth - 6\tabcolsep) * \real{0.2500}}@{}}
\toprule\noalign{}
\begin{minipage}[b]{\linewidth}\raggedright
Row
\end{minipage} & \begin{minipage}[b]{\linewidth}\raggedright
Basket
\end{minipage} & \begin{minipage}[b]{\linewidth}\raggedleft
\(P_3\)
\end{minipage} & \begin{minipage}[b]{\linewidth}\raggedright
Triple
\end{minipage} \\
\midrule\noalign{}
\endhead
\bottomrule\noalign{}
\endlastfoot
A1 & \(7(1,2)+(1,4)+(1,5)+(5,12)\) & 1 & \((3,4,5)\) \\
A2 & \((1,2)+(1,5)+(2,7)+(10,21)\) & 1 & \((3,5,6)\) \\
A3 & \((1,2)+(1,5)+3(4,9)\) & 1 & \((4,5,7)\) \\
A4 & \(6(1,2)+4(1,3)+2(1,4)+(1,5)\) & 2 & \((3,5,5)\) \\
A5 & \(5(1,2)+(1,3)+2(1,4)+(2,5)+(3,10)\) & 2 & \((3,5,6)\) \\
A6 & \(4(1,2)+(1,3)+(1,4)+(1,5)+2(3,8)\) & 2 & \((4,5,6)\) \\
A7 & \(3(1,2)+3(1,3)+(1,4)+(1,5)+(5,12)\) & 2 & \((4,5,6)\) \\
A8 & \(2(1,2)+(1,3)+2(1,4)+4(2,5)\) & 2 & \((4,4,6)\) \\
A9 & \(2(1,2)+(1,4)+(2,5)+(3,10)+(5,12)\) & 2 & \((4,5,7)\) \\
A10 & \((1,2)+(1,3)+(1,5)+(4,13)+(6,13)\) & 2 & \((4,6,7)\) \\
\end{longtable}

By \eqref{eq:whole-member-trace}, the coefficient-one moving component
of a general fourth-canonical member satisfies
\[
P_{10}-P_6\ge P_2(F)=3+\chi(\mathcal O_F)=7.
\]
For A3 and A4, Riemann--Roch gives \(P_6=5\) and \(P_{10}=11\).
Thus these two cases are impossible, and we consider the remaining
eight baskets.

The marked divisor construction in \hyperlink{section-3.2}{Section 3.2} gives \[
\operatorname{ord}_T(s)=1,\quad
\operatorname{ord}_T\operatorname{Fix}|4K_W|=2,\quad
\dim\varphi_4(T)=1,\qquad 3/5<a:=\operatorname{coeff}_TQ=AD\le2/3.
\tag{8.6}\label{eq:8.7}
\] At equality, only the upper inequality changes. For clarity, the
coefficient-two assertion needs just \(4B=F_0+Z_0\),
\(\operatorname{Supp}F_0\subset\operatorname{Supp}B\), and
\(K_{S_0}F_0=2\): the unique component of \(B\) of positive canonical
degree has coefficient and degree one, hence coefficient two in \(F_0\).
No additional assertion about the entire half-fibre is used here.

To test a possible marked point, we use the full index-cover characters
\((1,-1,b,0)\), omitting the last coordinate for a smooth cover. Write
\(s=f\nu\) with \(f\) of character \(-b\). If \(a=p/e\) is reduced, the
normalized pullback over \(T\) has inertia order \(e\) and \(v(f)=p\).
If inertia acts on a uniformizer with character \(c\), then
\(cp\equiv-b\pmod e\), and the coordinate lower valuations are the least
positive residues of \(c^{-1}(1,-1,b,0)\). A quartic coefficient
attaining the fixed order two has valuation \(4p-2e\) and full character
\(-4b\) modulo the actual index \(r\). Equivariant power-series lifts
give the same test for \(cA/r\) points. A one-dimensional image of \(T\)
requires at least two independent initial monomials.

For A1 and A2 all possible marked data are

\begin{longtable}[]{@{}
  >{\raggedright\arraybackslash}p{(\columnwidth - 12\tabcolsep) * \real{0.0500}}
  >{\raggedright\arraybackslash}p{(\columnwidth - 12\tabcolsep) * \real{0.1800}}
  >{\raggedright\arraybackslash}p{(\columnwidth - 12\tabcolsep) * \real{0.2000}}
  >{\raggedright\arraybackslash}p{(\columnwidth - 12\tabcolsep) * \real{0.1200}}
  >{\raggedleft\arraybackslash}p{(\columnwidth - 12\tabcolsep) * \real{0.0700}}
  >{\raggedright\arraybackslash}p{(\columnwidth - 12\tabcolsep) * \real{0.1800}}
  >{\raggedright\arraybackslash}p{(\columnwidth - 12\tabcolsep) * \real{0.2000}}@{}}
\toprule\noalign{}
\begin{minipage}[b]{\linewidth}\raggedright
Row
\end{minipage} & \begin{minipage}[b]{\linewidth}\raggedright
\((b,r)\); \(a\)
\end{minipage} & \begin{minipage}[b]{\linewidth}\raggedright
Lower valuations
\end{minipage} & \begin{minipage}[b]{\linewidth}\raggedright
Coefficient tested
\end{minipage} & \begin{minipage}[b]{\linewidth}\raggedleft
Target
\end{minipage} & \begin{minipage}[b]{\linewidth}\raggedright
Full character
\end{minipage} & \begin{minipage}[b]{\linewidth}\raggedright
Initial monomials
\end{minipage} \\
\midrule\noalign{}
\endhead
\bottomrule\noalign{}
\endlastfoot
A1 & \((5,12)\); \(2/3\) & \((2,1,1,3)\) & canonical & 2 & \(7\) mod
\(12\) & none \\
A2 & \((10,21)\); \(13/21\) & \((5,16,8,21)\) & quartic & 10 & \(2\) mod
\(21\) & \(x^2\) \\
A2 & \((10,21)\); \(2/3\) & \((1,2,1,3)\) & quartic & 2 & \(2\) mod
\(21\) & \(x^2\) \\
\end{longtable}

In the table we list every monomial whose lower weight is at most the
target and whose full character is the prescribed one. Its first row forbids \(s\) itself;
each other row gives an initial space of dimension at most one. Thus
A1 and A2 are impossible.

For the remaining rows we have \(P_3=2\). We write \(R_4=sR_3+\mathbb Cz_4\) and
\(d=1-a\ge1/3\). A tricanonical section has order \(v_T(g)+3d\). If
\(d>1/3\), this integer is at least two; then \(sR_3\) has order at
least three and only \(z_4\) contributes to the initial quartic system,
contrary to \eqref{eq:8.7}. Hence \(a=2/3\). Some tricanonical coefficient is a
unit, so its character is trivial and the actual local index divides
three. It equals three, excluding A9.

We therefore have \(2/3=AD\le5D^2\le20V=2/3\). \hyperlink{result-8.2}{Lemma 8.2} applies to a general
subpencil. Its pointwise index alternative excludes A5, A6, A7, and A10,
because of the denominators \(10,8,12,13\). For A8, the plt
residue sequence \[
0\longrightarrow\omega_Y\longrightarrow\mathcal O_Y(5K_Y)
\longrightarrow\omega_{F_Y}\longrightarrow0
\] is exact: take adjunction on the index-one cover and invariants.
Vanishing and \(h^1(\omega_Y)=1\) give \(p_g(F_Y)=P_5\). Rational
singularities give \(p_g(F_Y)=3\), excluding A8.

Together with the trace exclusions, this excludes all ten rows.

\needspace{4\baselineskip}

\(\square\)

\hypertarget{canonical-models-and-moduli-at-volume-one-thirtieth}{%
\needspace{10\baselineskip}\hypertarget{section-9}{}\section{9. Canonical models and moduli at volume one
thirtieth}\label{canonical-models-and-moduli-at-volume-one-thirtieth}}

\needspace{4\baselineskip}

\hypertarget{the-hilbert-series-and-the-fifth-map}{%
\hypertarget{section-9.1}{}\subsection{9.1. The Hilbert series and the fifth
map}\label{the-hilbert-series-and-the-fifth-map}}

We now assume \(p_g>0\) and \(V=1/30\) and determine the canonical
model. The bound for \(p_g\ge2\) \cite{ref8}
gives \(p_g=1\). The irregular bound in \hyperlink{section-3.1}{Section 3.1} gives \(q=0\), so
\(\chi=h^{2,0}\ge0\). \hyperlink{section-8}{Section 8} excludes the positive values.
Chen--Chen's equality classification \cite[Theorem 1.1(ii)]{ref7} therefore
gives \[
\chi=q=h^{2,0}=0,\quad (P_2,P_3,P_4,P_5)=(1,2,3,4),\qquad
B_Y=\{(1,2),(1,3),(1,5)\}.
\tag{9.1}\label{eq:9.1}
\] In particular \(30K_Y\) is Cartier by the terminal basket table
\cite[Remark 2.1]{refCH11}. Reid's formula \eqref{eq:2.1}, summing the periodic
corrections of periods \(2,3,5\), gives \[
H_R(T)=\sum_{m\ge0}P_mT^m
=\frac{1-T^{28}}{(1-T)(1-T^3)(1-T^4)(1-T^5)(1-T^{14})}.
\tag{9.2}\label{eq:9.2}
\] The periodic sums give this identity in every degree.

The fourth map has surface image: a pencil would have moving degree at
least two and canonical restriction would give \(V\ge1/18\). To exclude a surface image for the fifth map, we use the following
estimate for two ruling pencils.

\needspace{3\baselineskip}

\hypertarget{result-9.1}{}\ResultHeading{Lemma 9.1 (two rulings).} Suppose
\(Y\dashrightarrow\mathbb P^1\times\mathbb P^1\) has connected general
curve fibre, with ruling pencils \(F_{1,Y},F_{2,Y}\), and
\(mK_Y\sim F_{1,Y}+F_{2,Y}+E_Y\), \(E_Y\ge0\). Then \[
V\ge\frac8{m(m+1)^2}.
\tag{9.3}\label{eq:9.3}
\]

\ProofHeading{Proof.} We use the special model for the first pencil. The second
cuts a moving integral curve \(C\) on its smooth general surface \(F\),
possibly with fixed part. Effective adjunction gives \[
(m+1)D\sim_{\mathbb Q}K_F+C+B,\qquad B\ge0,
\] with \(B\) fixed as \(C\) varies. The curve is nef, \(AC\ge1\), and
its normalization has genus at least two. Thus \[
(m+1)AD\ge A^2+AC\ge2,\qquad
(m+1)DC\ge(K_F+C)C\ge2.
\] Intersecting again with \(D\) gives \((m+1)D^2\ge DA+DC\ge4/(m+1)\).
Applying this estimate to each ruling pencil and intersecting the
original relation with \(K_Y^2\), we obtain the required bound. \(\square\)

\needspace{3\baselineskip}

\hypertarget{result-9.2}{}\ResultHeading{Proposition 9.2.} The fifth-canonical map is generically
finite.

\ProofHeading{Proof.} Suppose otherwise. Multiplication by \(s\), together
with the fourth-map calculation, gives a nondegenerate surface \(\Sigma\subset\mathbb P^3\). Let
\(W\to B\to\Sigma\) be its Stein factorization, \(C\) a connected
general curve fibre, and \(S\in|M_5|\) general. Set
\(a=\deg(B/\Sigma)\deg\Sigma\ge2\), \(\xi=QC\). Thus \[
M_5|_S\equiv aC,\qquad a\xi=QM_5^2\le25V.
\] The numerical curve inequality \cite[Theorem 2.3]{ref7}, with
\(\beta=a/5\), gives \(\xi\ge2a/(6a+5)\). If \(a\ge3\), its choice
\(n=12\) has \((6-5/a)\xi>1\) and gives \(12\xi\ge4\). Hence
\(V\ge a\xi/25\ge1/25\), impossible. Thus \(a=2\), \(\Sigma\) is a
normal quadric, and \(B\to\Sigma\) is an isomorphism.

For a quadric cone, its ruling gives \(5K_Y\sim2F_Y+E_Y\) and \[
\frac4{49}\le K_Y^2F_Y\le\frac1{12}.
\] \cite[Lemma 3.12]{ref1} makes \(30K_Y^2F_Y\) integral, whereas the interval
\([120/49,5/2]\) contains no integer. For a smooth quadric, \hyperlink{result-9.1}{Lemma 9.1}
gives \(V\ge2/45\). In either case we obtain a contradiction to \(V=1/30\). \(\square\)

\needspace{4\baselineskip}

\hypertarget{the-canonical-ring-and-all-pluricanonical-maps}{%
\hypertarget{section-9.2}{}\subsection{9.2. The canonical ring and all pluricanonical
maps}\label{the-canonical-ring-and-all-pluricanonical-maps}}

\ProofHeading{Proof of Theorem 1.2.} We apply the recognition theorem of Chen--Hu--Jiang
\cite[Theorem 8.1]{ref9} to \((Z,H)=(Z,K_Z)\), with
\((a_1,a_2,a_3,a_4)=(3,4,5,14)\) and \(d=28\).
The numerical hypotheses are \(H^3=1/30\) and the Hilbert series
\eqref{eq:9.2}. The fourth map is not composed with a pencil, \hyperlink{result-9.2}{Proposition 9.2}
makes the fifth map generically finite, and \(Z\) is not rational
because it is of general type. Finally, the fourteenth map is birational
by \cite[Theorem 1.1]{ref1}. A well-formed quasismooth hypersurface with
these weights and degree is exhibited in \eqref{eq:9.8}. Thus all the hypotheses
of the recognition theorem hold.

We obtain homogeneous generators \(s,x,y,z,w\), of degrees
\(1,3,4,5,14\), for which
\[
A=\mathbb C[s,x,y,z],\qquad H_R(T)=(1+T^{14})H_A(T),
\tag{9.4}\label{eq:9.4}
\]
and an isomorphism
\[
R\simeq\mathbb C[s,x,y,z,w]/(w^2-F_{28}(s,x,y,z)).
\tag{9.5}\label{eq:9.5}
\]
Here \(A\) is a polynomial ring and \(S=\mathbb C[s,x,y,z,w]\) denotes
the ambient weighted polynomial ring. On \(s\ne0\), we put
\(a=x/s^3\), \(b=y/s^4\), \(c=z/s^5\), and \(e=w/s^{14}\).
The identification \eqref{eq:9.5} is therefore an application of the published
recognition theorem; the additional geometric step needed here was
\hyperlink{result-9.2}{Proposition 9.2}.

Conversely, a normal hypersurface \eqref{eq:9.5} with canonical singularities
has \(K_Z=\mathcal O_Z(1)\) and
\(K_Z^3=28/(1\cdot3\cdot4\cdot5\cdot14)=1/30\). Adjunction may be
performed on the weighted projective stack and descended reflexively;
the ambient singular locus has dimension at most one. The hypersurface
ring is Cohen--Macaulay of dimension four, so \[
H^1(Z,\mathcal O_Z(k))=H^2(Z,\mathcal O_Z(k))=0
\quad(k\in\mathbb Z).
\tag{9.6}\label{eq:9.6}
\] The same graded local-cohomology calculation identifies its ring with
the full canonical section ring; see \cite{int30,int24}. Thus \(p_g=1\),
\(q=h^{2,0}=0\), and every such hypersurface is a canonical model.
Quasismoothness has not been assumed.

The subring \(A\) is intrinsic, since it is generated in degrees at most
five, and \(R=A\oplus A(-14)w\). It defines a finite double cover \[
\rho:Z\longrightarrow\mathbb P(1,3,4,5)
\tag{9.7}\label{eq:9.7}
\] with canonical involution \(\iota(w)=-w\). In codimension one the
branch divisor is \((F_{28}=0)\); at quotient points the construction
uses reflexive sheaves. We can now read the canonical maps from
\(\mathbb C(Z)=\mathbb C(a,b,c,e)\), \(e^2=F_{28}(1,a,b,c)\):

\needspace{12\baselineskip}
\begin{longtable}[]{@{}
  >{\raggedright\arraybackslash}p{(\columnwidth - 4\tabcolsep) * \real{0.4200}}
  >{\raggedleft\arraybackslash}p{(\columnwidth - 4\tabcolsep) * \real{0.1700}}
  >{\raggedright\arraybackslash}p{(\columnwidth - 4\tabcolsep) * \real{0.4100}}@{}}
\toprule\noalign{}
\begin{minipage}[b]{\linewidth}\raggedright
System
\end{minipage} & \begin{minipage}[b]{\linewidth}\raggedleft
Image dimension
\end{minipage} & \begin{minipage}[b]{\linewidth}\raggedright
Image field or degree
\end{minipage} \\
\midrule\noalign{}
\endhead
\bottomrule\noalign{}
\endlastfoot
\(\varphi_1,\varphi_2\) & 0 & A point \\
\(\varphi_3\) & 1 & \(\mathbb C(a)\) \\
\(\varphi_4\) & 2 & \(\mathbb C(a,b)\) \\
\(\varphi_m\), \(5\le m\le13\) & 3 & Degree two \\
\(\varphi_m\), \(m\ge14\) & 3 & Birational \\
\end{longtable}

We also determine the general fibre of the fourth map. Over
\(k=\mathbb C(a,b)\), its equation has degree at most five in \(c\). If
the quadratic extension split over \(\bar k\), its connected geometric
fibre would be rational, contradicting general type by the easy addition
inequality. Its geometric generic fibre is therefore an integral double
cover of the \(c\)-line. General type forces genus at least two, while
the degree bound gives genus at most two. Consequently its squarefree part has degree five. Since the polynomial
itself has degree at most five, it has degree five and is squarefree
over \(\bar k\). Its smooth completion
is a genus-two curve, with pole orders two and five for \(c,e\) at
infinity. This proves \hyperlink{result-1.2}{Theorem 1.2}. The quasismooth example below will also show
that the volume bound is sharp and that degree fourteen is the optimal
uniform birationality bound. \(\square\)

\needspace{4\baselineskip}

\hypertarget{general-members-and-deformation}{%
\hypertarget{section-9.3}{}\subsection{9.3. General members and
deformation}\label{general-members-and-deformation}}

\ProofHeading{Proof of Theorem 1.3.} The polynomial \[
w^2+s^{28}+y^7+x^8y+xz^5
\tag{9.8}\label{eq:9.8}
\] defines a quasismooth member: vanishing of its derivatives gives
\(s=w=0\), \(xz^4=0\), \(8x^7y+z^5=0\), \(7y^6+x^8=0\), hence the
origin. Removing \(w\) from the same calculation shows that the branch surface
has a nonempty quasismooth locus.

A quasismooth member has exactly the singularities \[
\frac12(1,1,1),\qquad\frac13(1,2,2),\qquad\frac15(1,4,4).
\tag{9.9}\label{eq:9.9}
\] The last two occur at the weight-three and weight-five points. The
first is the unique point cut on the ambient line \(\mathbb P(4,14)\).
Quasismoothness requires nonzero \(y^7,w^2\) coefficients and excludes
the weight-four and weight-fourteen points. These chart calculations
also prove terminality and index thirty; weighted Lefschetz \cite[Theorem 3.2.4]{int24} gives rational Picard rank one.

For a general member the canonical divisor \(D=(s=0)\) is the
quasismooth surface \(D_{28}\subset\mathbb P(3,4,5,14)\) with
singularities \(\frac1r(1,1)\) for \(r=2,3,5\). Adjunction and
cohomology give \[
p_g(\widetilde D)=q(\widetilde D)=0,\qquad K_D^2=2/15,
\qquad K_{\widetilde D}^2=\frac2{15}-\frac13-\frac95=-2,
\tag{9.10}\label{eq:9.10}
\] where \(\widetilde D\) is the minimal resolution, not the minimal
model of its birational class. The last equality uses the discrepancy
\(-(r-2)/r\) of the exceptional \((-r)\)-curve. These surface assertions
concern a general member only.

We next deform an arbitrary canonical member. Given \(w^2=F\), we
choose a general weighted
polynomial \(G\) of degree twenty-eight and consider \[
\mathcal Z:\quad w^2=(1-\tau)F+\tau G.
\tag{9.11}\label{eq:9.11}
\] Its graded ring is free of rank two over \(\mathbb C[\tau,s,x,y,z]\).
The corresponding quotient stack is flat with hypersurface fibres,
which are Cohen--Macaulay and in particular satisfy \(S_2\). Relative adjunction gives
\(\omega=\mathcal O(1)\), and there are no codimension-one stabilizers.
Graded twists, or equivalently exact finite-group invariants on the
weighted charts, give flat relative reflexive powers commuting with base
change. This is an admissible \(\mathbb Q\)-Gorenstein family \cite[Theorem 5.3.6]{int1}. Canonicity is open \cite{int1,int37}; after shrinking about
\(0\) all fibres are canonical. The general member is quasismooth, since
\eqref{eq:9.8} makes that open locus nonempty.

For simple connectedness, we follow the specialization argument of
Chen--Hu--Jiang \cite[Section 9.3]{ref9}. A quasismooth member is
simply connected by \cite[Theorem 3.2.4]{int24}. The
total space of \eqref{eq:9.11} is normal after shrinking, since it is flat with
normal fibres over a smooth curve. Specialization for a normal total
space with reduced central fibre gives
\(\pi_1(\mathcal Z_\tau)\twoheadrightarrow\pi_1(Z)\) \cite[Lemma 2.8.1 and Theorem 2.12]{int38}. Thus \(\pi_1(Z)=1\). Takayama's resolution theorem
for klt varieties \cite[Theorem 1.1]{int44} gives the same fundamental group
on a resolution, and smooth projective birational invariance completes
\hyperlink{result-1.3}{Theorem 1.3}. \(\square\)

\needspace{4\baselineskip}

\hypertarget{the-parameter-quotient}{%
\hypertarget{section-9.4}{}\subsection{9.4. The parameter quotient}\label{the-parameter-quotient}}

\ProofHeading{Proof of Theorem 1.4.} We put \(G=\operatorname{Aut}_{\rm gr}S\), \(H=G/\mathbb G_m^{\rm wt}\),
and let \(U\subset\mathbb P(S_{28})\) be the open locus of normal
canonical hypersurfaces with nonzero \(w^2\) coefficient. It is nonempty
and irreducible. Every equality model occurs in \(U\) by \eqref{eq:9.5}.

An isomorphism of canonical models preserves their canonical rings. Each
generator degree \(1,3,4,5,14\) occurs once, and all are below the first
relation. Hence every isomorphism lifts to an invertible triangular
graded substitution of \(S\), taking one equation to a scalar multiple
of the other. The kernel on weighted projective space is
\(\mathbb G_m^{\rm wt}\) \cite[Lemma 1.3]{int26}. Consequently the fibres of
\(U\to\mathcal M_{1/30,1}\) are exactly the \(H\)-orbits.

Since the weights are distinct, \(G\) is connected and solvable: each variable of
weight \(a_i\) maps to a nonzero multiple of itself plus a polynomial of
weight \(a_i\) in lower-weight variables. Monomial counting gives \[
\dim S_{28}=137,\qquad
(\dim S_1,\dim S_3,\dim S_4,\dim S_5,\dim S_{14})=(1,2,3,4,24).
\tag{9.12}\label{eq:9.12}
\] Thus \(\dim G=34\) and \(\dim H=33\). Canonically polarized varieties
have finite automorphism groups, so every orbit has dimension
thirty-three and the coarse space has dimension \(136-33=103\).

Substitution is a linear representation on the fixed-degree space
\(S_{28}\). Lie--Kolchin triangularizes it, and Miyata's theorem \cite[Theorem 1]{int42} makes \(\mathbb C(S_{28})^G\) purely transcendental.
Weighted scalars act as all nonzero scalars on \(S_{28}\), via
\(\lambda\mapsto\lambda^{28}\); therefore \[
\mathbb C(\mathcal M_{1/30,1})
=\mathbb C(\mathbb P(S_{28}))^H
=\mathbb C(S_{28})^G.
\tag{9.13}\label{eq:9.13}
\] Thus the coarse moduli space is rational. We use the invariant function
field here, so the argument does not require \(H\) to be reductive.

\needspace{4\baselineskip}

\hypertarget{all-admissible-deformations-and-the-coarse-local-rings}{%
\hypertarget{section-9.5}{}\subsection{9.5. All admissible deformations and the coarse local
rings}\label{all-admissible-deformations-and-the-coarse-local-rings}}

We next show that this quotient description also controls infinitesimal
families, including deformations of nonquasismooth members. Let \(B\) be a local Artinian
\(\mathbb C\)-algebra and \(\mathcal Z/B\) an admissible deformation of
\(Z\). By \eqref{eq:9.6}, flatness of every reflexive canonical power and
induction over small extensions give finite free modules \[
R_{B,m}=H^0(\mathcal Z,\omega_{\mathcal Z/B}^{[m]})
\] of ranks \(P_m\), commuting with base change. We lift the five canonical
generators. Applying Nakayama's lemma in each degree gives a surjection \[
S_B=B[s,x,y,z,w]\longrightarrow R_B.
\] Its homogeneous kernel \(I\) reduces to \((w^2-F)\), because each
\(R_{B,m}\) is free. We lift the degree-twenty-eight relation to
\(f_B\in I_{28}\). In every degree the reduction of
\(S_{B,m-28}\xrightarrow{f_B}I_m\) is surjective; Nakayama's lemma therefore gives
\(I=(f_B)\). The \(w^2\) coefficient is a unit, so normalization and
completion of the square give \[
R_B\simeq B[s,x,y,z,w]/(w^2-F_B(s,x,y,z)).
\tag{9.14}\label{eq:9.14}
\] Relative Proj recovers \(\mathcal Z\), since its canonical divisor is
relatively ample. Conversely, arbitrary coefficient lifts of \(F_B\)
over a small extension of \(B\) define a flat family, and the
quotient-chart argument for \eqref{eq:9.11} verifies admissibility. Thus every
admissible deformation is induced by the equation and the deformation
functor is unobstructed.

The same argument lifts isomorphisms: the five generator images have unique
polynomial representatives below degree twenty-eight, and their
triangular leading coefficients are units. This proves the local family
presentation by \([U/H]\). Its tangent space at \(Z\) is the cokernel \[
\operatorname{Lie}(H)\longrightarrow
T_{[f]}\mathbb P(S_{28}).
\tag{9.15}\label{eq:9.15}
\] The map is injective, since a canonically polarized variety in
characteristic zero has no infinitesimal automorphisms. Its source and
target have dimensions \(33\) and \(136\). Hence the semiuniversal
admissible deformation space is smooth of dimension \(103\) at every
member. This argument retains the deformation directions of singular
index-one covers and uses no right-exact normal sequence for a singular
hypersurface.

The canonical-model moduli stack is therefore a smooth Deligne--Mumford
stack along this locus. The local quotient theorem \cite[Theorem 3.2]{int2}
gives \[
\widehat{\mathcal O}_{\mathcal M_{1/30,1},[Z]}
\simeq\mathbb C[[u_1,\ldots,u_{103}]]^{\operatorname{Aut}(Z)}.
\tag{9.16}\label{eq:9.16}
\] Equivalently, its analytic germ is a finite quotient of a smooth
\(103\)-dimensional germ. Such singularities are normal, rational, and
klt. The finite action may be ineffective. Together with the surjection
from the irreducible space \(U\) and \eqref{eq:9.13}, this proves the
corresponding assertions in \hyperlink{result-1.4}{Theorem 1.4}.

\needspace{4\baselineskip}

\hypertarget{automorphisms}{%
\hypertarget{section-9.6}{}\subsection{9.6. Automorphisms}\label{automorphisms}}

The intrinsic subring \(A\) and double cover \eqref{eq:9.7} give \[
1\longrightarrow\langle\iota\rangle\longrightarrow\operatorname{Aut}(Z)
\longrightarrow\operatorname{Aut}(\mathbb P(1,3,4,5),(F=0))
\longrightarrow1.
\tag{9.17}\label{eq:9.17}
\] An automorphism of the base preserving \((F=0)\) lifts because it
takes \(F\) to a scalar multiple, whose square root determines the
action on \(w\).

We apply Esser's theorem \cite[Theorem 4.1]{int26} to the branch surface. Its
degree satisfies \(28\ge5\max(1,3,4,5)\) and its quasismooth locus is
nonempty. For a very general branch surface its ambient automorphism
group is toric and central. To see that such an element is trivial, we normalize its
action on \(s\), and commute with the shears \(x\mapsto x+s^3\),
\(y\mapsto y+s^4\), \(z\mapsto z+s^5\). All diagonal scalars are then
one. Thus a very general \(Z\) has
\(\operatorname{Aut}(Z)=\langle\iota\rangle\simeq\mathbb Z/2\).

For an arbitrary member, \eqref{eq:9.14} extends \(w\mapsto-w\) over the identity of
every infinitesimal base. Hence the canonical involution acts trivially
on the semiuniversal deformation space. At a very general point the
quotient \eqref{eq:9.16} is therefore smooth; the smooth locus is open and
dense. This completes \hyperlink{result-1.4}{Theorem 1.4}. \(\square\)

The moduli space here parametrizes normal canonical models with Kollár's
admissible family condition. No assertion about the nonnormal or
strictly log-canonical boundary of its stable compactification is
needed.

\needspace{4\baselineskip}

\hypertarget{appendix-a.-an-exact-polynomial-calculation}{%
\section{Appendix A. An exact polynomial calculation}\label{appendix-a.-an-exact-polynomial-calculation}}

We include the exact calculation used in \hyperlink{section-4.4}{Section 4.4}. The following code constructs the linear operator on each weighted homogeneous space in the range established there. We record a monomial by its exponents in \(t,L,z\), allowing negative exponents of \(t\). Thus the assertions check all four parameters and all twenty-two weights using rational arithmetic.

\begin{verbatim}
from fractions import Fraction as Q

def rank(rows):
    a=[list(map(Q,row)) for row in rows]
    if not a:return 0
    i=0
    for j in range(len(a[0])):
        k=next((k for k in range(i,len(a)) if a[k][j]),None)
        if k is None:continue
        a[i],a[k]=a[k],a[i]
        v=a[i][j];a[i]=[x/v for x in a[i]]
        for k in range(i+1,len(a)):
            if a[k][j]:
                v=a[k][j];a[k]=[x-v*y for x,y in zip(a[k],a[i])]
        i+=1
        if i==len(a):break
    return i

count=0
for b in (Q(35,24),Q(5,8),Q(1,8),Q(-1,24)):
    for w in range(4,26):
        basis=[(j,k) for j in range(8) for k in range(8-j) if (w-3*(j+k))%2==0]
        target=[(j,k) for j in range(9) for k in range(9-j) if (w+3-3*(j+k))%2==0]
        idx={v:i for i,v in enumerate(target)}
        mat=[[Q(0) for _ in basis] for _ in target]
        for c,(j,k) in enumerate(basis):
            i=Q(w-3*(j+k),2)
            entries=[((j,k+1),i+Q(k,2)-2)]
            if j:entries.append(((j-1,k),Q(j)))
            if k:entries += [((j,k-1),k*(2+b)/3),((j+2,k-1),-k*b/2)]
            for term,v in entries:mat[idx[term]][c]+=v
        r=rank(mat); nullity=len(basis)-r
        assert nullity==(1 if w==4 else 0),(b,w,nullity)
        if w==4:
            c=basis.index((0,0));assert all(row[c]==0 for row in mat)
        count+=1
print(f'PASS: {count} exact rational matrices; only weight4 has kernel, generated by t^2.')
\end{verbatim}

\needspace{4\baselineskip}

\hypertarget{appendix-b.-exact-basket-computations}{%
\section{Appendix B. Exact basket computations}\label{appendix-b.-exact-basket-computations}}

\subsection{B.1. The open volume interval}\label{b-open-interval}

We now give the finite recursion used to obtain the numerical tables. It uses exact rational arithmetic and enumerates every basket satisfying the specified numerical conditions. In particular, its input is the numerical data and the proved denominator bound, rather than a database of known threefolds.

Write a terminal basket as a multiset of coprime pairs \((b,r)\) with \(0<b\le r/2\), and put
\[
\gamma(b,r)=\min\{r-2b,b\},\qquad
\delta(b,r)=b-\gamma(b,r),\qquad
\epsilon(b,r)=\frac b2-\frac{b^2}{r}-\frac{\gamma(b,r)}6.
\]
For prescribed \(\chi,P_2,P_3,P_4\), Riemann--Roch gives
\[
\begin{split}
S:=\sum b&=10\chi+5P_2-P_3,\\
D:=\sum\delta&=P_4+6P_2+5\chi-4P_3,\\
G:=\sum\gamma&=S-D,\\
\sum\epsilon&=A-V,\qquad
A=2P_2+6\chi-\frac S2-\frac G6.
\end{split}
\tag{B.1}\label{eq:B.1}
\]
The pairs \((1,2)\) and \((1,3)\) are the only pairs with \(\epsilon=0\), and they use respectively one unit of \(D\) and one unit of \(G\).

Here is an explicit index bound. If \(r<3b\), then
\[
r=2\delta+3\gamma\le2D+3G=2S+G.
\]
If \(r\ge3b\), then \(\gamma=b\le G\) and \(\epsilon=b/3-b^2/r\). For \(V\ge L\), provided \(b/3-A+L>0\), this gives
\[
r\le\frac{b^2}{b/3-A+L}.
\tag{B.2}\label{eq:B.2}
\]
These denominators are positive in every computation below. The displayed bounds produce a finite list of positive-gap pairs. Choose them in nondecreasing order, subtract their two budgets and their gaps, and fill the remaining budgets by \((1,3)\) and \((1,2)\). This describes every admissible multiset exactly once before identical outputs are merged.

For \((\chi,P_2,P_3,P_4)=(1,1,1,2)\), the budgets are
\[
(S,D,G,A)=(14,9,5,1/6).
\]
With \(1/48\le V<1/30\), \eqref{eq:B.2} and the first bound allow \(r\le33\). The recursion yields thirty-eight baskets. Seven have \(r_XV=1\) and violate the integral canonical-divisor test. The seven numerical witnesses are displayed in \hyperlink{section-6.1}{Section 6.1}. The remaining thirty-one rows, sorted first by volume and then lexicographically by their expanded baskets, are precisely the table in \hyperlink{section-6}{Section 6}. Every entry in that table is printed by the following standard-library Python program. The initial values \(P_0=P_1=1\) correspond to \(p_g=1\).

\begin{verbatim}
from fractions import Fraction as Q
from math import gcd, lcm

def plurigenera(B, chi=1, P2=1, bound=50):
    V = 2*P2 + 6*chi - sum(Q(b*(r-b), r) for b,r in B)
    P, ell = [1, 1], Q(0)
    for m in range(2, bound+1):
        ell += sum(Q(((m-1)*b % r)*(r-(m-1)*b % r), 2*r)
                   for b,r in B)
        val = Q(m*(m-1)*(2*m-1), 12)*V-(2*m-1)*chi+ell
        assert val.denominator == 1
        P.append(int(val))
    return V, P

def enumerate_baskets(chi, P2, P3, P4, lo, hi, closed=False):
    S = 10*chi + 5*P2 - P3
    D = P4 + 6*P2 + 5*chi - 4*P3
    G = S-D
    A = 2*P2 + 6*chi - Q(S,2) - Q(G,6)
    assert min(S,D,G) >= 0 and 0 < lo <= hi
    assert lo < hi or closed
    bounds = [2*S+G]
    for b in range(1, G+1):
        den = Q(b,3)-A+lo
        assert den > 0, "This index-bound argument is unavailable."
        bounds.append(int(Q(b*b)/den))
    rmax, atoms = max(bounds), []
    for b in range(1, S+1):
        for r in range(2*b, rmax+1):
            if gcd(b,r) != 1:
                continue
            g = min(r-2*b,b)
            d = b-g
            e = Q(b,2)-Q(b*b,r)-Q(g,6)
            if 0 < e <= A-lo and g <= G and d <= D:
                atoms.append((b,r,g,d,e))
    answers = {}
    def visit(start, g, d, used, chosen):
        V = A-used
        if lo <= V < hi or (closed and V == hi):
            B = tuple(sorted(chosen+((1,3),)*g+((1,2),)*d))
            v, P = plurigenera(B,chi,P2)
            assert v == V and P[2:5] == [P2,P3,P4]
            answers[B] = (V,P)
        for j in range(start,len(atoms)):
            b,r,ga,da,e = atoms[j]
            if ga <= g and da <= d and used+e <= A-lo:
                visit(j,g-ga,d-da,used+e,chosen+((b,r),))
    visit(0,G,D,Q(0),())
    return sorted(answers.items(),key=lambda row:(row[1][0],row[0]))

def domain_witness(B,V,P,chi=1):
    if lcm(*(r for b,r in B))*V != 1:
        return None
    h = [0]+[P[m]-P[m-1]+chi*int(m==2) for m in range(1,51)]
    for i in range(1,50):
        for j in range(1,51-i):
            if h[i]>0 and h[j]>0 and h[i+j] < h[i]+h[j]-1:
                return (i,j,h[i],h[j],h[i+j])
    return None

if __name__ == "__main__":
    raw = enumerate_baskets(1,1,1,2,Q(1,48),Q(1,30))
    deleted, retained = [], []
    for B,(V,P) in raw:
        w = domain_witness(B,V,P)
        (deleted if w else retained).append((B,V,P,w))
    assert (len(raw),len(deleted),len(retained)) == (38,7,31)
    for label,(B,V,P,w) in enumerate(retained,1):
        print(label,B,V,tuple(P[m] for m in (5,6,7,8,9,10,12,13)))
    print("Domain deletions:")
    for B,V,P,w in deleted:
        print(B,V,w)

    from collections import Counter
    for chi,lo in [(1,Q(1,40)),(1,Q(1,35))]:
        rows = enumerate_baskets(chi,1,1,2,lo,Q(1,30))
        kept = [P for B,(V,P) in rows if P[5]>=3 and P[6]==4
                and P[7]>=5 and all(P[m]>=P[m-1] for m in range(2,51))]
        print(chi,lo,len(rows),len(kept),
              Counter(tuple(P[5:8]) for P in kept))
        traced = [P for P in kept if P[10]-P[6] >= 6]
        assert len(traced) == (8 if lo == Q(1,40) else 4)
        print('Horikawa trace:',lo,len(kept),'->',len(traced))
\end{verbatim}

The same function verifies the Horikawa-fibre numerical reductions, without any extra index cutoff. Retain \(P_5\ge3\), \(P_6=4\), \(P_7\ge5\), and \(P_m\ge P_{m-1}\) for \(2\le m\le50\). The resulting counts are
\[
\begin{array}{c|c|r|r|c}
\chi&[L,1/30)&\text{raw}&\text{retained}&(P_5,P_6,P_7)\\ \hline
1&[1/40,1/30)&21&15&(3,4,5)\ \text{only}\\
1&[1/35,1/30)&9&7&(3,4,5)\ \text{only}.
\end{array}
\tag{B.3}\label{eq:B.3}
\]
The further trace inequality \(P_{10}-P_6\ge6\) leaves eight and four
baskets, respectively; the latter four are those in Proposition~4.1.
\hypertarget{b.2.-the-closed-endpoint}{%
\subsection{B.2. The closed endpoint}\label{b.2.-the-closed-endpoint}}

The same recursion, with \texttt{lo=hi=Q(1,30)} and
\texttt{closed=True}, gives the endpoint list in \hyperlink{section-8.3}{Section 8.3}. The bound
\eqref{eq:B.2} is valid without a strict upper inequality. For
\((h,P_2,P_3,P_4)=(1,1,1,1)\), \eqref{eq:B.1} gives \(A-V=-1/30\), so that case
is empty before enumeration.

The following continuation checks the endpoint counts of Section 8.3 and the numerical properties used in Section 8.4. Proposition~8.3 proves \(h=1\) before enumeration. Monotonicity needs only finitely many tests, since
\[
P_{m+1}-P_m=\frac{m^2}{60}-2+
\sum_{(b,r)}\frac{\overline{mb}(r-\overline{mb})}{2r}>0
\qquad(m\ge11).
\]
The domain test is only a necessary inequality. The full-character calculation enumerates all lower weights at most the target and retains the character modulo the full local index.

\begin{verbatim}
# Continuation of the preceding basket program: the closed endpoint.
def endpoint_rows():
    cases = [(1,1,2),(1,1,3),(1,2,2),(1,2,3)]
    expected = [(4,3,0,0),(11,9,3,3),(6,5,2,0),(24,12,7,7)]
    keep = []
    for (chi,p3,p4), wanted in zip(cases, expected):
        rows = enumerate_baskets(chi,1,p3,p4,Q(1,30),
                                 Q(1,30),closed=True)
        counts = [len(rows),0,0,0]
        for B,(V,P) in rows:
            if P[5]>4 or any(P[m+1]<P[m] for m in range(1,50)):
                continue
            counts[1] += 1
            if domain_witness(B,V,P,chi):
                continue
            counts[2] += 1
            if (p3,p4)==(2,2):
                marks = [(b,r,Q(a,r)) for b,r in sorted(set(B))
                         for a in range(1,r)
                         if Q(5,14)<=Q(a,r)<=Q(2,5)]
                assert marks == [(7,20,Q(2,5))]
                assert B.count((7,20)) == 1
                assert not any((u-v+7*w)%20 == 13
                    for u in range(3) for v in range(3-u)
                    for w in range(3-u-v))
                continue
            counts[3] += 1
            rY = lcm(*(r for b,r in B))
            assert rY > 30 and rY % 30 == 0
            keep.append((B,P))
        assert tuple(counts) == wanted
        print('Endpoint:',(chi,p3,p4),counts)
    assert len(keep) == 10
    for n,(basket,P) in enumerate(keep,1):
        print('A'+str(n),basket,tuple(P[3:8]))
    trace_excluded = [n for n,(B,P) in enumerate(keep,1)
                      if P[10]-P[6] < 7]
    assert trace_excluded == [3,4]
    print('Endpoint whole-member trace:',trace_excluded,
          'remaining:',len(keep)-len(trace_excluded))
    bad = lambda row: any(4%r and 15%r for b,r in row[0])
    assert all(bad(keep[n-1]) for n in (5,6,7,10))
    assert not any(r==3 for b,r in keep[8][0])
    assert keep[7][1][5]==4
    return keep

def initial_monomials(b,r,a,canonical=False):
    p,e = a.numerator,a.denominator
    c = (-b*pow(p,-1,e))%e
    inv = pow(c,-1,e)
    chars = [1,-1,b,0]
    weights = [(inv*q)%e or e for q in chars]
    target = p if canonical else 4*p-2*e
    character = -b if canonical else -4*b
    ans = []
    def visit(i,left,exponents,char):
        if i==4:
            if (char-character)%r==0:
                ans.append(exponents)
            return
        for k in range(left//weights[i]+1):
            visit(i+1,left-k*weights[i],exponents+[k],
                  char+k*chars[i])
    visit(0,target,[],0)
    return ans

if __name__ == '__main__':
    endpoint_rows()
    assert initial_monomials(5,12,Q(2,3),True)==[]
    for b,r,a in [(10,21,Q(13,21)),(10,21,Q(2,3))]:
        assert initial_monomials(b,r,a)==[[2,0,0,0]]
    print('PASS: endpoint arithmetic and full-character tests.')
\end{verbatim}

The initial counts are \(45\to29\to12\to10\). The ten candidates are precisely the rows of \hyperlink{section-8.4}{Section 8.4}. The trace test excludes A3 and A4; Proposition~8.4 excludes the other eight. The program certifies the arithmetic coverage and the numerical trace tests.

\needspace{4\baselineskip}

\end{document}